\documentclass[10pt, reqno]{amsart}

\usepackage{amssymb, mathrsfs, bbold}
\usepackage{amsmath, amsthm, yhmath}
\usepackage{accents}
\usepackage{enumerate}
\usepackage{dsfont}
\usepackage{color}
\usepackage[a4paper]{geometry}
\usepackage[utf8]{inputenc}   % To use accentuated letters \UTF{008E}, \UTF{008F}... directly in the Latex file
\usepackage{stmaryrd}
\usepackage{titlesec, wasysym}
\usepackage{appendix}
\usepackage{todonotes, fancyhdr}
\usepackage{chngcntr}
\usepackage{cancel, stmaryrd}
\usepackage{tikz-cd}
\usepackage{hyperref} %\RequirePackage[colorlinks,citecolor=purple,urlcolor=blue]{hyperref}

\usepackage{setspace}  \renewcommand{\baselinestretch}{0.99} % \renewcommand{\baselinestretch}{0.99}

\usepackage[all]{xy}
\numberwithin{subsection}{section}
\numberwithin{subsubsection}{subsection}
\numberwithin{equation}{section}   % Equation (p.q) is the q-st equation of Section p.

\renewcommand{\theenumi}{\alph{enumi}}
\renewcommand{\labelenumi}{\textsf{(\theenumi)}}

\newenvironment{Dem}[1][\unskip]{%
    \begin{list}{\hspace{1cm}{{\textit{ Proof #1 --}}}}{%
        \setlength{\topsep}{0pt}%
        \setlength{\leftmargin}{0pt}%
        \setlength{\rightmargin}{0pt}%
        \setlength{\listparindent}{0pt}%
        \setlength{\itemindent}{0pt}%
        \setlength{\parsep}{0pt}%
        \addtolength{\leftmargin}{0pt}  %\addtolength{\leftmargin}{20pt}% 
        \addtolength{\rightmargin}{0pt}%
    } \item }{\hfill $\rhd$\end{list}\smallskip}

\newenvironment{Dem*}[1][\unskip]{%
    \begin{list}{\hspace{0cm}{\sf \textbf{{\small Proof #1 --}}}}{%
        \setlength{\topsep}{0pt}%
        \setlength{\leftmargin}{0pt}%
        \setlength{\rightmargin}{0pt}%
        \setlength{\listparindent}{0pt}%
        \setlength{\itemindent}{0pt}%
        \setlength{\parsep}{0pt}%
        \addtolength{\leftmargin}{20pt}%
        \addtolength{\rightmargin}{0pt}%
    } \item }{\hfill $\rhd$\end{list}\smallskip}

\renewcommand\thesection       {\arabic{section}}
\renewcommand\thesubsection    {\thesection{\boldmath $.$}\arabic{subsection}}
\renewcommand\thesubsubsection    {\thesection{\boldmath $.$}\arabic{subsection}{\boldmath $.$}\arabic{subsubsection}}

\titleformat{\section}[block]
{\filcenter\normalfont\bfseries}  % {\filcenter\normalfont\sffamily\bfseries}
{{\hspace{-0cm}}\thesection \hspace{0.2em} --\vspace{0.1cm}}{0.5em}{}  %{ {\hspace{-0.7cm}}\thesection \hspace{0.2em} --\vspace{0.1cm}}{0.5em}{}

\titleformat{\subsection}[runin]
{\filcenter\normalfont\bfseries} % {\filcenter\normalfont\sffamily\bfseries}
{{\hspace{-0cm}}\thesubsection \hspace{0.2em} \hspace{-0.2cm}}{0.5em}{}

\titleformat{\subsubsection}[runin]
{\filcenter\normalfont\bfseries}
{{\hspace{-0cm}}{\thesubsubsection} \hspace{-0.1em} \hspace{-0.2cm}}{0.5em}{}  % {{\hspace{-1cm}}{\sffamily\thesubsubsection} \hspace{-0.1em} \hspace{-0.2cm}}{0.5em}{}

\newtheoremstyle{mystyle}
{3pt}               %space above
{3pt}               %space below
{\it }                      %bodyfont
{}                      %indent
{\rmfamily}      % {\sffamily\bfseries}       %headfont
{}                      %punctuation
{0.5em}                 %space after head
{\hspace{0cm}{\textbf{{#2}}\hspace{0cm} --} {\hspace{-0cm}}{\textbf{{#1.}}}}   %   {#1 #2{\hspace{0.2cm}--\hspace{-0.2cm}}  }   %{\llap{#2 }#1{\hspace{0.2cm}--}}

\theoremstyle{mystyle}

\newtheorem{thm}{Theorem}
\newtheorem*{thm*}{Theorem}

\newtheorem{cor}[thm]{\hspace{-0.15cm}  {Corollary}}% [section]
\newtheorem{lem}[thm]{\hspace{-0.14cm}  {Lemma}}%[section]
\newtheorem{prop}[thm]{\hspace{-0.13cm} {Proposition}}%[chapter]
\newtheorem{defn}[thm]{ \hspace{-0.32cm} {Definition}}%[chapter]

\newtheoremstyle{mystyle2}
{3pt}               %space above
{3pt}               %space below
{\it }                      %bodyfont
{}                      %indent
{\bfseries}             %headfont
{}                      %punctuation
{0.5em}                 %space after head
{\llap{#2 }\textbf{{#1}}{\hspace{0.2cm}--}}

\theoremstyle{mystyle2}

\newtheorem*{definition*}{Definition}
\newtheorem*{theorem*}{Theorem}
\newtheorem*{Remark*}{Remark}
\newtheorem*{rem*}{\hspace{-0.15cm} {Remark}}
\newtheorem*{lem*} {Lemma}
\newtheorem*{defn*} {Definition}
\newtheorem*{prop*} {Proposition}
\newtheorem*{cor*} {Corollary}

\newtheoremstyle{mystyle3}
{3pt}               %space above
{3pt}               %space below
{\it }                      %bodyfont
{}                      %indent
{\bfseries}             %headfont
{}                      %punctuation
{0.5em}                 %space after head
{\llap{#2 }\textbf{\textit{#1}}{\hspace{0.2cm}--}}

\theoremstyle{mystyle2}

\newtheorem*{lem**} {Lemma}

\newcommand{\ssk}{\smallskip}

\newcommand{\eps}{\varepsilon}

\newcommand\bbE{\mathbb{E}}

\newcommand{\bbH}{\mathbb{H}}

\newcommand{\bbK}{\mathbb{K}}
\newcommand\bbN{\mathbb{N}}
\newcommand\bbR{\mathbb{R}}
\newcommand\bbP{\mathbb{P}}
\newcommand\bbQ{\mathbb{Q}}
\newcommand{\bbT}{\mathbb{T}}

\newcommand\bbZ{\mathbb{Z}}
\newcommand{\bbOmega}{\mathbb{\Omega}}

\newcommand{{\mcB}}{\mathcal{B}}

\newcommand{\mcE}{\mathcal{E}}
\newcommand{\mcF}{\mathcal{F}}

\newcommand{\mcI}{\mathcal{I}}
\newcommand{\mcJ}{\mathcal{J}}
\newcommand{\mcK}{\mathcal{K}}
\newcommand{\mcL}{\mathcal{L}}
\newcommand{\mcM}{\mathcal{M}}
\newcommand{\mcN}{\mathcal{N}}

\newcommand{\mcQ}{\mathcal{Q}}
\newcommand\mcS{\mathcal{S}}

\newcommand{\mcR}{\mathcal{R}}

\newcommand{\scrR}{\ensuremath{\mathscr{R}}}

\newcommand{\strA}{\emph{\textrm{A}}}

\makeatletter
\newcommand*{\defeq}{\mathrel{\rlap{%
                     \raisebox{0.3ex}{$\m@th\cdot$}}%
                     \raisebox{-0.3ex}{$\m@th\cdot$}}%
                     =}
\makeatother

\makeatletter
\newcommand*{\eqdef}{=\mathrel{\rlap{\raisebox{0.3ex}{$\m@th\cdot$}}%
					                     \raisebox{-0.3ex}{$\m@th\cdot$}}%
                     }
\makeatother

\newcommand*{\mybullet}{\mathrel{\scalebox{0.5}{$\bullet$}}}

\makeatletter
\newcommand{\xdashrightarrow}[2][]{\ext@arrow 0359\rightarrowfill@@{#1}{#2}}
\def\rightarrowfill@@{\arrowfill@@\relax\relbar\rightarrow}
\def\arrowfill@@#1#2#3#4{%
  $\m@th\thickmuskip0mu\medmuskip\thickmuskip\thinmuskip\thickmuskip
   \relax#4#1
   \xleaders\hbox{$#4#2$}\hfill
   #3$%
}
\makeatother

\usepackage{tikz-cd}
\usetikzlibrary{calc}
\usetikzlibrary{shapes.misc}
\usetikzlibrary{shapes.symbols}
\usetikzlibrary{snakes}
\usetikzlibrary{decorations}
\usetikzlibrary{decorations.markings}

\makeatletter
\pgfdeclareshape{crosscircle}
{
  \inheritsavedanchors[from=circle] % this is nearly a circle
  \inheritanchorborder[from=circle]
  \inheritanchor[from=circle]{north}
  \inheritanchor[from=circle]{north west}
  \inheritanchor[from=circle]{north east}
  \inheritanchor[from=circle]{center}
  \inheritanchor[from=circle]{west}
  \inheritanchor[from=circle]{east}
  \inheritanchor[from=circle]{mid}
  \inheritanchor[from=circle]{mid west}
  \inheritanchor[from=circle]{mid east}
  \inheritanchor[from=circle]{base}
  \inheritanchor[from=circle]{base west}
  \inheritanchor[from=circle]{base east}
  \inheritanchor[from=circle]{south}
  \inheritanchor[from=circle]{south west}
  \inheritanchor[from=circle]{south east}
  \inheritbackgroundpath[from=circle]
  \foregroundpath{
    \centerpoint%
    \pgf@xc=\pgf@x%
    \pgf@yc=\pgf@y%
    \pgfutil@tempdima=\radius%
    \pgfmathsetlength{\pgf@xb}{\pgfkeysvalueof{/pgf/outer xsep}}%  
    \pgfmathsetlength{\pgf@yb}{\pgfkeysvalueof{/pgf/outer ysep}}%  
    \ifdim\pgf@xb<\pgf@yb%
      \advance\pgfutil@tempdima by-\pgf@yb%
    \else%
      \advance\pgfutil@tempdima by-\pgf@xb%
    \fi%
    \pgfpathmoveto{\pgfpointadd{\pgfqpoint{\pgf@xc}{\pgf@yc}}{\pgfqpoint{-0.707107\pgfutil@tempdima}{-0.707107\pgfutil@tempdima}}}
    \pgfpathlineto{\pgfpointadd{\pgfqpoint{\pgf@xc}{\pgf@yc}}{\pgfqpoint{0.707107\pgfutil@tempdima}{0.707107\pgfutil@tempdima}}}
    \pgfpathmoveto{\pgfpointadd{\pgfqpoint{\pgf@xc}{\pgf@yc}}{\pgfqpoint{-0.707107\pgfutil@tempdima}{0.707107\pgfutil@tempdima}}}
    \pgfpathlineto{\pgfpointadd{\pgfqpoint{\pgf@xc}{\pgf@yc}}{\pgfqpoint{0.707107\pgfutil@tempdima}{-0.707107\pgfutil@tempdima}}}
  }
}
\makeatother

\colorlet{symbols}{black}    
\colorlet{testcolor}{green!60!black}
\colorlet{supcolor}{red!60!black}

\tikzset{
	root/.style={circle, fill=testcolor!70, draw=testcolor, inner sep=1pt, minimum size=0.5mm},
	dot/.style={circle, draw=black, fill=black, inner sep=0pt, minimum size=0.2mm},
	noise/.style={circle, draw=black, fill=white, inner sep=0pt, minimum size=1.4mm},
	noiseblue/.style={circle, fill=blue!40, draw=blue, inner sep=0pt, minimum size=1mm},
	noisegray/.style={circle, fill=gray!40, draw=gray, inner sep=0pt, minimum size=1mm},
	blackdot/.style={circle, draw=black, fill=black, inner sep=0pt, minimum size=1.2mm},
	K/.style= {semithick, shorten >=0pt,shorten <=0pt,-},%{semithick,densely dashed,shorten >=1pt,shorten <=1pt,->},
	DK/.style={thick, densely dotted, shorten >=0pt,shorten <=0pt},   % 	DK/.style={semithick, densely dashed, shorten >=1pt,shorten <=1pt},
	}

\newcommand{\fr}{\mathfrak}
\newcommand{\id}{\mathop{\text{\rm id}}}
\newcommand{\spa}{\mathop{\text{\rm span}}}

\newcommand{\scT}{\mathscr{T}}
\newcommand{\scV}{\mathscr{V}}
\newcommand{\scW}{\mathscr{W}}

\newcommand{\mfs}{\mathfrak{s}}

\newcommand{\lp}{\llparenthesis\,}
\newcommand{\rp}{\,\rrparenthesis}
\newcommand{\tri}{|\!|\!|}
\newcommand{\lb}{\llbracket}
\newcommand{\rb}{\rrbracket}

\newcommand{\dt}[1]{\accentset{\text{\scriptsize$\mybullet$}}{#1}}

\newcommand{\res}{\hspace{-0.03cm}:\hspace{-0.03cm}}

\allowdisplaybreaks

\newtheoremstyle{mystyle4}
{3pt}               %space above
{3pt}               %space below
{\it }                      %bodyfont
{}                      %indent
{\bfseries}             %headfont
{}                      %punctuation
{0.5em}                 %space after head
{\llap{#2 }{\hspace{-0cm}}{\textbf{{#1.}}}}

\theoremstyle{mystyle4}

\newtheorem*{lem11}{\ref{lem:Biposi}' -- Lemma}
\newtheorem*{lem12}{\ref{lem:Binega}' -- Lemma}
\newtheorem*{lem2}{\ref{lem:Pitog}' -- Lemma}
\newtheorem*{lem3}{\ref{lem:goal}' -- Lemma}

\newcommand{\rtdb}[2]
{\begin{tikzpicture}[baseline=#1pt]
\coordinate (A1) at (0,0);
#2
\end{tikzpicture}
}

\newcommand{\rtds}[2]
{\begin{tikzpicture}[scale=#1]
\coordinate (A1) at (0,0);
#2
\end{tikzpicture}
}

\newcommand{\rtdbs}[3]
{\begin{tikzpicture}[baseline=#1pt, scale=#2]
\coordinate (A1) at (0,0);
#3
\end{tikzpicture}
}

\newcommand{\pol}[3]
{\coordinate (A#1) at ($(A#2)+0.4*({cos(#3)},{sin(#3)})$);}

\newcommand{\drc}[1]
{\foreach \n in {#1} \filldraw[white] (A\n) circle (2pt);
\foreach \n in {#1} \draw (A\n) circle (2pt);}

\newcommand{\drb}[1]
{\foreach \n in {#1} \fill (A\n) circle (2pt);}

\newcommand{\drch}[1]
{\foreach \n in {#1} \filldraw[white] (A\n) circle (2pt);
\foreach \n in {#1} \draw (A\n) circle (2pt);
\foreach \n in {#1} \filldraw (A\n) circle (0.5pt);}

\newcommand{\drl}[2]
{\foreach \n in {#2} \draw[thick] (A#1)--(A\n);}

\newcommand{\drdt}[2]
{\foreach \n in {#2} \draw[thick, densely dotted] (A#1)--(A\n);}
\newcommand{\drdh}[2]
{\foreach \n in {#2} \draw[double distance=0.7pt, thick, densely dotted] (A#1)--(A\n);}

\begin{document}

\begin{center}
{\LARGE{\textbf{Random models on regularity-integrability structures   \vspace{0.5cm}}}}
\end{center}

\begin{center}
{\sc I. Bailleul \& M. Hoshino}
\end{center}

\vspace{0.5cm}

\begin{center}
\begin{minipage}{0.7\textwidth}
\begin{center}\textbf{Abstract}\end{center} \vspace{0.1cm} \renewcommand\baselinestretch{0.9} \footnotesize
We prove a convergence result for a large class of random models that encompasses the case of the BPHZ models used in the study of singular stochastic PDEs. We introduce for that purpose a useful variation on the notion of regularity structure called a regularity-integrability structure. It allows to deal in a single elementary setting with models on a usual regularity structure and their first order Malliavin derivative.
\end{minipage}
\end{center}

\vspace{0.6cm}

{\it 
\begin{center}
\begin{minipage}[t]{13cm}
\baselineskip =0.35cm
{\scriptsize 

\center{\textrm{Contents}}

\vspace{0.1cm}

	\textrm{1.~Introduction\dotfill  \pageref{SectionIntroduction}}

	\textrm{2.~Functional setting and notations\dotfill \pageref{SectionFunctionalSetting}}

	\textrm{3.~A convergence result\dotfill \pageref{SectionResult}}

	\textrm{\hspace{0.35cm}3.1~Decorated trees\dotfill \pageref{SubsectionDecoratedTrees}}

	\textrm{\hspace{0.35cm}3.2~Differentiable sectors\dotfill \pageref{SubsectionDifferentialSectors}}

	\textrm{\hspace{0.35cm}3.3~Regularity-integrability structures and their models\dotfill \pageref{SubsectionRIS}}

	\textrm{\hspace{0.35cm}3.4~Renormalized models\dotfill \pageref{SubsectionRenormalizedModels}}

	\textrm{\hspace{0.35cm}3.5~Spectral gap inequality\dotfill \pageref{SubsectionSpectralGap}}
	
	\textrm{\hspace{0.35cm}3.6~A convergence result\dotfill \pageref{SubsectionMainResult}}

\textrm{4.~The mechanics of convergence\dotfill \pageref{SectionMechanics}}

	\textrm{\hspace{0.35cm}4.1~Structure of the proof\dotfill \pageref{SubsectionStructureProof}}
	
	\textrm{\hspace{0.35cm}4.2~The initial case \emph{$\{\textsf{bd}(\scW,1,p)\}_{2\le p\le\infty}$}\dotfill \pageref{sec:ind0}}

	\textrm{\hspace{0.35cm}4.3~Step 1: From \emph{$\textsf{bd}(\scW,i,\infty)$} to \emph{$\textsf{bd}(\scV,i)$}\dotfill \pageref{sec:ind1}}

	\textrm{\hspace{0.35cm}4.4~Step 2: From \emph{$\textsf{bd}(\scW,i,p)$} and \emph{$\textsf{bd}(\scV,i)$} to the $\sf g$-part of \emph{$\textsf{bd}(\scW,i+1,p)$}\dotfill \pageref{sec:ind2}}
	
	\textrm{\hspace{0.35cm}4.5~Step 3: From \emph{$\big\{\textsf{bd}(\scW,i,p)\big\}_{2\le p\le\infty}$} and \emph{$\textsf{bd}(\scV,i)$} to the $\sf \Pi$-part of \emph{$\big\{\textsf{bd}(\scW,i+1,p)\big\}_{2\le p\le\infty}$}  \dotfill \pageref{sec:ind3}}

	\textrm{\hspace{0.35cm}4.6~From uniform boundedness to convergence results\dotfill \pageref{sec:ind4}}

\textrm{5.~Proofs of the lemmas\dotfill  \pageref{SectionLemmas}}

	\textrm{\hspace{0.35cm}5.1~Trees conforming to rules\dotfill \pageref{ProofStronglyConform}}
	
	\textrm{\hspace{0.35cm}5.2~Proofs of algebraic identities\dotfill \pageref{ProofLemmaComparisonFormula}}
	
	\textrm{\hspace{0.35cm}5.3~Proof of the Lemma \ref{*lem:BesovTaylor}\dotfill \pageref{SectionBesovTaylor}}
	
	\textrm{\hspace{0.35cm}5.4~Proof of the lemmas of Section \ref{sec:ind1}\dotfill \pageref{SectionLastProofsLemmas}}
	
	\textrm{\hspace{0.35cm}5.5~Proofs of the Lemmas of Section \ref{sec:ind2}\dotfill \pageref{SectionProofLemmaSection33}}
	
	\textrm{\hspace{0.35cm}5.6~Proofs of the lemmas of Section \ref{sec:ind3}\dotfill \pageref{SectionProofLemmas34}}

\textrm{6.~Proofs of Theorem \ref{ThmContinuity} and Theorem \ref{ThmLipschitzContinuity}\dotfill  \pageref{SectionThmContinuity}}

	\textrm{\hspace{0.35cm}6.1~Proof of Theorem \ref{ThmContinuity}\dotfill \pageref{SectionThmContinuity1}}
	
	\textrm{\hspace{0.35cm}6.2~Theorem \ref{ThmLipschitzContinuity} and random Fourier series\dotfill \pageref{SectionThmContinuity2}}

\textrm{A.~Reconstruction and multilevel Schauder estimates in regularity-integrability structures\dotfill  \pageref{SectionAppendix}}

}\end{minipage}
\end{center}
}   \vspace{1cm}

\vspace{0.3cm}

%----------------------------%
\section{Introduction}
\label{SectionIntroduction}
%----------------------------%

{\it \S1. The four pillars.} The introduction by M. Hairer of the theory of regularity structures opened a new era in the domain of stochastic partial differential equations (PDEs). It provided in particular a robust solution theory for a number of equations, called `singular', whose study is beyond the range of the methods based on stochastic calculus. The singular feature of these equations is related to the fact that their formulations involve some ill-defined products. The development of this theory was done in several steps. The analytic core was developed in Hairer' seminal work \cite{Hai14}. Its algebraic backbone was deepened in Bruned, Hairer \& Zambotti's work \cite{BHZ}. The specific task of dealing with the ill-defined products of a singular stochastic PDE is called the renormalisation problem. This problem has a dynamics side and a probabilistic side. The dynamic meaning of the BPHZ renormalisation procedure of \cite{BHZ} was studied by Bruned, Chandra, Chevyrev \& Hairer in \cite{BCCH}, and lead to the identification of a solution to a singular stochastic PDE as the limit of solutions to renormalised versions of the initial equations, that is equations driven by a regularized noise with additional counterterms that typically diverge as the regularization parameter vanishes. The analytic machinery needs as an input an equation-dependent finite family of quantities built from a regularized noise. A systematic proof of probabilistic convergence of these quantities using the BPHZ renormalisation rule was given by Chandra \& Hairer in \cite{CH}; this is the probabilistic side of the renormalisation problem.

Altogether the four works \cite{Hai14, BHZ, BCCH, CH} form an automated blackbox for the study of a well identified large class of equations, with prominent examples coming as scaling limits of some microscopic discrete systems of statistical mechanics. This is the case of the (KPZ) equation from continuous interface growth models, of the parabolic Anderson model equation giving the scale limit of branching particle systems, or of the $\Phi^4_3$ equation from Euclidean quantum field theory. While the works \cite{Hai14, BHZ, BCCH} are now well understood by a growing community this is not the case of the work \cite{CH}. The latter uses ideas from the multiscale expansion method developed by Feldman, Magnen, Rivasseau \& S\'en\'eor in \cite{FMRS}, for the study of divergent Feynman integrals, to analyse the convergence problem of an equation-dependent finite collection of iterated integrals. The sophistication of their analysis and the very general assumptions on the law of the noise adopted in \cite{CH} make their work very challenging.

\bigskip

{\it \S2 The inductive mechanics of Linares, Otto, Tempelmayr \& Tsatsoulis.} Meanwhile, Otto developed with a number of co-authors a variant of the theory of regularity structures tailor made for the study of a certain class of singular quasilinear stochastic PDEs. Its analytic machinery was constructed in the works \cite{OW, OSSW} with Weber, and Sauer \& Smith. The algebraic machinery was described in the work \cite{LOT} with Linares \& Tempelmayr. Importantly, they were able to identify in \cite{OSSW} a renormalisation procedure with a similar dynamic meaning as the BPHZ renormalisation process. Linares, Otto, Tempelmayr \& Tsatsoulis proved in \cite{LOTT} the convergence result corresponding in their setting to the convergence result of \cite{CH}. Most interestingly, the authors of \cite{LOTT} used a set of assumptions and tools different from \cite{CH}, trading assumptions on cumulants and questions on iterated integrals for a spectral gap assumption on the law of the random noise and an iterative control of the stochastic objects. (We will introduce formally the spectral gap inequality in Section \ref{SubsectionSpectralGap}.) Their approach bypasses in particular the intricate algebraic content of the BPHZ strategy. We note that the idea of differentiation with respect to the noise that is involved in the spectral gap assumption was used in a different form in the early 80s by Caswell \& Kennedy \cite{CaswellKennedy} in their approach to perturbative renormalisation of quantum field theories. The results \cite{OSSW, LOT, LOTT} are not directly applicable to the study of semilinear subcritical singular stochastic PDEs. Hairer \& Steele \cite{HS23} gave very recently an improved and simplified version of the convergence result of \cite{LOTT}, in the original regularity structure setting. Their general convergence result for the BPHZ renormalisation procedure provides an alternative to the result of \cite{CH} of similar scope for practical purposes. We provide in the present work an alternative proof of their convergence result that holds for a larger class of renormalisation procedures containing the BPHZ procedure of \cite{HS23} as a particular example.

\medskip

Like \cite{LOTT}, the convergence proof of \cite{HS23} is done by induction. The objects to control are renormalized models on a regularity structure. They are built from a regularized noise and come under the form of a family of distributions ${\sf\Pi}^n_x\tau$ indexed by the points $x$ of the state space and a finite, equation-dependent, family of symbols $\tau$. The integer $n$ accounts here for the regularization parameter. The stochastic convergence of these models as the regularization is removed is mainly controlled by the $L^p(\Omega)$ convergence of some real-valued quantities of the form
\begin{equation} \label{EqTermToControlHS23}
\lambda^{-r(\tau)}({\sf \Pi}^n_x\tau)(\varphi^\lambda_x) 
\end{equation}
where the  smooth test functions $\varphi^\lambda_x$ behave like Dirac masses at $x$ as $\lambda$ goes to $0$ and $r(\tau)$ is some real number. The spectral gap assumption on the law of the random noise allows to control the $L^p(\Omega)$ norm of \eqref{EqTermToControlHS23} by its expectation and the quantity
\begin{equation} \label{EqExpectationDerivative}
\bbE\Big[\sup_{\Vert h\Vert_H\leq 1} \big\vert \big(d_\omega({\sf \Pi}^n_x\tau)(h)\big) (\varphi^\lambda_x) \big\vert^p\Big],
\end{equation}
where $d_\omega(\cdot)(h)$ stands for the G\^ateaux derivative in the direction $h$, for $h\in H$ in some space $H$. It turns out that a good control on the expectation of \eqref{EqTermToControlHS23} for $\lambda=1$ can be propagated by induction to any $0<\lambda\leq 1$ and all symbols. Building on the insight of \cite{LOTT} Hairer \& Steele show that $d_\omega({\sf \Pi}^n_x\tau)(h)$ can be represented as the reconstruction of a modelled distribution defined on an extended regularity structure that contains an extra noise symbol -- a placeholder for a generic $h$. This representation comes with estimates that play a crucial role in the inductive procedure. To apply this strategy Hairer \& Steele introduced a concept of {\it pointed modelled distribution} that allows to harvest the benefits associated with the improved regularity of the functions $h$ involved in the spectral gap assumption, compared to the regularity of the noise, and get as a consequence a good scaling bound for \eqref{EqExpectationDerivative}. This notion is closely related to the analytical mechanism developed in \cite{LOTT}. In both approaches, the Malliavin derivative of a model component is described as a modelled distribution relative to the original model, and its improved behaviour near the point at which the model is evaluated is encoded through localized two-scale estimates. Hairer \& Steele formulate this mechanism as an abstract theory of pointed Besov modelled distributions and prove a corresponding pointed reconstruction theorem. One then needs to extend the analytic core of the theory of regularity structures to the setting of pointed modelled distributions; a non-trivial task. Further, the construction of a pointed modelled distribution associated with the derivative $d_\omega({\sf \Pi}^n_x\tau)(h)$ of the renormalized model is only done in \cite{HS23} for BPHZ-like renormalisation procedures. 

\ssk

We use a different strategy to prove the convergence of a larger family of renormalized models. These models are built from a class of maps that act on the linear space spanned by the symbols $\tau$, called {\it preparation maps}. They were introduced by Bruned in \cite{Bru18} as a fundamental brick in the inductive construction of a large class of admissible models. BPHZ-like renormalisation procedures correspond to particular examples of preparation maps. The dynamic meaning of the renormalisation procedure associated with (strong) preparation maps was studied by Bailleul \& Bruned in \cite{BB21}. As in the BPHZ setting it involves some renormalized equations that include some additional counterterms. Since those counterterms depend on the regularization parameter $n$, we have to consider a family of preparation maps indexed by $n$. A family $(R_n)_{n\geq 0}$ of preparation maps is called a {\it renormalization procedure}.

\bigskip

{\it \S3. Our main results.} We trade in our setting the testing operation \eqref{EqTermToControlHS23} against a scaled centered function for a testing operation against some kernel $Q_t(x,\cdot)$ and we aim at getting some probabilistic bounds on quantities of the form $t^{-r(\tau)/\ell} \mcQ_t(x,{\sf \Pi}_x^n(\tau))$, where ${\sf \Pi}_x^n(\tau)$ is associated with a regularized noise $\xi_n$ and an $n$-dependent preparation map. At the informal and simplified level of this introduction, our main result, Theorem \ref{thm:main}, reads as follows. 

\ssk

\begin{thm} \label{ThmMainIntro}
Let $\xi$ be an $\mcS'(\bbR^d)$-valued random distribution whose law is centered, translation-invariant and satisfies an $H^{-s_0}$-spectral gap inequality, for some real exponent $s_0$. Then $\xi$ has any H\"older regularity strictly less than $-d/2-s_0$. Assume as well that we are given a regularity structure whose only symbol with degree less than or equal to $-d/2$ is the noise symbol. Then for any approximation of the identity $(\varrho_n)_{n\geq 0}$, for any renormalization procedure $(R_n)_{n\geq 0}$ such that the quantities $\bbE\big[\mcQ_1(0,{\sf \Pi}_0^n\tau)\big]$ converge for all the trees in the regularity structure with non-positive degree, the renormalized models associated with the regularized noises $\xi_n=\varrho_n*\xi$ and the preparation maps $R_n$ converge in $L^q(\bbP)$ for any $1\leq q<\infty$. The limit model is independent of the approximation of the identity and depends on the renormalization procedure only via the finite collection of numbers $\big\{\lim_{n\to\infty} \bbE\big[\mcQ_1(0, {\sf\Pi}_0^n\tau)\big] \big\}$, where $\tau$ runs over the set of symbols of the regularity structure with non-negative degree.
\end{thm}

\ssk

A non-Euclidean setting is required when studying parabolic equations. The dimension $d$ in the above statement is then replaced by the size $|\mfs|$ of the relevant scaling vector $\mfs$ -- see \eqref{*eq:interpolateHOmega}. Given a parabolic type subcritical stochastic PDE of the form \eqref{*eq:cor:generalSPDE} below, Hairer \cite{Hai14} provided a continuous map $\mcS$, with arguments an initial value $u_0$ of the equation and a model $\sf M$, that gives a modelled distribution solution of some regularity structure lift of Equation \eqref{*eq:cor:generalSPDE}. Therefore, for a converging sequence of renormalized models $({\sf M}^n)_{n\geq 0}$, we also have the convergence of the ${\sf M}^n$-reconstructions $(u_n)_{n\geq 0}$ of the modelled distributions $(\mcS(u_0,{\sf M}^n))_{n\geq 0}$. Since Bailleul \& Bruned \cite{BB21} derived an explicit equation for $u_n$ in terms of the preparation maps $R_n$, the combination of the \cite{BB21} and the present work give back the following known result. Below $\nabla$ and $\Delta$ denote the spatial gradient and Laplacian respectively. Note that the role of $d$ in Theorem \ref{ThmMainIntro} is played in a parabolic setting by the size $d+2$ of the scaling vector $(2,1,\dots,1)$.

\ssk

\begin{cor} \label{CorMainThm}
Let $\xi$ be a random noise as in Theorem \ref{ThmMainIntro}.
Consider a parabolic equation of the form
\begin{equation}\label{*eq:cor:generalSPDE}
(\partial_t-\Delta)u=F(u,\nabla u,\xi),\qquad t>0,\ x\in\bbT^d,
\end{equation}
where the unknown $u$ is a parabolic function/distribution and $F$ is a smooth function which is affine with respect to $\xi$. Suppose that we can construct the regularity structure associated with the above equation and all the assumptions of Theorem \ref{ThmMainIntro} are satisfied for the renormalized models associated with some smooth noises $\xi_n$ and a renormalization procedure $(R_n)_{n\geq 0}$. Then there exists a sequence of smooth functions $C_n$ such that the local in time solutions $u_n$ of the equations
$$
(\partial_t-\Delta)u_n = F(u_n,\nabla u_n,\xi_n) + C_n(u_n,\nabla u_n),
$$
with common, regular enough, initial function/distribution $u_n|_{t=0}$, converge in probability as $n$ goes to $\infty$.
\end{cor}

\ssk

The counterterm $C_n$ is explicitly known from $F$ and the preparation map $R_n$ -- see Theorem 4.1 of \cite{BB21}. The technical details on the space of initial values and the topology of convergence are the same as in \cite{Hai14}. They are typically (spatial or spacetime) H\"older--Besov norms $C^\alpha$ with some $\alpha\in\bbR$ determined by $F$, $d$, and the regularity of the noise $\xi$. 

\medskip

The class of renormalized models from \cite{BHZ} is built from a sub-class of preparation maps. Within that sub-class there is a unique choice of preparation maps such that $\bbE[({\sf \Pi}^n\tau)(x)]=0$ for all $\tau$ with negative degree, $x\in\bbR^d$ and $n\geq 1$ -- this is the content of Theorem 6.18 of \cite{BHZ}. We will see that Theorem \ref{ThmMainIntro} implies that the limit of the renormalized models $({\sf M}^n)_{n\geq 0}$ exists and depends only on the noise $\xi$ but does not depend on the approximation procedure. The limit model ${\sf M}(\xi)$ is called the {\it BPHZ model}. If $\xi$ has a law $\bbP$ we also denote by ${\sf M}(\bbP)$ the law of ${\sf M}(\xi)$ in the space of models. The next statement expresses a continuity property of the law of ${\sf M}(\bbP)$ in the class of probability measures that have the same spectral gap. Tempelmayr \cite{Tempelmayr} obtained recently a similar result in a different setting. The probability space $\Omega$ is here the sample space of the random noise, typically a Besov space of negative regularity exponent with a polynomial weight. 

\ssk

\begin{thm} \label{ThmContinuity}
Let $(\bbP_j)_{j\in\bbN}$ be a sequence of centered and translation-invariant probability measures on $\Omega$ that converges weakly to a limit probability measure $\bbP$. If all the $\bbP_j$ satisfy a spectral gap inequality with the same $H^{-s_0}$ space and the same constant then the law of ${\sf M}(\bbP_j)$ converges weakly to the law of ${\sf M}(\bbP)$.
\end{thm}

\ssk

This qualitative result ensures for instance the independence of the BPHZ model with respect to the regularization procedure used to define it within a large class of regularization procedures. Theorem \ref{ThmContinuity} has a quantitative version which involves a notion of distance $\|{\sf M}_1 \hspace{-0.03cm}:\hspace{-0.03cm} {\sf M}_2\|_{\textbf{\textsf{M}}(\scW)_{w_c}}$ on the space $\textbf{\textsf{M}}(\scW)_{w_c}$ of models we will work with. This space is properly introduced in Section \ref{SubsectionRIS} and we only give here an informal statement which will be fully detailed in Section \ref{SectionThmContinuity}.

\ssk

\begin{thm} \label{ThmLipschitzContinuity}
Let $(\xi_j)_{j\in\bbN}$ be a sequence of random noises whose laws are centered, translation-invariant and satisfy a spectral gap inequality with the same $H^{-s_0}$ space and the same constant. Assume the $\xi_j$ converge to a limit random distribution $\xi$ in $L^r(\Omega)$ for any $1\leq r<\infty$. Then the BPHZ models $\emph{\textsf{M}}(\xi_j)$ converge to the BPHZ model $\emph{\textsf{M}}(\xi)$ in $L^q\big(\Omega,\textbf{\textsf{M}}(\scW)_{w_c}\big)$ for any $1\leq q<\infty$ at speed $\Vert \xi_j-\xi\Vert_{L^r(\Omega)}$ with some $q< r<\infty$.
\end{thm}

\ssk

We illustrate this statement in Section \ref{SectionThmContinuity} with the example of some random Fourier series.

\bigskip

{\it \S4 What is new here?} Of course, the convergence result of the BPHZ models was first proved by Chandra \& Hairer in \cite{CH} assuming the noise satisfies some quantitative estimates on its cumulants. These assumptions are neither stronger nor weaker than the spectral assumption involved in the work \cite{HS23} of Hairer \& Steele. Our work shares the same inductive structure as \cite{HS23}, inherited from the seminal work \cite{LOTT} of Linares, Otto, Tempelmayr \& Tsatsoulis. The algebraic and analytic details are different though. 

To deal with models on a usual regularity structure and their first order Malliavin derivative in a single setting we introduce a useful variant of the notion of regularity structure that we call a {\it regularity-integrability structure}. It provides a setting where to measure each quantity ${\sf \Pi}_x\tau$ in a $\tau$-dependent way. Its symbol space is in particular graded by a subset of $\bbR \times [1,\infty]$, with the first component accounting for a regularity exponent and the second component accounting for an integrability exponent. An associated notion of modelled distribution, their reconstruction and Schauder estimates for some integration operator, can be developed in the regularity-integrability setting as in the classical setting. This is done in the companion work \cite{SemigroupMasato} and all the analytic tools proved in that work are summarized in Appendix \ref{SectionAppendix}.

We work in the flexible analytical setting of regularity-integrability structures introduced in \cite{SemigroupMasato} while \cite{HS23} work with pointed model distributions. 
Both approaches exploit the same regularity-integrability compensation associated with Cameron--Martin directions, but encode it differently. Hairer and Steele represent $d_\omega({\sf\Pi}_x\tau)(h)$ as the reconstruction of a modelled distribution pointed at $x$; the pointed norm records the improved local scaling near this distinguished point. We instead attach both a regularity exponent and an integrability exponent to each symbol. The pointed approach requires pointed reconstruction and Schauder estimates as non-trivial tasks. Our approach avoids distinguished-point-dependent spaces and treats the model and its first derivative in a common framework with a slight modification that $B_{\infty,\infty}^r$-type estimates are replaced with $B_{p,\infty}^r$-type ones.
Our use of preparation maps fits ideally the inductive nature of the proof of convergence and allows to prove it for an a priori very large class of renormalization procedures, not only for the BPHZ procedure. The qualitative stability result of Theorem \ref{ThmContinuity} was already proved in \cite{HS23}; its quantitative version in Theorem \ref{ThmLipschitzContinuity} is new in this setting.

\medskip

{\it \S5 The take home picture.}   It is possible to give a non-technical picture of the structure of the argument that we use to prove Theorem \ref{ThmMainIntro} in its proper form, Theorem \ref{thm:main}. Anticipating slightly over what follows, let us say here that we will work with two kinds of (algebraic) trees: the classical trees involving the noise symbol $\ocircle$ only, and a class of trees involving at most one noise derivative symbol $\odot$ in addition to the noise symbol. We collect the first trees in a finite collection $V$ and the second set of trees in a collection $W$. We will see that this can be done in a recursive way, building gradually some sets $W_1,V_1$ then $W_2,V_2$, etc. with $V=V_\ell$ and $W=W_\ell$ for some finite $\ell$. Using a familiar pictorial notation that we will re-introduce below, we would for instance have the following sets 
\begin{align*} 
\begin{aligned}
&W_1 = \spa\big(\{\odot\}\cup\{X^{k}\}_{{k}\in\bbN^d}\big), 
&&V_1 = \spa\big(\{\ocircle\}\cup\{X^{k}\}_{{k}\in\bbN^d}\big),   
\\
&W_2 = \spa\big(\big\{\ocircle, \odot, \begin{tikzpicture}[scale=0.3,baseline=0.05cm] \node at (0,0)  [noise] (1) {}; \node at (0,1.1)  [noise] (2) {}; \node at (0,1.1)  [dot] {}; \draw[K] (1) to (2); \end{tikzpicture}\hspace{0.03cm},\hspace{0.03cm} \begin{tikzpicture}[scale=0.3,baseline=0.05cm] \node at (0,0)  [noise] (1) {}; \node at (0,1.1)  [noise] (2) {}; \node at (0,0)  [dot] {}; \draw[K] (1) to (2); \end{tikzpicture}\hspace{0.03cm}\big\}\cup\{X^{k}\}_{{k}\in\bbN^d}\big),
&&V_2 = \spa\big(\big\{\ocircle, \begin{tikzpicture}[scale=0.3,baseline=0.05cm] \node at (0,0)  [noise] (1) {}; \node at (0,1.1)  [noise] (2) {}; \draw[K] (1) to (2); \end{tikzpicture}\hspace{0.03cm}\big\}
\cup\{X^{k}\}_{{k}\in\bbN^d}
\big),
\\
&W_3 = \spa\Big(\big\{\ocircle, \begin{tikzpicture}[scale=0.3,baseline=0.05cm] \node at (0,0)  [noise] (1) {}; \node at (0,1.1)  [noise] (2) {}; \draw[K] (1) to (2); \end{tikzpicture}\hspace{0.03cm}, \odot, \begin{tikzpicture}[scale=0.3,baseline=0.05cm] \node at (0,0)  [noise] (1) {}; \node at (0,1.1)  [noise] (2) {}; \node at (0,1.1)  [dot] {}; \draw[K] (1) to (2); \end{tikzpicture}\hspace{0.03cm},\hspace{0.03cm} \begin{tikzpicture}[scale=0.3,baseline=0.05cm] \node at (0,0)  [noise] (1) {}; \node at (0,1.1)  [noise] (2) {}; \node at (0,0)  [dot] {}; \draw[K] (1) to (2); \end{tikzpicture}\hspace{0.03cm}, \begin{tikzpicture}[scale=0.3,baseline=0.05cm] \node at (0,0)  [noise] (1) {}; \node at (-0.5,1.1) [noise] (2) {}; \node at (0.5,1.1) [noise] (3) {}; \node at (-0.5,1.1)  [dot] {}; \draw[K] (1) to (2); \draw[K] (1) to (3); \end{tikzpicture}\hspace{0.03cm} , \begin{tikzpicture}[scale=0.3,baseline=0.05cm] \node at (0,0)  [noise] (1) {}; \node at (-0.5,1.1) [noise] (2) {}; \node at (0.5,1.1) [noise] (3) {}; \node at (+0.5,1.1)  [dot] {}; \draw[K] (1) to (2); \draw[K] (1) to (3); \end{tikzpicture}\hspace{0.03cm}, \begin{tikzpicture}[scale=0.3,baseline=0.05cm] \node at (0,0)  [noise] (1) {}; \node at (-0.5,1.1) [noise] (2) {}; \node at (0.5,1.1) [noise] (3) {}; \node at (0,0)  [dot] {}; \draw[K] (1) to (2); \draw[K] (1) to (3); \end{tikzpicture}\hspace{0.03cm}\big\}\cup\{X^{k}\}_{{k}\in\bbN^d}\Big),  \quad
&&V_3 = \spa\Big(\big\{\ocircle, \begin{tikzpicture}[scale=0.3,baseline=0.05cm] \node at (0,0)  [noise] (1) {}; \node at (0,1.1)  [noise] (2) {}; \draw[K] (1) to (2); \end{tikzpicture}\hspace{0.03cm}, \begin{tikzpicture}[scale=0.3,baseline=0.05cm] \node at (0,0)  [noise] (1) {}; \node at (-0.5,1.1) [noise] (2) {}; \node at (0.5,1.1) [noise] (3) {}; \draw[K] (1) to (2); \draw[K] (1) to (3); \end{tikzpicture}\hspace{0.03cm}\big\}\cup\{X^{k}\}_{{k}\in\bbN^d}\Big).
\end{aligned}
\end{align*}
These sets would be used if we were studying the three-dimensional parabolic Anderson model equation
$$
(\partial_t-\Delta)u=f(u)\xi,\qquad t>0,\ x\in\bbT^d,
$$
where $f$ is a sufficiently regular function.
The existence of a limit model in Theorem \ref{ThmMainIntro} follows from some quantitative probabilistic estimates on the models ${\sf M}^n$. We show these estimates gradually by considering what happens on $W_1$ and $V_1$ first, then on $W_2$ and $V_2$, etc. In this inductive procedure, from the estimates on the models ${\sf M}^n$ on $W_i$ we can infer these estimates on $V_i$. The joint estimates on $W_i$ and $V_i$ then imply the estimates on $W_{i+1}$, and we can proceed inductively.

Different types of arguments are involved in each step. We use the spectral gap property of the law of $\xi$ to get the estimates on $V_i$ from the estimates on $W_i$. This is a truly probabilistic step. On the other hand, proving the estimates on $W_{i+1}$ from the estimates on $W_i$ and $V_i$ follows from some deterministic, pathwise, estimates. There is one subtle point in this second step. We happen indeed to be working with some regularity structures and some models indexed by some parameter $2\leq p\leq \infty$. The parameter of interest is $p=\infty$. We first prove the estimates on $W_{i+1}$ for the models with $p=2$ by an elementary reconstruction argument. Knowing in addition by induction that the estimates hold for all $2\leq p\leq \infty$ on $W_i$ allows then to obtain them on $W_{i+1}$ for all $2\leq p\leq \infty$.

\bigskip

{\it \S6 Reading guide.} We refer the reader to the reviews \cite{ChandraWeber, CorwinShen} of Chandra \& Weber and Corwin \& Shen for some non-technical introductions to the domain of semilinear singular stochastic PDEs. One can refer to the books \cite{FrizHairer, Berglund} of Friz \& Hairer and Berglund for mildly technical introductions to regularity structures, and to Bailleul \& Hoshino's Tourist's Guide \cite{RSGuide} for a thorough tour of the analytic and algebraic sides of the theory. Hairer's lecture notes \cite{HairerBPHZ,HairerTakagi} are centered on the problems of renormalisation in the setting of Feynmann graphs and in the setting of singular stochastic PDEs, respectively. The lecture notes \cite{OST, OttoStFlour} of Otto \& co. give a gentle introduction to the tree-free approach \cite{LOTT} to the renormalisation of the random models that are involved in the analytic and algebraic settings of \cite{OW,OSSW, LOT}. The present work is independent of any of these works.

\ssk

\noindent \textbf{{Organisation of the work.}} Section \ref{SectionFunctionalSetting} introduces a number of notations and objects used throughout the text. Section \ref{SectionResult} sets the scene for our main convergence result for random models. We specify the spectral gap assumption on the law of the random noise in Section \ref{SubsectionSpectralGap}. We introduce algebraic structures of decorated trees that depend on a parameter $p\in[2,\infty]$ in Section \ref{SubsectionDecoratedTrees}. A notion of differentiable sector is introduced in Section {{\ref{SubsectionDifferentialSectors}}}. It specifies a setting where one can talk of a sector that is stable by a natural noise-derivative operator. Section {{\ref{SubsectionRIS}}} introduces regularity-integrability structures. Since the degree of a tree depends on $p$ so does a regularity-integrability structure, that is, the algebraic rules for making local expansions depend on $p$. About models on regularity-integrability structures the main point to get is that we use different Lebesgue spaces to measure some quantities indexed by trees depending on whether or not there is a derivative noise symbol in that tree. The fundamental results about modelled distributions in the setting of regularity-integrability structures are proved in the companion work \cite{SemigroupMasato} and stated in Appendix \ref{SectionAppendix}. Preparation maps and their renormalized models are introduced in Section {{\ref{SubsectionRenormalizedModels}}}. Lemma \ref{LemComparisonFormula}, in this section, is important: It provides an explicit comparison for the fundamental quantities ${\sf\Pi}_x^n\tau$ for two different values of $p$, that is, when the local expansion rules are possibly different. We state our main result, Theorem {{\ref{thm:main}}, in Section {{\ref{SubsectionMainResult}}}.} The remaining two sections are dedicated to the proof of Theorem {{\ref{thm:main}}}. We present the inductive core of the proof in Section {{\ref{SectionMechanics}}} and defer to Section {{\ref{SectionLemmas}}} the proof of a number of lemmas used in the induction. In a nutshell, we first introduce an order to construct inductively a limit model on an increasing finite sequence of spaces $V_i, W_i$. The trees of $V_i$ have no derivative noise while the trees of $W_i$ may have one derivative noise. The induction proceeds in three steps after proving the convergence result for the base case. In Step $1$ we prove that the probabilistic convergence of the renormalized models on $W_i$ implies its probabilistic convergence on $V_i$. We use for that purpose the reconstruction theorem and the spectral gap inequality. In Steps $2$ and $3$ we prove that the convergence on $W_i$ and $V_i$ implies the convergence on $W_{i+1}$. In Step 2 we use the multilevel Schauder estimate to control the $\sf g$-part of the renormalized model on the elements of $W_{i+1}$. In Step 3 we use the comparison lemma, Lemma \ref{LemComparisonFormula}, to control the $\sf \Pi$-part of the model. Indeed one gets for free some analytic estimates when $p=2$. The comparison lemma then allows to compare the $\sf \Pi$-part for an arbitrary $p$ to its counterpart for $p=2$. It turns out that the difference between the two quantities involves only some terms whose control is provided by the induction mechanics. We prove Theorem \ref{ThmContinuity} and Theorem \ref{ThmLipschitzContinuity} in Section \ref{SectionThmContinuity}.

\medskip

We give in this work our convergence result in a situation where there is only one noise and one integration operator. The modifications needed to accommodate a situation with different noises and different integration operators, as in \cite{HS23}, are standard and left to the reader.

%--------------------------------------------------%
\section{Functional setting and notations}
\label{SectionFunctionalSetting}
%--------------------------------------------------%

 We fix the functional settings and some notations in this preliminary section. Throughout this paper, we fix an integer $d\ge1$ and the \emph{scaling} $\mfs=(\mfs_1,\dots,\mfs_d)\in[1,\infty)^d$. Set
$$
\vert\mfs\vert \defeq \mfs_1+\cdots+\mfs_d.
$$
Note that, in applications, $d$ stands for the dimension of spacetime. We denoted the time variable by $t$ and the spatial variable by $x$ in \eqref{*eq:cor:generalSPDE}, but from now on we rather write $x$ for the spacetime variable.
For any multiindex ${k}=(k_j)_{j=1}^d\in\bbN^d$ we define
$$
|{k}|_{\mfs}\defeq\sum_{j=1}^d\mfs_j k_j,\qquad k!\defeq\prod_{j=1}^dk_j!,\qquad
\partial^k\defeq\prod_{j=1}^d\partial_j^{k_j},
$$
where $\partial_j$ is the partial derivative with respect to the $j$-th variable $x_j$ on $\bbR^d$.
For $k,l\in\bbN^d$, we write $k\le l$ if $k_j\le l_j$ for any $j\in\{1,\dots,d\}$. Then define $\binom{l}{k}\defeq\prod_{j=1}^d\binom{l_j}{k_j}$.
Also for every $x=(x_j)_{j=1}^d\in\bbR^d$ and ${k}=(k_j)_{j=1}^d\in\bbN^d$ we define
$$
x^{k} \defeq \prod_{j=1}^d x_j^{k_j},  \qquad  \|x\|_\mfs\defeq\sum_{j=1}^d |x_j|^{\frac{1}{\mfs_j}}.
$$

We define the Besov spaces associated with a heat semigroup $(\mcQ_t)_{t>0}$ of an (anisotropic) elliptic operator. This $t$ is not a physical time variable, but is used only to define the regularities of (spacetime) functions/distributions.
Throughout this paper, we fix an exponent $\ell > \max_{1\le j\le d}\mfs_j$ and the differential operator $\mcL=P(\partial_1,\dots,\partial_d)$, where
\begin{equation}\label{sec1:ellipticop}
P(\lambda_1,\dots,\lambda_d)=\sum_{|{k}|_\mfs\le \ell} a_{k} \lambda^{k}
\end{equation}
is a polynomial with real constant coefficients which satisfies
\begin{equation}\label{sec1:ellipticity}
P(i\lambda_1,\dots, i\lambda_d) \leq -\delta \|\lambda\|_{\fr s}^\ell
\end{equation}
for some fixed $\delta>0$ and for any $\lambda\in\bbR^d$.
A typical example is the spacetime elliptic polynomial
$$
P(\lambda)=\lambda_1^2-(|\lambda'|^2-1)^2\defeq\lambda_1^2-(\lambda_2^2+\cdots+\lambda_d^2-1)^2,
$$
where $\lambda_1$ and $\lambda'=(\lambda_2,\dots,\lambda_d)$ represent time and space variables, respectively.
The factorization $P(\lambda)=(\lambda_1-|\lambda'|^2+1)(\lambda_1+|\lambda'|^2-1)$ gives a representation of the inverse of the heat operator $\partial_1-\Delta+1$, where $\Delta\defeq\partial_2^2+\cdots+\partial_d^2$, by
$$
\big(\partial_1-\Delta+1\big)^{-1} = -\int_0^\infty (\partial_1+\Delta-1) e^{t\mcL}dt.
$$
We denote by $(\mcQ_t=e^{t\mcL})_{t>0}$ the operator semigroup associated with $\mcL$ and by $Q_t(x)$ the integral kernel of $\mcQ_t$ determined by the relation

$$
(\mcQ_tf)(x) = \int_{\bbR^d}Q_t(x-y)f(y)dy
$$
for any compactly supported smooth functions $f$. We also write $(\mcQ_tf)(x)=\mcQ_t(x,f)$ later. In Appendix A of \cite{BHK22} it was proved that, for any $k\in\bbN^d$, $(\partial^kQ_t)_{t>0}$ satisfies the upper `Gaussian' estimate
$$
|\partial^kQ_t(x)|\lesssim \frac1{t^{|k|_\mfs/\ell}}G_t(x)\defeq\frac1{t^{(|\mfs|+|k|_\mfs)/\ell}}\exp\Bigg\{-c_1\sum_{j=1}^d\bigg(\frac{|x_j|}{t^{\mfs_j/\ell}}\bigg)^{\frac{\ell}{\ell-\mfs_j}}\Bigg\}
$$
for any $t>0$ and $x\in\bbR^d$, where $c_1>0$ is a constant.
We also define the family of weight functions $(w_c)_{c\ge0}$ on $\bbR^d$ by
$$
w_c(x) \defeq (1+\|x\|_\mfs)^{-c}.
$$
It is elementary to show that for any $c\ge0$ and $a\ge0$, the inequality
\begin{align}\label{ineq:Gw}
\big\Vert G_t(x) w_c(x)^{-1} \|x\|_\mfs^a\big\Vert_{L_x^p(\bbR^d)}
\lesssim t^{-\frac{|\mfs|}\ell(1-\frac1p)+\frac{a}\ell}
\end{align}
holds uniformly over $p\in[1,\infty]$ and $t\in(0,1]$. This estimate will be used in Section \ref{SectionLemmas}. For $c\ge0$ and $p\in[1,\infty]$ we define the weighted $L^p$ norm by
$$
\|f\|_{L^p(w_c)}\defeq\|fw_c\|_{L^p(\bbR^d)}.
$$
We follow \cite{BB16} for the definition of Besov spaces.

\ssk

\begin{defn}
For any nonnegative integer $a$ we set $\mcQ_t^{(a)}=(-t\mcL)^a\mcQ_t$. For every $r<a\ell$ and $p,q\in[1,\infty]$ we define the \textbf{{Besov space}} $B_{p,q}^{r,Q}(w_c)$ as the completion of $C(\bbR^d)\cap L^p(w_c)$ under the norm
\begin{equation} \label{EqDefnQBesovSpaces}
\|f\|_{B_{p,q}^{r,Q}(w_c)}^{(a)}\defeq\|\mcQ_1f\|_{L^p(w_c)}+\big\|t^{-r/\ell}\|\mcQ_t^{(a)}f\|_{L^p(w_c)}\big\|_{L^q((0,1];\frac{dt}t)}.
\end{equation}
\end{defn}

\ssk

The topological space $B_{p,q}^{r,Q}(w_c)$ is defined independently to the choice of $a$ as long as $a>r/\ell$, since the norms $\|\cdot\|_{B_{p,q}^{r,Q}(w_c)}^{(a)}$ and $\|\cdot\|_{B_{p,q}^{r,Q}(w_c)}^{(b)}$ are equivalent if $a,b>r/\ell$. We omit the proof of this fact since it is not actually used in this paper and is essentially the same as the proofs of Proposition 4.5 of \cite{BDY12} and Proposition A.3 of \cite{BB16}. See \cite{BDY12}, \cite{BB16}, or Section 2 of \cite{SemigroupMasato} for other detailed properties of Besov spaces associated with $(Q_t)_{t>0}$. Especially, the continuous embedding result
\begin{equation}\label{BesovEmbedding}
B_{p_1,q}^{r,Q}(w_c)\hookrightarrow B_{p_2,q}^{r-|\mfs|(\frac1{p_1}-\frac1{p_2}),Q}(w_c)
\end{equation}
for any $r\in\bbR$, $p_1,p_2,q\in[1,\infty]$ with $p_1\le p_2$, and $c\ge0$, is important. Write 
$$
H^{r,Q}(w_c) \defeq B_{2,2}^{r,Q}(w_c), \qquad C^{r,Q}(w_c) \defeq B_{\infty,\infty}^{r,Q}(w_c).
$$

Finally we fix the regularizing operator.
Pick $0<\ell_1<\ell$ and some real constants $(b_{l})_{\vert{l}\vert_{\frak{s}}<\ell_1}$ and set 
$$
\widetilde{K}_t(x) \defeq \sum_{\vert{l}\vert_{\frak{s}}\le\ell_1} b_{l} \, \partial_x^{l}Q_t(x).
$$
Then $(\widetilde{K}_t)_{0<t\le1}$ is an $(\ell-\ell_1)$-regularizing kernel in the terminology of Section 5 of \cite{SemigroupMasato}. As mentioned before a typical example is $\widetilde{K}_t=-(\partial_1+\Delta-1)e^{t(\partial_1^2-(\Delta-1)^2)}$. For technical reasons we fix a compactly supported smooth function $\chi$ which is equal to $1$ on a neighborhood of $0$ and we consider the modified kernels
$$
K_t = \big( 1-\chi(\partial) \big) \widetilde{K}_t,
$$
where 
$$
\chi(\partial)f\defeq\mcF^{-1}\big(\chi \, \mcF f\big)
$$ 
is a Fourier multiplier operator defined by Fourier transform $\mcF$ and its inverse $\mcF^{-1}$. This modification ensures that
$$
\int_{\bbR^d}x^k\partial^lK_t(x)dx = 0
$$
for any $k,l\in\bbN^d$. Since $\chi(\partial)$ is an operator with sufficient regularizing effects, the discussions in Section 5 of \cite{SemigroupMasato} still work with small modifications. Finally, for any $f\in C(\bbR^d)\cap L^p(w_c)$ and ${k}\in\bbN^d$ define
\begin{equation} \label{EqQSpacetimeOperator}
\partial^{k}\mcK(x,f) \defeq \int_0^1 \hspace{-0.15cm}\int_{\bbR^d} \partial_x^{k} K_t(x-y)f(y) \, dydt.
\end{equation}
When $\widetilde{K}_t = -(\partial_1+\Delta-1) \, e^{t(\partial_1^2-(\Delta-1)^2)}$, the operator $\mcK$ coincides with $(\partial_1-\Delta+1)^{-1}$ modulo some sufficient regularizing operators.

\ssk

For any normed vector space $X$ we denote by $L^q(\Omega;X)=L^q(\Omega,\bbP;X)$ the space of $q$-integrable $X$-valued random variables, for any $1\leq q\leq \infty$.

%\vfill \pagebreak

%-----------------------------------%
\section{A convergence result}
\label{SectionResult}
%-----------------------------------%

We introduce in this section all the notions involved in the formal statement and proof of Theorem \ref{ThmMainIntro}, given in Theorem \ref{thm:main} below. Of course, some familiarity with the practice of regularity structures is needed to take full profit of the present work; it is nonetheless formally self-contained. Theorem \ref{ThmMainIntro} involves a number of ingredients, the first of which being decorated trees. They serve as label set for the quantities that define a model over a regularity structure. We recall shortly in Section \ref{SubsectionDecoratedTrees} an ad hoc formalism to talk about decorated trees and some algebraic structures defined on this set that depend on two parameters $\varepsilon>0$ and $1\leq p\leq \infty$. The trees that we consider are the algebraic representation of some polynomial functionals of a random noise. In the setting we are working in, these polynomial functions of the noise can be differentiated in some sense. The algebraic counterpart of this differentiation operation is presented in Section \ref{SubsectionDifferentialSectors}, that deals with the notion of differentiable sector. 

The parameter $p$ is part of the data that define a regularity-integrability structure. These structures differ mainly from regularity structures by the fact that the index set grading the different sets of trees does not consist anymore of real numbers but rather of pairs $(r,i)$ of real numbers. The introduction of this structure is motivated by the fact that we want to measure the size of the analytical objects associated with each tree in some tree-dependent Besov space of $B^r_{i,\infty}$ type. Regularity-integrability structures and their associated notion of model were introduced in \cite{SemigroupMasato}; we recall in Section \ref{SubsectionRIS} what we need here from this setting. The convergence result of Theorem \ref{thm:main} involves a class of models on a regularity-integrability structure that are built recursively from some interpretation operators and some preparation maps. This construction is recalled in Section \ref{SubsectionRenormalizedModels}, where we also state an important lemma which makes clear what happens to these models when we vary the parameter $p$.

All this concerns the analytic and algebraic study of some models built from a smooth noise, random or not. The essential ingredient  in Theorem \ref{ThmMainIntro} that allows to extend this analysis to some non-smooth {\it random} noise is the spectral gap inequality from Section \ref{SubsectionSpectralGap}. It is a kind of Poincar\'e inequality satisfied by the law of the random noise of interest. We also state a fundamental lemma which relates the differentiation operation with respect to the noise to its algebraic counterpart.

We can then finally give in Section \ref{SubsectionMainResult} the statement of Theorem \ref{thm:main} and give a few comments about it. The proof of Theorem \ref{thm:main} is deferred to Section \ref{SectionMechanics} and Section \ref{SectionLemmas}.

\medskip

%%-----------------------%%
\subsection{Decorated trees$\boldmath .$   \hspace{0.03cm}}
\label{SubsectionDecoratedTrees}
%%-----------------------%%

Anticipating over Section \ref{SubsectionSpectralGap} where we introduce formally the spectral gap inequality, we construct in the present section the space of decorated trees corresponding to noise functionals following \cite{BHZ}, but we also include some trees corresponding to the first order Malliavin derivatives of these noise functionals. 
As will be defined more precisely in Section \ref{SubsectionSpectralGap}, $\Omega$ denotes a separable Banach space in which the noise $\xi$ takes its values, and $H$ denotes a separable Hilbert space continuously and densely embedded into $\Omega$. In this paper, derivatives of functionals on $\Omega$ along directions in $H$ play an important role.
We introduce the set of three edge labels 
$$
{\fr L}=\{\bbOmega,\bbH,\bbK\},
$$
whose elements will play in the sequel the role of an element of $\Omega$, an element of $H$, and the integral operator $\mcK$. 
We denote by $\widetilde{\bf T}$ the set of all 4-tuples $(\tau,{\fr t},{\fr n},{\fr e})$ of a non-planar rooted tree $\tau$ with vertex set $N_\tau$ and edge set $E_\tau$, a label map ${\fr t}:E_\tau\to{\fr L}$, and two decoration maps ${\fr n}:N_\tau\to\bbN^d$ and ${\fr e}:E_\tau\to\bbN^d$. Each element of $\widetilde{\bf T}$ is sometimes simply denoted by $\tau$. In the pictures, edges with labels $\bbOmega,\bbH$, and $\bbK$ are represented by dotted lines, double dotted lines and solid lines respectively. The $\fr n$- and $\fr e$-decorations assigned to nodes and edges are indicated by writing the values directly next to the corresponding nodes or edges. If a decoration is zero, it is omitted. Here is an example
$$
\rtds{1.5}{\pol{2}{1}{90}\pol{3}{2}{120}\pol{4}{2}{60}\pol{5}{3}{90}\pol{6}{4}{90}\pol{7}{1}{45}
\drl{1}{2}\drl{2}{3,4}\drdh{3}{5}\drdt{4}{6}\drdt{1}{7}\drb{1,2,3,4,5,6,7}
\node at ($(A7)+(0.13,0)$) {\tiny $k$};
\node at ($0.5*(A4)+0.5*(A6)+(0.06,0)$) {\tiny $l$};}
$$

Following \cite{BHZ} we define the algebraic structure on the linear space $\widetilde{T}$ spanned by $\widetilde{\bf T}$. First, the tree product of $\tau,\sigma\in\widetilde{\bf T}$, denoted by $\tau\sigma$, is obtained by identifying the roots of $\tau$ and $\sigma$ in the disjoint union $\tau\sqcup\sigma$, whose decorations are inherited from those of $\tau$ and $\sigma$ except that the $\fr n$-decoration at the root of $\tau\sigma$ is the sum of $\fr n$-decorations at the roots of $\tau$ and $\sigma$. Then $\widetilde{T}$ becomes an algebra with the linearly extended tree product. Second, the linear map 

$$
\widetilde{\Delta}:\widetilde{T}\to \widetilde{T}\otimes \widetilde{T}
$$ 
is defined by the following formula. Precisely, since $\widetilde{\Delta}$ may produce infinite sums, we have to consider the `bigraded space' spanned by $\widetilde{\bf T}$ and the tensor product of bigraded spaces. See Section 2.3 of \cite{BHZ} for the details. Since we later consider the truncation of $\widetilde{\Delta}$ into finite sums, we do not touch the details here.

\ssk

\begin{defn}
For each $(\tau,{\fr t},{\fr n},{\fr e})\in\widetilde{\bf T}$, the $\widetilde{\Delta}(\tau,{\fr t},{\fr n},{\fr e})\in\widetilde{T}\otimes\widetilde{T}$ is defined as the infinite sum
\begin{equation}\label{*eq:graphicalcoprod}
\sum_\sigma\sum_{{\fr n}_\sigma,{\fr e}_{\partial\sigma}}\frac1{{\fr e}_{\partial\sigma}!}\binom{{\fr n}}{{\fr n}_\sigma} \big(\sigma,{\fr t}\vert_{E_\sigma},{\fr n}_\sigma+\pi{\fr e}_{\partial\sigma},{\fr e}\vert_{E_\sigma}\big) 
\otimes
\Big(\tau/\sigma,{\fr t}\vert_{E_\tau\setminus E_\sigma},[{\fr n}-{\fr n}_\sigma]_\sigma,{\fr e}\vert_{E_\tau\setminus E_\sigma}+{\fr e}_{\partial\sigma}\Big),
\end{equation}
where
\begin{itemize}
	\item[--] $\sigma$ runs over all subtrees of $\tau$ which contain the root of $\tau$. Then the quotient tree $\tau/\sigma$ is obtained by identifying all nodes of $\sigma$ in the tree $\tau$. The edge set of $\tau/\sigma$ is $E_\tau\setminus E_\sigma$.   \vspace{0.1cm}

	\item[--] ${\fr n}_\sigma$ runs over all maps $N_\sigma\to\bbN^d$ such that ${\fr n}_\sigma(v)\le{\fr n}(v)$ for any $v\in N_\sigma$. The map 
	$$
	[{\fr n}-{\fr n}_\sigma]_\sigma:N_{\tau/\sigma}\to\bbN^d
	$$ 
	is defined by $[{\fr n}-{\fr n}_\sigma]_\sigma(v)={\fr n}(v)$ for non-root $v\in N_{\tau/\sigma}$ and by $[{\fr n}-{\fr n}_\sigma]_\sigma(\varrho)=\sum_{v\in N_\sigma}({\fr n}(v)-{\fr n}_\sigma(v))$ for the root $\varrho$ of $\tau/\sigma$. Moreover one sets 
	$$
	\binom{{\fr n}}{{\fr n}_\sigma}\defeq\prod_{v\in N_\sigma}\binom{{\fr n}(v)}{{\fr n}_\sigma(v)}.
	$$ %  \vspace{0.1cm}

	\item[--] $\partial\sigma$ is the boundary of $\sigma$, that is, the set of all edges $(u,v)\in E_\tau$ such that $u\in N_\sigma$ and $v\in N_\tau\setminus N_\sigma$. (Every edge is represented as an ordered pair $(u,v)$ of two nodes, where $u$ is nearer to the root than $v$.) ${\fr e}_{\partial\sigma}$ runs over all maps $\partial\sigma\to\bbN^d$.  For any $u\in N_\sigma$, $\pi{\fr e}_{\partial\sigma}(u)\defeq\sum_{(u,v)\in\partial\sigma}{\fr e}_{\partial\sigma}((u,v))$. Moreover one sets 
	$$
	{\fr e}_{\partial\sigma}!\defeq\prod_{e\in\partial\sigma}{\fr e}_{\partial\sigma}(e)!
	$$
\end{itemize}
\end{defn}

\ssk

As shown in Proposition 3.23 of \cite{BHZ}, $\widetilde{T}$ is a Hopf algebra (as a bigraded space) with the tree product and the above coproduct. In the present paper, the recursive representation of $\widetilde{\Delta}$ is more convenient. A 4-tuple $(\bullet,\varnothing,k,\varnothing)$ of a single node with decoration ${\fr n}(\bullet)=k$ is denoted by $X^k$. In particular, we write ${\bf1}\defeq X^0$ and $X_j\defeq X^{e_j}$, where $e_j=(0,\dots,1,\dots,0)$ is the $j$-th standard unit normal vector of $\bbR^d$. We write $\mcI_{k}^{\fr l}(\tau)$ for the planted tree obtained from $\tau\in\widetilde{\bf T}$ by grafting it to a new root with $\frak{n}$-decoration $0$, along an edge with label ${\fr l}\in{\fr L}$ and $\frak{e}$-decoration ${k}\in\bbN^d$. Then $\widetilde{\Delta}:\widetilde{T}\to \widetilde{T}\otimes \widetilde{T}$ is calculated by the following formulas.
\begin{equation}\label{*eq:recursivecoprod}
\begin{aligned}
&\widetilde{\Delta}(X^{k}) = \sum_{l\le k} \binom{k}{l}X^{l}\otimes X^{k-l},   
\qquad\widetilde{\Delta}(\tau\sigma)=(\widetilde{\Delta}\tau)(\widetilde{\Delta}\sigma),\\
&\widetilde{\Delta}(\mcI_{k}^{\fr l}(\tau)) = (\mcI_{k}^{\fr l}\otimes\id)\widetilde{\Delta}\tau + \sum_{{l}\in\bbN^d}\frac{X^{l}}{{l}!}\otimes\mcI_{{k}+{l}}^{\fr l}(\tau).
\end{aligned}
\end{equation}
In the last equality the map $\mcI_k^{\fr l}:\widetilde{T}\to \widetilde{T}$ is regarded as the linear extension. (See Proposition 4.17 of \cite{BHZ} for a proof.)

For technical reasons, we do not consider trees that have a $\bbK$-labeled leaf. A {\it leaf} is an edge $e=(u,v)$ such that there is no outgoing edge from $v$. This influences the fact that, $\mcK$-admissible models ${\sf\Pi}$ defined in Section \ref{SubsectionRenormalizedModels} later satisfy ${\sf\Pi}(\mcI_k^{\bbK}X^l)(x)=\partial^k\mcK(x,(\cdot)^l)=0$ for any $k,l\in\bbN^d$ by the modification of $K_t$ given in Section \ref{SectionFunctionalSetting}. Because of this, we actually consider the quotient Hopf algebra
$$
T=\widetilde{T}/I,
$$
where $I$ is the Hopf ideal spanned by trees that have a $\bbK$-labeled leaf. 
We identify $T$ with the subspace spanned by $\bf T$, the set of all trees without $\bbK$-labeled leaves. Then the tree product $T\otimes T\to T$ and the coproduct 
$$
\Delta\defeq(\pi\otimes\pi)\circ\widetilde{\Delta}\circ\iota:T\to T\otimes T
$$
are well-defined, where $\iota:T\to\widetilde{T}$ is an inclusion map and $\pi:\widetilde{T}\to T$ is a canonical surjection.

\ssk

As in \cite{Hai14}, each tree is assigned an exponent representing the regularity of the random distribution it corresponds to. Then one might think that the label $\bbH$ should be assigned the exponent $-s_0$ as an element of $H=H^{-s_0,Q}(w_0)$, but this is not sufficient for the proof of the main result. One of the most important ideas in this paper is to give $\bbH$ the role as an element of the space $B_{p,\infty}^{r_0+\frac{|\mfs|}p,Q}(w_c)$ for any $p\in[2,\infty]$. By the Besov embedding \eqref{BesovEmbedding}, we see that this space is actually an intermediate space between $H$ and $\Omega$
\begin{equation}\label{*eq:interpolateHOmega}
H=H^{-s_0,Q}(w_0)\hookrightarrow B_{2,\infty}^{r_0+\frac{|\mfs|}2,Q}(w_c)\hookrightarrow B_{p,\infty}^{r_0+\frac{|\mfs|}p,Q}(w_c)\hookrightarrow C^{r_0,Q}(w_c)=\Omega,
\end{equation}
for any $r_0<-\frac{|\mfs|}2-s_0$. We define the family of degree maps as follows. The parameter $\varepsilon$ will later express an infinitesimal loss of regularity at each step of some inductive reasoning.

\ssk

\begin{defn}
Fix $s_0\geq-\frac{|\mfs|}2$, $r_0<-\frac{|\mfs|}2-s_0$, and $\beta_0\in(0,\ell-\ell_1)$. Assume that 
$$
r_0\notin\bbQ[\beta_0,\mfs] = \Big\{ r\beta_0+\sum_{j=1}^dq_j\mfs_j\,;\,r,q_1,\dots,q_d\in\bbQ \Big\}.
$$
For any parameters $\varepsilon\ge0$ and $p\in[2,\infty]$ we define the degree map $r_{\varepsilon,p}:{\bf T}\to\bbR$ by
\begin{align*}
&r_{\varepsilon,p}(\bbOmega) \defeq r_0-\varepsilon,\qquad r_{\varepsilon,p}(\bbH) \defeq r_0-\varepsilon+\frac{|{\fr s}|}p, \qquad r_{\varepsilon,p}(\bbK) \defeq \beta_0,
\\
&r_{\varepsilon,p}(\tau_{\fr e}^{\fr n}) \defeq \sum_{v\in N_\tau}|{\fr n}(v)|_{\fr s} + \sum_{e\in E_\tau} \big(r_{\varepsilon,p}({\fr t}(e))-|{\fr e}(e)|_{\fr s}\big).
\end{align*}
\end{defn}

\ssk

The letter $r$ is chosen for the word `regularity' as the degree of a tree is intimately related to the regularity of some analytic object associated with it. 
The technical constraint that $r_0\notin\bbQ[\beta_0,\mfs]$ is used in Lemma \ref{*epsilonuseless}, which ensures the rigidity of the structure when $\varepsilon$ varies slightly.
For each $\varepsilon\ge0$ and $p\in[2,\infty]$ we introduce the projection map $P^+_{\varepsilon,p}$ from $T$ to the subalgebra $T^+_{\varepsilon,p}$ spanned by the symbols   %\vspace{-0.2cm}
\begin{align} \label{baseofT+2}
X^{k}\prod_{i=1}^n\mcI_{{k}_i}^{{\fr l}_i}(\tau_i)
\end{align}
with $n\in\bbN$, ${k},{k}_i\in\bbN^d$, ${\fr l}_i\in{\fr L}$ and $\tau_i\in{\bf T}$ such that $r_{\varepsilon,p}(\tau_i)+r_{\varepsilon,p}({\fr l}_i)>|{k}_i|_{\fr s}$ for each $i$. We define
\begin{equation} \begin{split} \label{EqDefnDeltap}
&\Delta_{\varepsilon,p} \defeq(\id\otimes \mathop{P^+_{\varepsilon,p}})\Delta:T\to T\otimes T_{\varepsilon,p}^+,   \\
&\Delta^+_{\varepsilon,p} \defeq \big(P^+_{\varepsilon,p}\otimes P^+_{\varepsilon,p}\big)\Delta:T_{\varepsilon,p}^+\to T_{\varepsilon,p}^+\otimes T_{\varepsilon,p}^+.
\end{split} \end{equation}
The ranges are algebraic tensor products since $P_{\varepsilon,p}^+$ restricts the choices of ${\fr e}_{\partial\sigma}$ in \eqref{*eq:graphicalcoprod}. By a similar argument to \cite{Hai14, BHZ} we can see that $T^+_{\varepsilon,p}$ is a Hopf algebra with coproduct $\Delta^+_{\varepsilon,p}$ and $T$ has a right comodule structure with coaction $\Delta_{\varepsilon,p}$. Denote by $S^+_{\varepsilon,p}$ the antipode of $(T^+_{\varepsilon,p},\Delta^+_{\varepsilon,p})$. 

\ssk

We provide an elementary example describing the $p$-dependence of $\Delta_{\varepsilon,p}$. Under the choices $\mfs=(1,1,1)$, $r_0=-\frac32$, and $\beta_0=2$ we have

\begin{align}\label{ExL6-}
\Delta_{\varepsilon,p}
\Bigg(\, \rtdbs{7}{1.3}{\pol{2}{1}{90}\pol{3}{1}{30}\pol{4}{2}{30}
\drl{1}{2}\drdt{1}{3}\drdh{2}{4}\drb{1,2,3,4}}\, \Bigg)
= \left\{\begin{aligned}
&\, \rtdbs{5}{1.3}{\pol{2}{1}{90}\pol{3}{1}{30}\pol{4}{2}{30}
\drl{1}{2}\drdt{1}{3}\drdh{2}{4}\drb{1,2,3,4}}\,  \otimes\bullet
- \, \rtdbs{5}{1.3}{\pol{2}{1}{90}\drdt{1}{2}\drb{1,2}}\, 
\otimes \, \rtdbs{5}{1.3}{\pol{2}{1}{90}\pol{3}{2}{30}
\drl{1}{2}\drdh{2}{3}\drb{1,2,3}}\,, 
& \quad \textrm{ for } p\geq\frac{6}{1+2\varepsilon},   \\
&\, \rtdbs{5}{1.3}{\pol{2}{1}{90}\pol{3}{1}{30}\pol{4}{2}{30}
\drl{1}{2}\drdt{1}{3}\drdh{2}{4}\drb{1,2,3,4}}\,  \otimes\bullet
- \, \rtdbs{5}{1.3}{\pol{2}{1}{90}\drdt{1}{2}\drb{1,2}}\, 
\otimes \, \rtdbs{5}{1.3}{\pol{2}{1}{90}\pol{3}{2}{30}
\drl{1}{2}\drdh{2}{3}\drb{1,2,3}}\,
-\sum_{i=1}^3 \rtdbs{5}{1.3}{\pol{2}{1}{90}\drdt{1}{2}\drb{1,2}\node at ($(A1)-(0.17,0)$) {\tiny $e_i$};}\, 
\otimes \, \rtdbs{5}{1.3}{\pol{2}{1}{90}\pol{3}{2}{30}
\drl{1}{2}\drdh{2}{3}\drb{1,2,3}\node at ($0.5*(A1)+0.5*(A2)+(0.13,0)$) {\tiny $e_i$};}\,
&  \quad \textrm{ for } p<\frac{6}{1+2\varepsilon}.
\end{aligned} \right.
\end{align}
The number $p=\frac6{1+2\varepsilon}$ is the unique solution of the equation $r_{\varepsilon,p}\Big(\, \rtdbs{5}{1}{\pol{2}{1}{90}\pol{3}{2}{30} \drl{1}{2}\drdh{2}{3}\drb{1,2,3}\node at ($0.5*(A1)+0.5*(A2)+(0.13,0)$) {\tiny $e_i$};}\,\Big)=\frac3p-\frac{1+2\varepsilon}2=0$. In general, for fixed $\tau$ and $\varepsilon$, the number of terms in the decomposition of $\Delta_{\varepsilon,p}\tau$ changes at $p$ such that $r_{\varepsilon,p}(\mu)=0$ for some $\mu$ appearing in the right part of the tensor products in \eqref{*eq:graphicalcoprod}. Moreover $\Delta_{\varepsilon,p}\tau$ does not change between such `phase transition' points. This observation will be useful in Lemma \ref{LemComparisonFormula} below.

\medskip

To help the reader figuring out easily the nature of some sets we will systematically use below some bold capital letters to denote a finite collection and some non-bold capital letter to denote a vector space.

\ssk
%\medskip

%%----------------------------------%%
\subsection{Differentiable sectors$\boldmath .$ \hspace{0.05cm}}
\label{SubsectionDifferentialSectors}
%%----------------------------------%%

We define the formal Malliavin differentiation map. We denote by ${\bf T}^{(n)}$ be the set of all $\tau\in{\bf T}$ that have exactly $n$ edges with label $\bbH$. For each $\tau\in{\bf T}^{(0)}$ and $e\in E_\tau$ with label $\bbOmega$, we define $D_e\tau\in{\bf T}^{(1)}$ as the same tree as $\tau$ except that the label at $e$ has been replaced with $\bbH$. Also we define the linear map $D:\spa({\bf T}^{(0)})\to \spa({\bf T}^{(1)})$ by setting for any $\tau\in{\bf T}^{(0)}$
\begin{equation} \label{EqDefnDerivation}
D\tau \defeq \sum_{e\in E_\tau,\, {\fr t}(e)=\bbOmega} D_e\tau
\end{equation}
if ${\fr t}^{-1}(\bbOmega)\neq\emptyset$, and $D\tau\defeq0$ otherwise. 
Equivalently, the map $D$ is defined by the following formulas.
$$
DX^k=0,\qquad
D(\mcI_k^\bbOmega\tau)=\mcI_k^\bbH\tau,\qquad
D(\mcI_k^\bbK\tau)=\mcI_k^\bbK(D\tau),\qquad
D(\tau\sigma)=(D\tau)\sigma+\tau(D\sigma).
$$
Here is an example without $\fr n$ and $\fr e$ decorations
$$
D\Bigg(\, \rtdb{13}{\pol{2}{1}{90}\pol{3}{2}{120}\pol{4}{2}{60}\pol{5}{3}{90}\pol{6}{4}{90}\pol{7}{1}{45}
\drl{1}{2}\drl{2}{3,4}\drdt{3}{5}\drdt{4}{6}\drdt{1}{7}\drb{1,2,3,4,5,6,7}}\, \Bigg)
=
2\,\rtdb{13}{\pol{2}{1}{90}\pol{3}{2}{120}\pol{4}{2}{60}\pol{5}{3}{90}\pol{6}{4}{90}\pol{7}{1}{45}
\drl{1}{2}\drl{2}{3,4}\drdt{3}{5}\drdh{4}{6}\drdt{1}{7}\drb{1,2,3,4,5,6,7}}\,
+\,
\rtdb{13}{\pol{2}{1}{90}\pol{3}{2}{120}\pol{4}{2}{60}\pol{5}{3}{90}\pol{6}{4}{90}\pol{7}{1}{45}
\drl{1}{2}\drl{2}{3,4}\drdt{3}{5}\drdt{4}{6}\drdh{1}{7}\drb{1,2,3,4,5,6,7}}\,.
$$
Note that $D$ preserves the $r_{\varepsilon,\infty}$-degree of trees ($r_{\epsilon,\infty}(\tau) = r_{\epsilon,\infty}(D\tau)$) because $r_{\varepsilon,\infty}(\bbH)=r_{\varepsilon,\infty}(\bbOmega)$, but $D$ may change the $r_{\varepsilon,p}$-degree when $p$ is finite.

\ssk

For any subset ${\bf S}\subset{\bf T}^{(0)}$ we define $\dt{\bf S}\subset{\bf T}^{(1)}$ by
$$
\dt{\bf S} \defeq \big\{D_e\tau\ ;\, \tau\in {\bf S},\ e\in{\fr t}^{-1}(\bbOmega)\big\}.
$$
We summarize the properties of the sector of $T$ that we will study in the present paper.

\ssk

\begin{defn} \label{asmpV}
For a finite set ${\bf B}\subset{\bf T}^{(0)}$ we define
\begin{equation*} \begin{split}
V&\;\defeq\spa({\bf B}),   \\
W&\;\defeq\spa\hspace{-0.07cm}\big({\bf B}\cup\dt{\bf B}\big).
\end{split} \end{equation*}
For each $(\varepsilon,p)\in[0,\infty)\times[2,\infty]$ we define $V_{\varepsilon,p}^+$ as the subalgebra of $T_{\varepsilon,p}^+$ generated by the symbols 
$$
{\bf V}_{\varepsilon,p}^+\defeq\{X_j\}_{j=1}^d\cup\{\mcI_k^{\bbK}(\tau)\}_{k\in\bbN^d,\, \tau\in{\bf B}\setminus\{X^l\}_l,\, r_{\varepsilon,p}(\tau)+\beta_0>|k|_\mfs},
$$
and we define $W_{\varepsilon,p}^+$ as the subalgebra of $T_{\varepsilon,p}^+$ generated by the symbols 
$$
{\bf W}_{\varepsilon,p}^+\defeq\{X_j\}_{j=1}^d \cup \big\{\mcI_k^{\bbH}(\bullet)\big\}_{k\in\bbN^d,\, r_{\varepsilon,p}(\bbH)>|k|_\mfs} \cup \big\{\mcI_k^{\bbK}(\tau)\big\}_{k\in\bbN^d,\, \tau\in({\bf B}\cup\dt{\bf B})\setminus\{X^l\}_l,\, r_{\varepsilon,p}(\tau)+\beta_0>|k|_\mfs}.
$$
We call ${\bf B}$ a \textbf{differentiable sector} if it satisfies the following properties.
\newcounter{mymemory}
\begin{enumerate} \setlength{\itemsep}{0.1cm}
\renewcommand{\labelenumi}{\textrm{\bf (\alph{enumi})}}
	\item\label{asmpV0} ${\bf B}=\{X^k\}_{|k|_\mfs<L}\cup{\bf B}_{\ocircle}$ for a fixed $L>0$ and a set ${\bf B}_{\ocircle}$ of trees that have at least one $\bbOmega$-labeled edge. The tree $\mcI_0^{\bbOmega}(\bullet)=\, \rtdb{1}{\pol{2}{1}{90}\drdt{1}{2}\drb{1,2}}\, $ is contained in ${\bf B}_{\ocircle}$.
	
	\item\label{asmpleaf} For each $(\tau,{\fr t},{\fr n},{\fr e})\in{\bf B}_\ocircle$, if $e=(u,v)\in {\fr t}^{-1}(\bbOmega)$, then $e$ is a leaf, ${\fr n}(v)={\fr e}(e)=0$, and there are no other edges $e'=(u,v')\in{\fr t}^{-1}(\bbOmega)$ with $v'\neq v$.
	
	\item\label{asmpV1} $\Delta_{0,2}(V) \subset V\otimes V_{0,2}^+$.
	
	\item\label{asmpV2} $\Delta_{0,2}(W) \subset W\otimes W_{0,2}^+$.
	\setcounter{mymemory}{\value{enumi}}
\end{enumerate}
\end{defn}

\medskip

In plain words, $V$ is the linear span of the collection $\bf B$ of noise functionals that we originally sought to estimate, whereas $W$ is the linear span of these functionals together with their first-order Malliavin derivatives.
We choose the letter ${\bf B}$ for `basis'. The number $L>0$ in {\bf (\ref{asmpV0})} is such that we can establish the solution theory of a given SPDE within the regularity structure constructed from $\bf B$. We do not keep track of the dependence of $V$ on ${\bf B}$ in the notation as ${\bf B}$ will be fixed in the sequel. The following observations are easy but important.
\begin{enumerate}\setlength{\itemsep}{0.1cm}
\setcounter{enumi}{\value{mymemory}}
\renewcommand{\labelenumi}{\textrm{\bf (\alph{enumi})}}

\item\label{asmpleaf'}
By the property {\bf (\ref{asmpleaf})}, for each $\tau\in{\bf B}\cup\dt{\bf B}$ all the $\bbOmega$- or $\bbH$- labeled edges can be contracted into labeled nodes. We write
$$
\ocircle=\mcI_0^{\bbOmega}(\bullet)=\,\rtdb{1}{\pol{2}{1}{90}\drdt{1}{2}\drb{1,2}}\, ,\qquad
\odot=\mcI_0^{\bbH}(\bullet)=\,\rtdb{1}{\pol{2}{1}{90}\drdh{1}{2}\drb{1,2}}\,.
$$
Here is an example without $\fr n$ and $\fr e$ decorations
$$
\rtdb{13}{\pol{2}{1}{90}\pol{3}{2}{120}\pol{4}{2}{60}\pol{5}{3}{90}\pol{6}{4}{90}\pol{7}{1}{45}
\drl{1}{2}\drl{2}{3,4}\drdt{3}{5}\drdh{4}{6}\drdt{1}{7}\drb{1,2,3,4,5,6,7}}\,
=
\,
\rtdb{13}{\pol{2}{1}{90}\pol{3}{2}{120}\pol{4}{2}{60}
\drl{1}{2}\drl{2}{3,4}\drc{1,3}\drb{2}\drch{4}}\,.
$$
The $\fr n$-decorations at the roots of some contracted edges are inherited. Some edges with $\bbH$-label and nonzero $\fr e$-decoration can appear in elements of ${\bf W}_{\varepsilon,p}^+$. We write then
$$
\odot_k=\mcI_k^{\bbH}(\bullet).
$$
Do not confuse the label $k$ with the $\fr n$-decoration. Here $\odot_k$ plays the role of $\partial^kh$, for some function $h$ in some Hilbert space, while the $\frak{n}$-decoration plays the role of a polynomial map.

\item\label{asmp_Veps}
Since $p$ has no influences on the trees in ${\bf T}^{(0)}$ the space $V_{\varepsilon,p}^+$ is independent of $p$. We write 
$$
{\bf V}_\varepsilon^+\defeq {\bf V}_{\varepsilon,p}^+\qquad \textrm{and}\qquad V_\varepsilon^+\defeq V_{\varepsilon,p}^+.
$$

\item 
Since $r_{\varepsilon,p}(\tau)$ is monotonically decreasing in terms of $\varepsilon$ and $p$, one has 
$$
P_{\varepsilon,p}^+V_{0,2}^+=V_{\varepsilon,p}^+,\qquad P_{\varepsilon,p}^+W_{0,2}^+=W_{\varepsilon,p}^+,  \qquad(\id\otimes P_{\varepsilon,p}^+)\Delta_{0,2}=\Delta_{\varepsilon,p}.
$$ 
Therefore the properties {\bf (\ref{asmpV1})} and {\bf (\ref{asmpV2})} still hold if we replace the index $(0,2)$ by any $(\varepsilon,p)$ with $\varepsilon>0$ and $2\leq p\leq \infty$.
\end{enumerate}

\ssk

We note also that the properties {\bf (\ref{asmpV1})} and {\bf (\ref{asmpV2})} imply the following facts.

\ssk

\begin{lem}
For a differentiable sector $\bf B$ one has the following properties for any $(\varepsilon,p)\in[0,\infty)\times[2,\infty]$.
\begin{itemize} \setlength{\itemsep}{0.1cm}
	\item[{\bf (\ref{asmpV1}')}] The set $V_{\varepsilon}^+$ is a sub-Hopf algebra of $T_{\varepsilon,p}^+$, and $V$ is a right comodule over $V_\varepsilon^+$ with coaction $\Delta_{\varepsilon,p}$.
	
	\item[{\bf (\ref{asmpV2}')}] The set $W_{\varepsilon,p}^+$ is a sub-Hopf algebra of $T_{\varepsilon,p}^+$, and $W$ is a right comodule over $W_{\varepsilon,p}^+$ with coaction $\Delta_{\varepsilon,p}$.
\end{itemize}
\end{lem}

\ssk

\begin{Dem}
Since $\Delta_{\varepsilon,p}(V)\subset V\otimes V_\varepsilon^+$, we have that $\Delta^+_{\varepsilon,p}(V_{\varepsilon}^+) \subset V_{\varepsilon}^+\otimes V_{\varepsilon}^+$ from the definition of $V_\varepsilon^+$ and the recursive definition \eqref{*eq:recursivecoprod} of $\Delta$. Moreover, since the antipode $S_{\varepsilon,p}^+$ of $\Delta_{\varepsilon,p}^+$ is determined by recursive formulas
\begin{align*}
S_{\varepsilon,p}^+(X^k)&=(-X)^k,\qquad
S_{\varepsilon,p}^+(\tau\sigma)=S_{\varepsilon,p}^+(\tau)S_{\varepsilon,p}^+(\sigma),\\
S_{\varepsilon,p}^+\mcI_k^{\fr l}(\tau)&=-\sum_{l\in\bbN^d}\frac{(-X)^l}{l!}\mcM(P_{\varepsilon,p}^+\mcI_{k+l}^{\fr l}\otimes S_{\varepsilon,p}^+)\Delta_{\varepsilon,p}(\tau)
\qquad(r_{\varepsilon,p}(\mcI_k^{\fr l}(\tau))>0)
\end{align*}
(see Proposition 4.18 of \cite{BHZ}), where $\mcM$ denotes the tree product, we also have that $S_{\varepsilon,p}^+V_\varepsilon^+\subset V_\varepsilon^+$ by an induction on the size of trees. Hence we have {\bf (\ref{asmpV1}')}. The proof of {\bf (\ref{asmpV2}')} is similar.
\end{Dem}

\ssk

Given a complete subcritical rule (see \cite[Section 5]{BHZ}), the set $\bf B$ of trees strongly conforming this rule is a differentiable sector -- see Proposition \ref{*prop:ExampleDifferentiable} below.
In the sequel we fix a differentiable sector $V$ and its basis ${\bf B}$. The defining properties of a differentiable sector ensure that the tuple 
$$
\scV_\varepsilon \defeq \big((V,\Delta_{\varepsilon,p}), (V_{\varepsilon}^+,\Delta^+_{\varepsilon,p})\big)
$$
defines a concrete regularity structure in the sense of \cite{RSGuide}. The structure of the $(\varepsilon,p)$-dependent tuple
$$
\scW_{\varepsilon,p} \defeq \big((W,\Delta_{\varepsilon,p}), (W_{\varepsilon,p}^+, \Delta^+_{\varepsilon,p})\big),
$$
is encapsulated in a useful variation on the notion of regularity structure that we call a {\it regularity-integrability structure}. It is described in Section \ref{SubsectionRIS}.

\ssk

We close this section with an observation on the $(\eps,p)$-dependence of $\Delta_{\eps,p}$. Recall that the structure of $\Delta_{\varepsilon,p}\tau$ changes when $r_{\varepsilon,p}(\tau)=0$ for some $\tau\in{\bf W}_{0,2}^+\setminus\{X_j\}_{j=1}^d$.
For each $\tau\in{\bf W}_{0,2}^+\setminus\{X_j\}_{j=1}^d$, since $r_{\eps,p}(\tau)=r_{0,\infty}(\tau)-n_\tau\eps+\frac{h_\tau}{p}$ for some positive integer $n_\tau$ and $h_\tau\in\{0,1\}$, the set 
$$
\ell_\tau = \Big\{ \Big(\eps,\frac1p\Big)\,;\,r_{\eps,p}(\tau)=0 \Big\}
$$ 
forms a line as in the following picture. If $\tau\in{\bf T}^{(1)}$, the line $\ell_\tau$ is upward-sloping. Otherwise, $\ell_\tau$ is vertical to the $\eps$-axis.
\begin{center}
\begin{tikzpicture}[scale=1]
\draw[->,>=stealth,semithick] (0,0)--(5,0) node[right]{$\eps$};
\draw[->,>=stealth,semithick] (0,0)--(0,4) node[left]{$\frac1p$};
\draw[semithick] (0,3)--(5,3);
\draw (0,0) node[left]{$p=\infty$};
\draw (0,3) node[left]{$p=2$};
\draw[semithick] (0,0.2)--(3,4) node[right]{$\ell_{\tau_1}$};
\draw[semithick] (0,0.8)--(5,4) node[right]{$\ell_{\tau_2}$};
\draw[semithick] (0,2)--(2.5,4) node[above]{$\ell_{\tau_3}$};
\draw[semithick] (2,0)--(5,3.2) node[right]{$\ell_{\tau_4}$};
\draw[semithick] (4,0)--(4,4) node[above]{$\ell_{\tau_5}$};
\draw[semithick, dotted] (0.9574,0)--(0.9574,3) node[above]{};
\fill (0.9574,1.413) circle (2pt);
\node[below] at (0.9574,0) {$\eps_0$};
\fill (1.25,3) circle (2pt);
\fill (2,0) circle (2pt);
\draw[semithick, dotted] (0.4,3)--(0.4,0) node[below]{$\eps$};
\fill (0.4,0.7067) circle (2pt);
\fill (0.4,1.056) circle (2pt);
\fill (0.4,2.32) circle (2pt);
\draw[semithick, dotted] (0.4,0.7067)--(0,0.7067) node[left]{\scriptsize$p_{\eps}(\tau_1)^{-1}$};
\draw[semithick, dotted] (0.4,1.056)--(0,1.056) node[left]{\scriptsize$p_{\eps}(\tau_2)^{-1}$};
\draw[semithick, dotted] (0.4,2.32)--(0,2.32) node[left]{\scriptsize$p_{\eps}(\tau_3)^{-1}$};
\end{tikzpicture}
\end{center}
When going across the line $\ell_\tau$ from left to right, the sign of $r_{\eps,p}(\tau)$ changes from positive to negative. The above figure gives some examples of lines $\ell_\tau$ for some different trees $\tau_i$. For each fixed $\eps$ we write
$$
{\sf I}_\varepsilon \defeq \Big\{ p_\eps(\mu)\,;\,\mu\in{\bf W}_{0,2}^+\cap{\bf T}^{(1)},\ r_{\eps,\infty}(\mu)\le 0\le r_{\eps,2}(\mu) \Big\},
$$
for the points of discontinuity of the function $p\mapsto P_{\eps,p}^+$, where $p_\eps(\mu)\in[2,\infty]$ is the unique $p$ for which $r_{\eps,p}(\mu)=0$. Since $r_{\eps,p}(\mu)=r_{\eps,\infty}(\mu)+\frac{|\mfs|}p$ we have
$$
p_\varepsilon(\mu)= \frac{|{\fr s}|}{-r_{\eps,\infty}(\mu)}.
$$
Similarly, for each fixed $p$ we write
$$
{\sf J}_p \defeq \Big\{ \eps\ge0\,;\,\mu\in{\bf W}_{0,2}^+\cap{\bf T}^{(1)},\ r_{\eps,p}(\mu)=0 \Big\}
$$
for the point of discontinuity of the function $\eps\mapsto P_{\eps,p}^+$. Since the set ${\bf W}_{0,2}^+$ is finite, the numbers of elements of ${\sf I}_\eps$ and ${\sf J}_p$ is uniformly bounded. These sets are used in Lemma \ref{LemComparisonFormula} and Theorem \ref{thm:main} below. Finally, we prove the following lemma which is used in Section \ref{SubsectionRenormalizedModels} in the definition of the renormalized models and in the proofs of Lemma \ref{lem:Binega} and Lemma \ref{lem:goal}. Its proof is not essential on a first reading.

\ssk

\begin{lem}\label{*epsilonuseless}
If $\bf B$ is a differentiable sector, there exists a constant $\eps_0=\eps_0({\bf B})$ satisfying the following properties.
\begin{itemize}
\setlength{\itemsep}{0.1cm}
\item[(1)] $(V_0^+,W_{0,p}^+,\Delta_{0,p},\Delta_{0,p}^+)=(V_\eps^+,W_{\eps,p}^+,\Delta_{\eps,p},\Delta_{\eps,p}^+)$ for any $\eps\in[0,\eps_0)$ and $p\in\{2,\infty\}$.
\item[(2)] For any $\tau,\sigma\in{\bf B}\cup\dt{\bf B}$ and $p\in\{2,\infty\}$, the following condition holds:
\begin{center}
`{\sl If $r_{0,p}(\tau)>r_{0,p}(\sigma)$, then $r_{\eps,p}(\tau)>r_{\eps,p}(\sigma)$ for any $\eps\in[0,\eps_0)$}.'
\end{center}
\item[(3)] There exists a choice of distinct elements $\tau_1,\dots,\tau_n\in{\bf W}_{0,2}^+$ such that 
$$
{\sf I}_\eps = \big\{ p_\eps(\tau_1)<\dots<p_\eps(\tau_n) \big\}
$$ 
and their order does not change within $\eps\in[0,\eps_0)$.
\end{itemize}
\end{lem}

\ssk

\begin{Dem}
Note that we can write $r_{\eps,\infty}(\tau)=n_\tau(r_0-\eps)+m_\tau\beta_0+\sum_{j=1}^dl_\tau^j\mfs_j$ for some integers $n_\tau,m_\tau,l_\tau^j$.
\begin{enumerate}
	\item[(1)] By the assumption that $r_0\notin\bbQ[\beta_0,\mfs]$ and by the definition of ${\bf W}_{0,2}^+$, we have $r_{0,\infty}(\tau)\neq0$ and $r_{0,2}(\tau)>0$ for any $\tau\in{\bf W}_{0,2}^+$. Therefore, $[0,\eps_0)$ does not intersect with ${\sf J}_\infty$ nor ${\sf J}_2$ for small $\eps_0>0$.
	\item[(2)] It is obvious because $\eps\mapsto r_{\eps,p}(\tau)$ is continuous.
	
	\item[(3)] For each $\tau,\sigma\in{\bf W}_{0,2}^+$, if $p_0(\tau)=p_0(\sigma)$, then since $r_{0,\infty}(\tau)=r_{0,\infty}(\sigma)$, we have $n_\tau=n_\sigma$ by the assumption that $r_0\notin\bbQ[\beta_0,\mfs]$. Then $p_\eps(\tau)=p_\eps(\sigma)$ for any $\eps\ge0$. 
This implies that, in the above picture, two lines intersect at $\frac1p$-axis are the same.
Therefore, there exists $\tau_1,\dots,\tau_n\in{\bf W}_{0,2}^+$ such that ${\sf I}_\eps = \big\{p_\eps(\tau_1)<\dots<p_\eps(\tau_n)\big\}$ for any $\eps\in[0,\eps_0)$, where $\eps_0>0$ is the first time when two of the lines $\big\{\ell_{\tau_1},\dots,\ell_{\tau_n},\{p=2\},\{p=\infty\}\big\}$ intersect.
\end{enumerate}
This concludes the proof of this lemma.
\end{Dem}

%%-------------------------------------------------------------------------%%
\subsection{Regularity-integrability structures and their models$\boldmath .$   \hspace{0.03cm}}
\label{SubsectionRIS}
%%-------------------------------------------------------------------------%%

We introduce this notion to deal with models on a usual regularity structure and their first order Malliavin derivative in a single setting. See the companion work \cite{SemigroupMasato} for more detailed descriptions. 

\medskip

{\it \S1 Regularity-integrability structures.} We define a strict partial order on $\bbR\times[1,\infty]$ by setting
$$
(r,i)<(s,j) \quad \overset{\text{\rm def}}{\Longleftrightarrow} \quad \Big\{r<s \ \textrm{and} \ \frac{1}{i}\leq \frac{1}{j}\Big\},
$$
and start with a general definition. We sometimes denote by ${\bf a}$ a generic element of $\bbR\times[1,\infty]$.

\ssk

\begin{defn}
A \textbf{regularity-integrability structure} $(\strA,U,{\sf G})$ consists of the following elements.
\begin{enumerate} \setlength{\itemsep}{0.1cm}
\renewcommand{\labelenumi}{\textrm{\bf (\alph{enumi})}}
	\item[--] The index set $\strA$ is a subset of $\bbR\times[1,\infty]$ such that for every $(s,j)\in\bbR\times[1,\infty]$ the set $\{(r,i)\in \strA\, ;\, (r,i)<(s,j)\}$ is finite.

	\item[--] The vector space $U=\bigoplus_{{\bf a}\in \strA}U_{\bf a}$ is an algebraic sum of Banach spaces $(U_{\bf a},\|\cdot\|_{\bf a})$.  

	\item[--] The structure group $\sf G$ is a group of continuous linear operators on $U$ such that one has for all $\Gamma\in {\sf G}$ and ${\bf a}\in \strA$
$$
(\Gamma-\id)(U_{\bf a})\subset\bigoplus_{{\bf a}'\in \strA,\, {\bf a}'<{\bf a}} U_{{\bf a}'}.
$$
\end{enumerate}
\end{defn}

\ssk

For $(r,i)\in \textit{\strA}$ the number $r$ is interpreted as a regularity exponent and $i$ as an integrability exponent. One says that the regularity-integrability structure has {\it regularity $r_0\in\bbR$} if $(r_0,\infty) < {\bf a}$ for any ${\bf a}\in \textit{\strA}$. Denoting by $P_{\bf a} : T\to T_{\bf a}$ the canonical projection, we set with a slight abuse of notations
$$
\|\tau\|_{\bf a} \defeq \|P_{\bf a}\tau\|_{\bf a}
$$
for any $\tau\in T$ and ${\bf a}\in \textit{\strA}$.

\ssk

We now return to the setting of Section {{\ref{SubsectionDifferentialSectors}}} with some fixed parameters $\varepsilon$ and $p$. For each $\tau\in{\bf T}^{(0)}\cup{\bf T}^{(1)}$ we define the integrability exponent
\begin{equation}\label{Eq:ip(tau)}
i_p(\tau)=\begin{cases}
\infty&\textrm{if } \tau\in{\bf T}^{(0)},\\
p&\textrm{if } \tau\in{\bf T}^{(1)}.
\end{cases}
\end{equation}
Setting
\begin{equation*} \begin{split}
\textit{\strA}_{\varepsilon,p} &\defeq \Big\{\big(r_{\varepsilon,p}(\tau),i_p(\tau)\big)\, ;\, \tau\in{\bf B}\cup\dt{\bf B}\Big\}   \\
&\;= 
\Big\{\big(r_{\varepsilon,\infty}(\tau),\infty\big)\, ;\, \tau\in{\bf B}\Big\} \cup \Big\{\big(r_{\varepsilon,p}(\dt{\tau}),p\big)\, ;\, \dt{\tau}\in\dt{\bf B}\Big\}
\end{split} \end{equation*}
and 
$$
W_{\bf a} \defeq \spa\hspace{-0.06cm}\Big\{\tau \in {\bf B}\cup\dt{\bf B}\, ;\, \big(r_{\varepsilon,p}(\tau),i_p(\tau)\big)={\bf a}\Big\}
$$
we then have the following grading on $W=V\oplus\dt{V}$
\begin{align*}
W &=\bigoplus_{{\bf a}\in \textit{\strA}_{\varepsilon,p}} W_{\bf a}.
\end{align*}
We further introduce the group ${\sf G}_{\varepsilon,p}^+$ of characters on the Hopf algebra $(W_{\varepsilon,p}^+,\Delta_{\varepsilon,p}^+)$. The group ${\sf G}_{\varepsilon,p}^+$ has a representation in $GL(W)$ where ${\sf g}\in {\sf G}_{\varepsilon,p}^+$ is mapped to $(\id\otimes \mathop{\sf g})\Delta_{\varepsilon,p}$. Denote by ${\sf G}_{\varepsilon,p}$ the image group. Then the triple $(\textrm{A}_{\varepsilon,p},W, {\sf G}_{\varepsilon,p})$ is a regularity-integrability structure; it is said to be associated with the \textbf{\textit{concrete regularity-integrability structure}} 
$$
\scW_{\varepsilon,p} = \big((W , \Delta_{\varepsilon,p}), (W_{\varepsilon,p}^+ , \Delta^+_{\varepsilon,p})\big).
$$ 
A \textbf{\textit{sub-structure}} $(\textrm{A}'_{\varepsilon,p}, W', {\sf G}'_{\varepsilon,p})$ of $(\textrm{A}_{\varepsilon,p}, W, {\sf G}_{\varepsilon,p})$ is a regularity-integrability structure where $\textrm{A}'_{\varepsilon,p}\subset \textrm{A}_{\varepsilon,p}, W'\subset W$ and each element of ${\sf G}'_{\varepsilon,p}$ is a restriction of an element of ${\sf G}_{\varepsilon,p}$ to $W'$. We also have an analogous notion of concrete regularity-integrability sub-structure.

\medskip

{\it \S2 Models on regularity-integrability structures.} Fix $(\varepsilon,p)\in[0,\infty)\times[2,\infty]$ and $c>0$. Assume we are given a pair of maps ${\sf M}=({\sf \Pi}, {\sf g})$ such that
$$
{\sf\Pi} : W \to C^{r_0,Q}(w_c),   \qquad   {\sf g} : \bbR^d \to {\sf G}_{\varepsilon,p}^+
$$
and $\sf\Pi$ is continuous and linear.
The map $\sf \Pi$ is called an \textit{\textbf{interpretation map}}. For any $\tau\in W$ and $\mu\in W_{\varepsilon,p}^+$ set
\begin{equation*} \begin{split}
{\sf\Pi}^{\varepsilon,p}_x(\tau) &\defeq \big({\sf \Pi}\otimes ({\sf g}_x\circ S^+_{\varepsilon,p})\big) \Delta_{\varepsilon,p}(\tau),   \\
{\sf g}_{yx}^{\varepsilon,p}(\mu) &\defeq \big({\sf g}_{y}\otimes ({\sf g}_x\circ S^+_{\varepsilon,p})\big) \Delta^+_{\varepsilon,p}(\mu).
\end{split} \end{equation*}
We talk of the map ${\sf\Pi}^{\varepsilon,p}_x$ as the \textbf{\textit{recentered interpretation map}}. To define the space of models, recall from Section \ref{SectionFunctionalSetting} the definition of the operator semigroup $(\mcQ_t)_{t>0}$.

\ssk

\begin{defn}
A pair of maps ${\sf M}=({\sf \Pi}, {\sf g})$ as above is called a \textbf{\textit{model on $\scW_{\varepsilon,p}$ (with weight $w_c$)}} if
\begin{align} \label{EqDefnSizePiTau}
\|{\sf\Pi}^{\varepsilon,p} \hspace{-0.03cm}:\hspace{-0.03cm} \tau\|_{w_c}
&\defeq 
\sup_{0<t\le1} t^{-r_{\varepsilon,p}(\tau)/\ell} \big\| \mcQ_t(x, {\sf\Pi}^{\varepsilon,p}_x(\tau)) \big\|_{L_x^{i_p(\tau)}(w_c)}<\infty
\end{align}
for any $\tau\in{\bf B}\cup\dt{\bf B}$ and
\begin{align} \label{EqDefnSizeGMu}
\|{\sf g}^{\varepsilon,p} \hspace{-0.03cm}:\hspace{-0.03cm} \mu\|_{w_c}
&\defeq 
\|{\sf g}_x(\mu)\|_{L_x^{i_p(\mu)}(w_c)} + \sup_{y\in\bbR^d\setminus\{0\}}\Bigg(w_c(y)\frac{\big\|{\sf g}^{\varepsilon,p}_{(x+y)x}(\mu)\big\|_{L_x^{i_p(\mu)}(w_c)}}{\|y\|_{\fr s}^{r_{\varepsilon,p}(\mu)}}\Bigg)<\infty
\end{align}
for any $\mu\in{\bf W}_{\varepsilon,p}^+$. We denote by $\textbf{\textsf{M}}(\scW_{\varepsilon,p})_{w_c}$ the set of all models and set
$$
\|{\sf M}\|_{\textbf{\textsf{M}}(\scW_{\varepsilon,p})_{w_c}} \defeq \max_{\tau\in{\bf B}\cup\dt{\bf B}} \|{\sf\Pi}^{\varepsilon,p}\hspace{-0.03cm}:\hspace{-0.03cm} \tau\|_{w_c}  +  \max_{\mu\in{\bf W}_{\varepsilon,p}^+} \|{\sf g}^{\varepsilon,p} \hspace{-0.03cm}:\hspace{-0.03cm} \mu\|_{w_c}
$$
for any ${\sf M}\in\textbf{\textsf{M}}(\scW_{\varepsilon,p})_{w_c}$. We define a metric $\|{\sf M}_1 \hspace{-0.03cm}:\hspace{-0.03cm} {\sf M}_2\|_{\textbf{\textsf{M}}(\scW_{\varepsilon,p})_{w_c}}$ on $\textbf{\textsf{M}}(\scW_{\varepsilon,p})_{w_c}$ setting
$$
\|{\sf M}_1 \hspace{-0.03cm}:\hspace{-0.03cm} {\sf M}_2\|_{\textbf{\textsf{M}}(\scW_{\varepsilon,p})_{w_c}} \defeq \max_{\tau\in{\bf B}\cup\dt{\bf B}} \|{\sf\Pi}_1^{\varepsilon,p}, {\sf\Pi}_2^{\varepsilon,p} \hspace{-0.03cm}:\hspace{-0.03cm} \tau\|_{w_c}  +  \max_{\mu\in{\bf W}_{\varepsilon,p}^+} \|{\sf g}_1^{\varepsilon,p},{\sf g}_2^{\varepsilon,p} \hspace{-0.03cm}:\hspace{-0.03cm} \mu\|_{w_c}
$$
for all ${\sf M}_1, {\sf M}_2\in \textbf{\textsf{M}}(\scW_{\varepsilon,p})_{w_c}$. The quantities in the right hand side are defined in the same way as \eqref{EqDefnSizePiTau} and \eqref {EqDefnSizeGMu} but with $({\sf\Pi}_1)_x^{\varepsilon,p}\tau - ({\sf\Pi}_2)_x^{\varepsilon,p}\tau$, $({\sf g}_1)_{x}(\mu) - ({\sf g}_2)_{x}(\mu)$ and $({\sf g}_1)_{yx}^{\varepsilon,p}(\mu) - ({\sf g}_2)_{yx}^{\varepsilon,p}(\mu)$ in place of ${\sf\Pi}^{\varepsilon,p}_x\tau$, ${\sf g}_x(\mu)$ and ${\sf g}^{\varepsilon,p}_{yx}(\mu)$, respectively.
\end{defn}

\ssk

Since the exponents $\varepsilon$ and $p$ are involved not only in the definition of the recentered maps ${\sf\Pi}_x^{\varepsilon,p}$ and ${\sf g}_{yx}^{\varepsilon,p}$ but also in the definition of the norms \eqref{EqDefnSizePiTau} and \eqref{EqDefnSizeGMu} via the degree map $r_{\varepsilon,p}$, it would be more proper to write $\|(\cdot)^{\varepsilon,p} \hspace{-0.03cm}:\hspace{-0.03cm} \tau\|_{\varepsilon,p;w_c}$ for the norms. For the sake of readability, we use the above lightened notations.

\ssk

In the setting of regularity-integrability structures one can prove some analogues of the reconstruction theorem and the multilevel Schauder estimate. They are stated in Appendix {{\ref{SectionAppendix}}} in Theorem \ref{thm:besovreconstruction} and Theorem \ref{thm:MS}. One can find the detailed self-contained proofs of these statements in \cite{SemigroupMasato}.

\ssk

In the proof of Theorem \ref{thm:main} we will consider some restricted bounds on some concrete regularity-integrability sub-structures of $\scW_{\varepsilon,p}$ of the form
$$
\scW_{\varepsilon,p}'=\big((U,\Delta_{\varepsilon,p}),(U^+,\Delta_{\varepsilon,p}^+)\big)
$$
where $U=\spa({\bf C})$ for some subset ${\bf C}\subset{\bf B}\cup\dt{\bf B}$ and $U^+$ is a subalgebra of $W_{\varepsilon,p}^+$ generated by some subset ${\bf C}^+$ of ${\bf W}_{\varepsilon,p}^+$. For any ${\sf M}\in\textbf{\textsf{M}}(\scW_{\varepsilon,p})_{w_c}$ it is useful to define the restricted quantity
\begin{equation} \label{EqNormModelsSubstructure}
\begin{aligned}
\|{\sf M}\|_{\textbf{\textsf{M}}(\scW_{\varepsilon,p}')_{w_c}} \defeq \|{\sf\Pi}^{\varepsilon,p} \res {\bf C}\|_{w_c}  + \|{\sf g}^{\varepsilon,p} \res {\bf C}^+\|_{w_c} \defeq \max_{\tau\in\bf C} \|{\sf\Pi}^{\varepsilon,p} \res \tau\|_{w_c}
+\max_{\mu\in{\bf C}^+} \|{\sf g}^{\varepsilon,p} \res \mu\|_{w_c}.
\end{aligned}
\end{equation}
In particular the restriction of any ${\sf M}\in\textbf{\textsf{M}}(\scW_{\varepsilon,p})_{w_c}$ to the concrete regularity sub-structure 
$$
\scV_{\varepsilon} = \big((V,\Delta_{\varepsilon,p}),(V_{\varepsilon}^+,\Delta^+_{\varepsilon,p})\big)
$$ 
is a model in the usual sense of \cite{Hai14}.

\medskip

%%-----------------------------------%%
\subsection{Renormalized models$\boldmath .$   \hspace{0.03cm}}
\label{SubsectionRenormalizedModels}
%%-----------------------------------%%

We now introduce a whole class of models built in a recursive way.   

\medskip

{\it \S1 Admissible models.} Recall from \eqref{EqQSpacetimeOperator} the definition of the operator $\mcK$ acting on functions over $\bbR^d$. An interpretation map $\sf \Pi$ is said to be \textit{{$\mcK$-admissible}} if it satisfies
$$
({\sf\Pi}X^{k})(x)=x^{k},\qquad{\sf\Pi}(\mcI_{k}^\bbK\tau) = \partial^{k}\mcK(\cdot,{\sf\Pi}\tau)
$$
for all ${k}$ and $\tau$. 
Denote by $C_\star^\infty$ be the set of all smooth functions over $\bbR^d$ whose derivatives of all orders are in the class $\bigcup_{c>0}L^\infty(w_c)$. By the property {\bf (\ref{asmpleaf'})} of the differentiable sector, for any $\xi,h\in C_\star^\infty$ there is a unique {\it multiplicative} $\mcK$-admissible map ${\sf\Pi}^{\xi,h}:W\to C_\star^\infty$ such that
$$
{\sf\Pi}^{\xi,h}(\ocircle)=\xi,\qquad {\sf \Pi}^{\xi,h}(\odot)=h.
$$
Then, given any weight $w_c$ with $c>0$, we define a model ${\sf M}^{\xi,h;\varepsilon,p} = ({\sf\Pi}^{\xi,h}, {\sf g}^{\xi,h;\varepsilon,p})$ on $\scW_{\varepsilon,p}$ from the recursive definition of $({\sf g}_x^{\xi,h;\varepsilon,p})^{-1}={\sf g}_x^{\xi,h;\varepsilon,p}\circ S_{\epsilon,p}^+$ given by

\begin{equation*}   % \label{EqConstructionGfromPi}
\begin{aligned}
\big({\sf g}_x^{\xi,h;\varepsilon,p}\big)^{-1}(X^k) &= (-x)^k,   \\
\big({\sf g}_x^{\xi,h;\varepsilon,p}\big)^{-1}(\tau\sigma) &= ({\sf g}_x^{\xi,h;\varepsilon,p})^{-1}(\tau)({\sf g}_x^{\xi,h;\varepsilon,p})^{-1}(\sigma),   %\\
\end{aligned}
\end{equation*}

\begin{equation} \label{EqConstructionGfromPi}
\begin{aligned}
\big({\sf g}_x^{\xi,h;\varepsilon,p}\big)^{-1}(\mcI_k^{\fr l}(\sigma)) &= - \sum_{l\in\bbN^d}
\, \frac{(-x)^{l}}{{l}!} {\sf f}_x^{\xi,h;\varepsilon,p}(\mcI_{k+l}^{\fr l}(\sigma)),\\
{\sf f}_x^{\xi,h;\varepsilon,p}(\tau)
&={\bf1}_{r_{\varepsilon,p}(\tau)>0}\times
\left\{
\begin{aligned}
&\partial^kh(x)
&&(\tau=\mcI_k^\bbH(\bullet)),\\
&\partial^k\mcK\big(x, {\sf\Pi}_x^{\xi,h;\varepsilon,p}(\sigma)\big)
&&(\tau=\mcI_k^\bbK(\sigma)).
\end{aligned}
\right.
\end{aligned}
\end{equation}
See e.g. Section 3 of Bruned's work \cite{Bru18}.

\ssk

\begin{defn}
We call ${\sf M}^{\xi,h;\varepsilon,p} = ({\sf\Pi}^{\xi,h}, {\sf g}^{\xi,h;\varepsilon,p})$ the \textbf{\textit{naive model}} associated to $(\xi,h)$.
\end{defn}

\ssk

{\it \S2 Preparation maps and their associated models.} Here is how to build a large family of $\mcK$-admissible interpretation maps from a given interpretation map. Recall from \cite{Bru18} and Bruned \& Nadeem's work \cite{BNDiagramFree} the following definition. We write $|\tau|_\ocircle$ for the number of edges with $\bbOmega$-label, or nodes with $\ocircle$-label in the contracted form, for $\tau\in{\bf B}\cup\dt{\bf B}$.

\ssk

\begin{defn}
A \textbf{\textit{preparation map}} is a linear map 
$$
R : W\to W
$$ 
which leaves stable the subspace $V$ and has the following properties.
\begin{enumerate} \setlength{\itemsep}{0.1cm}
\renewcommand{\labelenumi}{\textrm{\bf (\alph{enumi})}}
\item\label{Prep1} The map $R$ fixes the polynomials and the noises: $R\tau=\tau$ for $\tau\in\{X^k,\ocircle,\odot\}$.
\item\label{PrepTriangle} For each $ \tau \in {\bf B}\cup\dt{\bf B}$ there exist finitely many $\tau_i \in {\bf B}\cup\dt{\bf B}$ and constants $\lambda_i$ such that
$$
R \tau = \tau + \sum_i \lambda_i \tau_i, \quad\textrm{with}\quad r_{0,p}(\tau_i) > r_{0,p}(\tau) \ \text{for}\ p\in\{2,\infty\}\quad\textrm{and}\quad |\tau_i|_{\ocircle} < |\tau|_{\ocircle}.
$$
\item\label{EqRandI}
The map $R$ fixes planted trees: $R\,\mcI_{k}^\bbK = \mcI_{k}^\bbK$.
\item\label{EqCommutationRDelta}
The map $R$ commutes with the coproduct: $(R \otimes \id) \Delta_{0,2} = \Delta_{0,2} R$.
\item\label{Prep5}
The map $R$ commutes with the derivative map \eqref{EqDefnDerivation}: $RD=DR$.
\end{enumerate}
\end{defn}

\ssk

In the triangular property {\bf (\ref{PrepTriangle})}, the strict inequaltiy $r_{\eps,p}(\tau_i)>r_{\eps,p}(\tau)$ holds for any $(\eps,p)\in[0,\eps_0)\times[2,\infty]$. First, by Lemma \ref{*epsilonuseless}, the inequality $r_{\eps,p}(\tau_i)>r_{\eps,p}(\tau)$ holds for any $(\eps,p)\in[0,\eps_0)\times\{2,\infty\}$. Since $r_{\eps,p}(\cdot)$ is affine for $\frac1p$, if one has $r_{\eps,p}(\tau_i)>r_{\eps,p}(\tau)$ for $p\times\{2,\infty\}$ one also has this inequality for $p\in(2,\infty)$. We note that the commutation property {\bf (\ref{EqCommutationRDelta})} also hold for any exponents $(\varepsilon,p)\in[0,\varepsilon_0)\times[2,\infty]$. 
Indeed, since $\Delta_{\varepsilon,p}=(\id\otimes P_{\varepsilon,p}^+)\Delta_{0,2}$ by definition, we have
\begin{align*}
(R\otimes\id)\Delta_{\varepsilon,p}&=(\id\otimes P_{\varepsilon,p}^+)(R\otimes\id)\Delta_{0,2}
=(\id\otimes P_{\varepsilon,p}^+)\Delta_{0,2}R=\Delta_{\eps,p}R.
\end{align*}

A typical example of preparation map $R$ is defined from the choice of some constants $\{c(\sigma)\}_{\sigma\in{\bf B},\, r_{0,\infty}(\sigma)<0}$ by the extraction-contraction formula
\begin{align}\label{EqRdual}
(R_c-\id)\tau=(c\otimes\id)\Delta_{0,2}\tau\qquad(\tau\in {\bf B}\cup\dt{\bf B})
\end{align}
as in Corollary 4.5 of \cite{Bru18} -- see Proposition \ref{*prop:ExampleDifferentiable} below for the precise definition and the proof of the fact that it is indeed a preparation map.

\ssk

For any given preparation map $R$ we can construct the renormalized model basically following \cite{Bru18}. First, we can define the unique linear map
$$
\widehat{M}^R : W\to W
$$ 
by requiring that $\widehat{M}^R$ is {\it multiplicative}, fixes $X^k$, $\ocircle$, and $\odot$, and satisfies
$$
\widehat{M}^R(\mcI_{k}^\bbK\tau) = \mcI_{k}^\bbK\big(\widehat{M}^R (R\tau)\big)
$$
for all $k$ and $\tau$. Then the linear map
$$
M^{\!R} \defeq \widehat{M}^R R.
$$
is the \textbf{\textit{renormalization map associated to the preparation map $R$}}. It follows from the property {\bf(\ref{EqRandI})} of preparation maps that if $\sf \Pi$ is a $\mcK$-admissible interpretation map then so is the map 
$$
{\sf \Pi}^{\!R} \defeq {\sf \Pi}M^{\!R}.
$$
For the naive model ${\sf M}^{\xi,h;\varepsilon,p}$ associated to $\xi,h\in C_\star^\infty$, we can define the model
$$
{\sf M}^{\xi,h,R;\varepsilon,p} \defeq ({\sf \Pi}^{\xi,h,R}, {\sf g}^{\xi,h,R;\varepsilon,p})
$$
by ${\sf\Pi}^{\xi,h,R}=({\sf\Pi}^{\xi,h})^R$ and ${\sf g}^{\xi,h,R;\varepsilon,p}$ defined by the recursive rules \eqref{EqConstructionGfromPi} where $({\sf \Pi}^{\xi,h}, {\sf g}^{\xi,h;\varepsilon,p})$ is replaced with $({\sf \Pi}^{\xi,h,R}, {\sf g}^{\xi,h,R;\varepsilon,p})$.

\ssk

\begin{defn}
We call ${\sf M}^{\xi,h,R;\varepsilon,p}$ the \textbf{\textit{renormalized model}} of ${\sf M}^{\xi,h;\varepsilon,p}$ associated to $R$.
\end{defn}

\ssk

Recall that $R$ leaves stable the subspace $V$. Since the above discussions still work if $W$ is replaced by $V$, the model ${\sf M}^{\xi,h,R;\varepsilon,p}$ restricted on $V$ is independent of $h$ and $p$. Thus we are allowed to write
$$
{\sf \Pi}^{\xi,h,R}\vert_V={\sf\Pi}^{\xi,R},\qquad
{\sf g}^{\xi,h,R;\varepsilon,p}\vert_{V_\varepsilon^+}={\sf g}^{\xi,R;\varepsilon}.
$$

In the setting of singular stochastic PDEs we insist on renormalizing models with preparation maps that are {\it interpretable}. Denote by ${\bf u}_R$ the solution in a space of modelled distributions to the regularity structure formulation of an equation $Lu=f(u,\xi)$ associated with a model built from a preparation map $R$ and a smooth noise $\xi$. We say that the preparation map $R$ is interpretable if the reconstruction $u_R$ of ${\bf u}_R$ is the solution to an equation of the form $Lu_R = f(u_R,\xi) + c(u_R)$, for some counterterms $c(u_R)$ that may depend on $u_R$ and some of its derivatives. It was proved in Proposition 3.2 of \cite{BB21} that the preparation maps $R$ of the form \eqref{EqRdual} are indeed interpretable. 
More precisely, this map is a `strong preparation map' in the sense of \cite{BB21}. This stronger property is needed to derive the renormalized equation in the form stated in Corollary \ref{CorMainThm} -- see Theorem 4.1 of \cite{BB21}. The BPHZ renormalisation map from \cite{BHZ} corresponds to a particular choice of coefficients $c(\sigma)$ in the formula \eqref{EqRdual}. Our proof of Theorem \ref{thm:main} works for {\it any} preparation map, strong or not.

\medskip

{\it \S3 An important property of the renormalized models.} Before going to the main result we introduce an important algebraic identity satisfied by the renormalized models ${\sf M}^{\xi,h,R;\varepsilon,p}$. The proofs of this results is given in Section \ref{ProofLemmaComparisonFormula}. 

\ssk

We introduce a comparison formula between ${\sf\Pi}_x^{\xi,h,R;\varepsilon,p}$ and ${\sf\Pi}_x^{\xi,h,R;\varepsilon,2}$. Recall the notations in Lemma \ref{*epsilonuseless}. Since the set ${\sf I}_\eps$ is finite, so its associated `floor function' is defined as
$$
\lfloor p\rfloor_{{\sf I}_\varepsilon} \defeq \max\big\{q\in \{2\}\cup {\sf I}_\varepsilon\, ;\, q<p\big\}
$$
for any $p\in(2,\infty]$. The structure of $(W_{\eps,p}^+,\Delta_{\eps,p}^+)$ is stable within the intervals $\big[\lfloor p_\eps(\mu)\rfloor_{{\sf I}_\eps},p_\eps(\mu)\big)$.
Any tree of the form $\mcI_k^{\fr l}(\sigma)$ is called {\it planted}. 
For given inputs $\xi,h,R$, we define the linear map $\lambda_x^{\xi,h,R;\varepsilon,p}:W_{\varepsilon,2}^+\to\bbR$ by
\begin{align} \label{EqDefnCaracterH}
\lambda_x^{\xi,h,R;\varepsilon,p}(\mu)\defeq
\left\{
\begin{aligned}
&{\bf1}_{r_{\varepsilon,p}(\mu)\le0<r_{\varepsilon,2}(\mu)}
\, {\sf f}_x^{\xi,h,R;\varepsilon,\lfloor p_\varepsilon(\mu)\rfloor_{{\sf I}_\varepsilon}}(\mu)
&&(\text{$\mu$ is a planted tree}),\\
&0
&&(\text{$\mu$ is a non-planted tree}),
\end{aligned}
\right.
\end{align}
where ${\sf f}_x^{\xi,h,R;\varepsilon,p}(\mu)$ is defined by \eqref{EqConstructionGfromPi} with ${\sf\Pi}_x^{\xi,h,R;\varepsilon,p}$ in the role of ${\sf\Pi}_x^{\xi,h;\varepsilon,p}$. We have in particular $\lambda_x^{\xi,h,R;\varepsilon,p}(\mu)=0$ if $\mu\in{\bf T}^{(0)}$, since then $r_{\varepsilon,p}(\mu)=r_{\varepsilon,2}(\mu)$.

\ssk

\begin{lem} \label{LemComparisonFormula}
Let $\xi,h\in C_\star^\infty$ and let $R$ be a preparation map. For any $\tau\in\dt{\bf B}$ and any $(\varepsilon,p)\in[0,\infty)\times[2,\infty]$, one has
\begin{equation} \label{EqSyntheticDecompositionFormulaDerivativePi} \begin{split}
{\sf\Pi}_x^{\xi,h,R;\varepsilon,p}(\tau)
&= {\sf\Pi}_x^{\xi,h,R;\varepsilon,2}(\tau) + \big({\sf\Pi}^{\xi,h,R;\varepsilon,p}_x\otimes\lambda_x^{\xi,h,R;\varepsilon,p}\big)\Delta_{\varepsilon,2}(\tau).
\end{split} \end{equation}
\end{lem}

\ssk

Bruned \& Nadeem proved in Proposition 3.7 of \cite{BNDiagramFree} a statement with a similar flavour. The above formula plays a crucial role in the proof of Theorem \ref{thm:main}, especially in the proof of Lemma {{\ref{lem:goal}}}. Roughly speaking, ${\sf\Pi}_x^{\xi,h,R;\varepsilon,2}(\tau)$ will be seen to have an interesting property for some elementary reason. The map $\Delta_{\varepsilon,2}$ splitting $\tau$ into `smaller' pieces, the identity \eqref{EqSyntheticDecompositionFormulaDerivativePi} will allow us to implement an inductive procedure to see that ${\sf\Pi}_x^{\xi,h,R;\varepsilon,p}(\tau)$ also has the property of interest. In the example \eqref{ExL6-} with $\xi,h\in C_\star^\infty$, and with the identity as a preparation map giving the naive model, one has
\begin{align*}
{\sf\Pi}_x^{\xi,h;\varepsilon,p}
\big(\,\rtdb{3}{\pol{2}{1}{90}\drl{1}{2}\drc{1}\drch{2}}\,\big)
(\cdot)
= \left\{\begin{aligned}
&\big(\mcK(\cdot,h) - \mcK(x,h)\big)\xi(\cdot), & \quad \textrm{ for } p\geq\frac{6}{1+2\varepsilon},   \\
&\Big(\mcK(\cdot,h) - \mcK(x,h) - \sum_{i=1}^3\partial_i\mcK(x,h)\,(\cdot-x)_i\Big)\xi(\cdot), &  \quad \textrm{ for } p<\frac{6}{1+2\varepsilon}.
\end{aligned} \right.
\end{align*}
(We use the contracted from $\, \rtdbs{5}{1}{\pol{2}{1}{90}\pol{3}{1}{30}\pol{4}{2}{30}
\drl{1}{2}\drdt{1}{3}\drdh{2}{4}\drb{1,2,3,4}}\, =\,\rtdb{3}{\pol{2}{1}{90}\drl{1}{2}\drc{1}\drch{2}}\,$.)
The formula \eqref{EqSyntheticDecompositionFormulaDerivativePi} writes here
\begin{align*}
{\sf\Pi}_x^{\xi,h;\varepsilon,p}\big(\,\rtdb{3}{\pol{2}{1}{90}\drl{1}{2}\drc{1}\drch{2}}\,\big)(\cdot) 
&=  
{\sf\Pi}_x^{\xi,h;\varepsilon,2}\big(\,\rtdb{3}{\pol{2}{1}{90}\drl{1}{2}\drc{1}\drch{2}}\,\big)(\cdot) 
+\sum_{i=1}^3 \partial_i\mcK(x,h)\,(\cdot-x)_i \, \xi(\cdot)
\\
&=
{\sf\Pi}_x^{\xi,h;\varepsilon,2}\big(\,\rtdb{3}{\pol{2}{1}{90}\drl{1}{2}\drc{1}\drch{2}}\,\big)(\cdot) 
+\sum_{i=1}^3
\lambda_x^{\xi,h,R;\varepsilon,p}\big(\,\rtdb{3}{\pol{2}{1}{90}\drl{1}{2}\drb{1}\drch{2}\node at ($0.5*(A1)+0.5*(A2)+(0.15,0)$) {\tiny $e_i$};}\!\big)
{\sf\Pi}_x^{\xi,h;\varepsilon,p}(\ocircle^{e_i})(\cdot)
\end{align*}
for $p\geq\frac{6}{1+2\varepsilon}$.

\medskip

\noindent \textbf{Remark.} {\sl For the readers familiar with the work \cite{LOTT} of Linares, Otto, Tempelmayr \& Tsatsoulis, one can give a `dictionary' between the different terms that appear in \eqref{EqSyntheticDecompositionFormulaDerivativePi} and some quantities that play a fundamental role in \cite{LOTT}. The term ${\sf\Pi}_x^{\xi,h,R;\varepsilon,p}\tau$ is our equivalent of their term $\delta{\sf \Pi}_{x\beta}$, the term $({\sf\Pi}_x^{\xi,h,R;\varepsilon,2}\tau)(y)$ is the equivalent of $\big(\delta{\sf \Pi}_{x\beta} - (d{\sf \Gamma}^*_{xy}){\sf \Pi}_{y\beta}\big)(y)$, and $\big(\big({\sf\Pi}^{\xi,h,R;\varepsilon,p}_x\otimes\lambda_x^{\xi,h,R;\varepsilon,p}\big)\Delta_{\varepsilon,2}\tau\big)(y)$ is our equivalent of their term $(d{\sf \Gamma}^*_{xy}){\sf \Pi}_{y\beta}(y)$.}

\medskip

%%------------------------------------%%
\subsection{Spectral gap inequality$\boldmath .$   \hspace{0.03cm}}
\label{SubsectionSpectralGap}
%%------------------------------------%%

The last ingredient we need before we can introduce the main convergence result of this work is the spectral gap inequality. It is an infinite dimensional analogue of the Poincar\'e inequality that can be satisfied or not by a probability measure on a Banach space.

\ssk

Let $\Omega$ be a separable Banach space and let $H$ be a separable Hilbert space embedded continuously and densely into $\Omega$. A function $F:\Omega\to\bbR$ is said to be \textit{\textbf{(continuously) $H$-differentiable}} if there is a function $dF: \omega\in\Omega\mapsto d_\omega F\in H^*$ such that
$$
\frac{d}{dt}F(\omega+th)\big|_{t=0} = (d_\omega F)(h).
$$
Denote by $\|h^*\|_{H^*}=\sup_{\Vert h\Vert_H\leq 1}\vert h^*(h)\vert$ the operator norm on $H^*$. A Borel probability measure $\bbP$ on $\Omega$ is said to satisfy the \textbf{\textit{$H$-spectral gap inequality}} if there exists a constant $C>0$ such that
\begin{equation} \label{EqSGInequality}
\bbE\big[(F-\bbE[F])^2\big] \leq C\,\bbE\big[\|dF\|_{H^*}^2\big]
=C\,\bbE\bigg[\sup_{\Vert h\Vert_H\leq 1}|dF(h)|^2\bigg]
\end{equation}
for any $H$-differentiable $F\in L^2(\Omega)$ such that $dF\in L^2(\Omega;H^*)$. Note that the supremum over $h$ is indeed a random variable since $(d_\omega F)(h)$ is continuous in $h\in H$ and the supremum over $h$ can be replaced by a countable supremum. By replacing $F$ by $F^2$ and iterating, one further has
\begin{equation}\label{EqSGInequalityModified}
\bbE\big[F^{2^r}\big] \lesssim_r \vert\bbE[F]\vert^{2^r} + \bbE\big[\|dF\|_{H^*}^{2^r}\big]
\end{equation}
for all $r\in\bbN$ -- see Remark 2.21 of Hairer \& Steele's work \cite{HS23}. This version of the inequality plays a role later in Section \ref{sec:ind1}. The $H$-spectral gap inequality holds if $(\Omega,H,\bbP)$ is an abstract Wiener space -- see e.g. Exercise 2.11.1 of Nourdin \& Peccati's textbook \cite{NP}. 

\ssk

A typical example is the white noise measure on $\bbR^d$, which satisfies the $H$-spectral inequality with $H=L^2(w_0)=L^2(\bbR^d)$. In this work we consider the Hilbert space
$$
H \defeq H^{-s_0,Q}(w_0)
$$
for some $s_0\in\bbR$. For example, the noise $(1-\Delta)^{s_0/2}(\xi)$ defined from a space white noise $\xi$ satisfies the $H^{-s_0,Q}(w_0)$-spectral inequality. By the Besov embedding \eqref{BesovEmbedding} the space $H$ is continuously embedded into the Banach space
$$
\Omega \defeq C^{r_0,Q}(w_c)
$$
for any $c>0$ and $r_0<-\frac{|\mfs|}2-s_0$. Since the renormalization problems do not arise when $r_0>0$, we assume $\frac{|\mfs|}2+s_0\geq 0$ in this work. We have $H^*=H^{s_0,Q}(w_0)$ in our setting, so the bigger $s_0>0$ the weaker the inequality \eqref{EqSGInequality}.

\ssk

We finish this section with an elementary, but fundamental, result about renormalized models.

\ssk

\begin{lem}\label{Lem:dPi=PiD}
Let $\xi,h\in C_\star^\infty$ and let $R$ be a preparation map. For any $\tau\in{\bf B}$ one has
\begin{equation}\label{Eq:dPi=PiD}
d_\xi({\sf\Pi}_x^{\xi,R;\varepsilon}\tau)(h)\defeq\frac{d}{dt}({\sf\Pi}_x^{\xi+th,R;\varepsilon}\tau)\big|_{t=0}={\sf\Pi}_x^{\xi,h,R;\varepsilon,\infty}(D\tau).
\end{equation}
\end{lem}

\ssk 

Recall that we can omit the letters $h$ and $p$ from ${\sf\Pi}_x^{\xi,h,R;\varepsilon,p}\tau$ since $\tau\in{\bf B}$. The above fact is also proved in Proposition 4.1 of \cite{BNDiagramFree}. We remark that the integrability exponent in the right hand side of \eqref{Eq:dPi=PiD} is $\infty$, rather than any finite $p$. The identity \eqref{Eq:dPi=PiD} is wrong if $\infty$ is replaced by any finite $p$.

\ssk

%----------------------------------------%
\subsection{A convergence result$\boldmath .$   \hspace{0.03cm}}
\label{SubsectionMainResult} 
%----------------------------------------%

We now have all the ingredients to state the formal version of Theorem \ref{ThmMainIntro} under the form of Theorem \ref{thm:main} below.

\ssk

For a family of compactly supported smooth functions $\varrho_n\in C^\infty(\bbR^d)$ that converge weakly to a Dirac mass at $0$ as $n\in\bbN$ goes to $\infty$, denote by $*$ the convolution operator and define a random variable $\xi_n$ on $\Omega$ setting
$$
\xi_n(\omega) \defeq \varrho_n*\omega \in\Omega.
$$
Set as well
$$
h_n  \defeq \varrho_n*h
$$
for $h\in H$. For a given preparation map $R_n$, we define a random variable ${\sf M}^{\xi_n,h_n,R_n;\varepsilon,p}$ setting 
$$
{\sf M}^{\xi_n,h_n,R_n;\varepsilon,p}(\omega)\defeq{\sf M}^{\xi_n(\omega),h_n,R_n;\varepsilon,p}.
$$ 
The following result is a direct consequence of Proposition 3.16 of \cite{Bru18}.

\ssk

\begin{lem} \label{LmModelTargetSpace}
The random variable ${\sf M}^{\xi_n,h_n,R_n;\varepsilon,p}$ takes almost surely its values in $\textbf{\textsf{M}}(\scW_{\varepsilon,p})_{w_c}$.
\end{lem}

\ssk

The choice of a sequence $(R_n)_{n\geq 0}$ of preparation maps is called a {\sl renormalization procedure}. In \eqref{EqControlExpectationQ1} below, for $\tau\in\bf{B}$, we write ${\sf\Pi}_0^{\xi_n,R_n;0}(\tau)$ rather than ${\sf\Pi}_0^{\xi_n,h_n,R_n;0,p}(\tau)$ as this function does not depend on $h_n$ and $p$ in that case. We can now state our main result.

\ssk

\begin{thm} \label{thm:main}
Pick $s_0\geq-\frac{|\mfs|}2$ and $r_0<-\frac{|\mfs|}2-s_0$ and set $H = H^{-s_0,Q}(w_0)$. Let $\bbP$ be a Borel probability measure on $\Omega=C^{r_0,Q}(w_c)$, for some $c>0$, that is centered, translation-invariant and satisfies the $H$-spectral gap inequality. Let ${\bf B}$ be a differentiable sector for which all $\tau\in{\bf B}\setminus\{\ocircle\}$ satisfy
$$
r_{0,\infty}(\tau) + \frac{\vert{\fr s}\vert}{2} > 0.
$$
Assume we are given a family $(R_n)_{n\geq 0}$ of preparation maps on $W$ for which all the quantities 
\begin{equation} \label{EqControlExpectationQ1}
\bbE\big[\mcQ_1\big(0, {\sf\Pi}_0^{\xi_n,R_n;0}(\tau)\big)\big]
\end{equation}
converge as $n$ goes to $\infty$, for any $\tau\in{\bf B}$ with $r_{0,\infty}(\tau)\le0$. 
\begin{itemize}
	\item[--] There exist a large constant $c_0>0$ depending only on $\mfs$ and ${\bf B}$ such that for any $c>c_0$, $p\in[2,\infty]$, $\varepsilon\in(0,\varepsilon_0)\setminus {\sf J}_p$ and $q\in[1,\infty)$ one has
\begin{equation} \label{EqExpectationSupH}
\sup_{n\in\bbN}\bbE\Bigg[\sup_{\Vert h\Vert_H\leq 1} \Vert {\sf M}^{\xi_n,h_n,R_n;\varepsilon,p} \Vert_{\textbf{\textsf{M}}(\scW_{\varepsilon,p})_{w_c}}^q\Bigg]<\infty.
\end{equation}

	\item[--] The models ${\sf M}^{\xi_n,h_n,R_n;\varepsilon,p}$ converge in $L^q\big(\Omega,\bbP ; \textbf{\textsf{M}}(\scW_{\varepsilon,p})_{w_c}\big)$ as $n$ goes to $\infty$.   \vspace{0.1cm}

	\item[--] The limit model depends on the renormalization procedure $(R_n)_{n\geq 0}$ only via the finite collection of numbers $\big\{\lim_{n\to\infty} \bbE\big[\mcQ_1(0, {\sf\Pi}_0^{\xi_n,R_n;0}\tau)\big] \big\}_{\tau\in{\bf B},\,r_{0,\infty}(\tau)\le0}$. Further, it does not depend on which approximation of the identity $(\varrho_n)_{n\geq 1}$ is used.
\end{itemize}
\end{thm}

\ssk

As a direct consequence of Theorem \ref{thm:main}, the ($(h_n,p)$-independent) restrictions to the concrete regularity structure $\scV_{\varepsilon} = \big((V,\Delta_{\varepsilon,p}),(V_\varepsilon^+,\Delta^+_{\varepsilon,p})\big)$ of the models ${\sf M}^{\xi_n, h_n,R_n;\varepsilon,p}$ converge in $L^q\big(\Omega,\bbP ; \textbf{\textsf{M}}(\scV_{\varepsilon})_{w_c}\big)$ as $n$ goes to $\infty$ for any $q\in[1,\infty)$ and $c>0$.

\ssk

Proceeding as in the proof of Theorem 6.18 of \cite{BHZ} one can see that there is a unique preparation map $\overline{R}_n=\overline{R}_n(\bbP)$ of the form \eqref{EqRdual} such that the associated model satisfies
$$
\bbE\big[\mcQ_1(0,{\sf\Pi}_0^{\xi_n,R_n;0,\infty}\tau)\big] = 0
$$
for any $\tau\in{\bf B}$ of non-positive $r_{0,\infty}$-degree. This modelled is called the $\overline{\textrm{BPHZ}}$  model in \cite{HS23}. This model obviously satisfies the assumption of our result. There is also a unique preparation map $R_n=R_n(\bbP)$ of the form \eqref{EqRdual} such that the associated model satisfies
\begin{equation} \label{EqdeterminesBHZ}
\bbE\big[({\sf\Pi}^{\xi_n,R_n}\tau)(x)\big]=0
\end{equation}
for all $\tau\in{\bf B}$ of non-positive $r_{0,\infty}$-degree. This model is called the BPHZ model. Although the argument is not straightforward, it is possible to show that the BPHZ model also satisfies the same condition. See Section \ref{sec:ind6} for a proof.

\ssk

\begin{cor}\label{cor:BPHZ}
The BPHZ model converges in $L^q\big(\Omega,\bbP ; \textbf{\textsf{M}}(\scW_{\varepsilon,p})_{w_c}\big)$ as $n$ goes to $\infty$ under the same choices of $c,p,\eps,q$ as in Theorem \ref{thm:main}. The limit model does not depend on the choice of the approximation of the identity $(\varrho_n)_{n\geq 1}$.
\end{cor}

\ssk

The convergence results of \cite{HS23} are restricted to the BPHZ and $\overline{\textrm{BPHZ}}$ models. Theorem \ref{thm:main} is thus an extension of the main result of \cite{HS23} to a (much) larger class of renormalization procedures. 

\ssk

We close this section with two remarks on some direct extensions of Theorem \ref{thm:main} to some different settings.
We expect that analogous extensions could also be obtained by adapting
the method of \cite{HS23}, although they are not stated explicitly there.

\begin{itemize}
	\item[--] If the noise $\xi$ does not necessarily depend on all variables $x_1,\dots,x_d$, without loss of generality, we assume that $\xi$ depends on only the last $e$ variables $x_\circ=(x_{d-e+1},\dots,x_d)$. If the polynomials $X_1,\dots,X_{d-e}$ associated with the first $d-e$ variables $x_\bullet=(x_1,\dots,x_{d-e})$ do not appear in the regularity structure under consideration, then the model ${\sf M}^{\xi_n,h_n,R_n;\varepsilon,p}$ does not depend on $x_\bullet$ and can be written in terms of the partially integrated kernel
$$
K_t^\circ(x_\circ)=\int_{\bbR^{d-e}}K_t(x_\bullet,x_\circ)dx_\bullet
$$
instead of $K_t(x)$. Therefore by replacing $Q_t(x)$ and $\mfs$ with
$$
Q_t^\circ(x_\circ)=\int_{\bbR^{d-e}}Q_t(x_\bullet,x_\circ)dx_\bullet
\quad\text{and}\quad
\mfs_\circ=(\mfs_{d-e+1},\dots,x_d),
$$
we can obtain the same result as Theorem \ref{thm:main}. A typical example is the parabolic Anderson model
$$
(\partial_1-\Delta)u=f(u)\xi
$$
with the spatial white noise $\xi(x_2,\dots,x_d)$ for $d=2,3$.   \vspace{0.15cm}

	\item[--] Assume that instead of the elliptic operator \eqref{sec1:ellipticop} with constant coefficients we consider an operator
$$
P(x,\partial_x)=\sum_{|k|_\mfs\le\ell}a_k(x)\partial_x^k
$$
with bounded and H\"older continuous coefficients which satisfies the ellipticity \eqref{sec1:ellipticity} uniformly over $x$. The fundamental solution $Q_t(x,y)$ associated with $P(x,\partial_x)$ is not of the form $Q_t(x-y)$ in that case, but the `translation invariance' of the kernel $Q_t(x-y)$ is only used in the proof of Theorem \ref{thm:main} via Lemma \ref{lem:stationary} below. This implies that, even if we employ the non-translation invariant kernel $Q_t(x,y)$, we can show a weaker statement than Theorem \ref{thm:main}: For any $\tau\in{\bf B}$ with $r_{0,\infty}(\tau)\le0$, instead of the convergence of the quantities \eqref{EqControlExpectationQ1}, assume that the uniform bound
\begin{align*}
\sup_{n\in\bbN}\big\|\bbE\big[\mcQ_t(x,{\sf\Pi}^{\xi_n,R_n;0,\infty}_x\tau)\big]\big\|_{L_x^\infty(w_c)}\lesssim t^{\frac{r_{0,\infty}(\tau)}\ell}
\end{align*}
and the convergence
\begin{align*}
\lim_{m,n\to\infty}\sup_{0<t\le1}t^{-\frac{r_{0,\infty}(\tau)}\ell}
\big\|\bbE\big[\mcQ_t(x,{\sf\Pi}^{\xi_m,R_m;0,\infty}_x\tau-{\sf\Pi}^{\xi_n,R_n;0,\infty}_x\tau)\big]\big\|_{L_x^\infty(w_c)}=0
\end{align*}
hold for any $c>0$.
Then the same convergence result of the model ${\sf M}^{\xi_n,h_n,R_n;\varepsilon,p}$ as in Theorem \ref{thm:main} holds. This `conditional' convergence result significantly reduces the effort required for a potential proof of the BPHZ theorem in a non-translation invariant setting, compared to the arguments in \cite{Hai14, CH}. If further $\xi$ is a centered Gaussian then any polynomials of $\xi$ with odd degrees have vanishing expectations.
\end{itemize}

\vfill \pagebreak

%-------------------------------------------------%
\section{The mechanics of convergence}
\label{SectionMechanics}
%------------------------------------------------------%

%%-----------------------------------%%
\subsection{Structure of the proof$\boldmath .$   \hspace{0.03cm}}
\label{SubsectionStructureProof}
%%-----------------------------------%%

Our proof of Theorem \ref{thm:main} proceeds by a finite induction. We introduce an ad hoc structure for this induction in {\it \S1}. We describe in {\it \S2} the flow of the proof.

\bigskip

{\it \S1 The inductive structure.} We define a preorder $\preceq$ on ${\bf T}$ by setting
\begin{equation} \label{EqDefnOrder}
\sigma \preceq \tau \quad \overset{\text{def}}{\Longleftrightarrow} \quad \big(|\sigma|_\ocircle, |E_\sigma|,r_{0,\infty}(\sigma)\big) \leq_{\textrm{lex}} \big(|\tau|_\ocircle,|E_\tau|,r_{0,\infty}(\tau)\big)
\end{equation}
with the inequality $\leq_{\textrm{lex}}$ in the right hand side standing here for the lexicographical order. We also define $\sigma\prec\tau$ if $\sigma\preceq\tau$ and $\sigma\neq\tau$. The following structural property is an immediate consequence of the formula \eqref{EqDefnDeltap} defining $\Delta_{\varepsilon,p}$.

\ssk

\begin{lem} \label{leminduction*}
For each $\tau\in{\bf B}\cup\dt{\bf B}$ we set 
\begin{equation*} \begin{split} 
&{\bf C}_{\prec\tau}\defeq\{\sigma\in{\bf B}\cup\dt{\bf B}\,;\,\sigma\prec\tau\},   \\
&C_{\prec\tau}\defeq\spa({\bf C}_{\prec\tau}).
\end{split} \end{equation*}
Define as well the sub-algebra $C_{\prec\tau,\varepsilon,p}^+$ of $W_{\varepsilon,p}^+$ generated by the symbols 
$$
{\bf C}_{\prec\tau,\varepsilon,p}^+\defeq\{X_j\}_{j=1}^d\cup\{\odot_k\}_{|k|_\mfs<r_{\varepsilon,p}(\bbH)}\cup\{\mcI_k^{\bbK}(\sigma)\}_{k\in\bbN^d,\, \sigma\in({\bf B}\cup\dt{\bf B})\setminus\{X^l\}_{l\in\bbN^d},\,\sigma\prec\tau,\, r_{\varepsilon,p}(\sigma)+\beta_0>|k|_\mfs}.
$$
Then for each $\tau\in{\bf B}\cup\dt{\bf B}$ we have
$$
\big(\Delta_{\varepsilon,p}(\tau)-\tau\otimes{\bf1}\big) \in  C_{\prec\tau}\otimes C_{\prec\tau,\varepsilon,p}^+.
$$
\end{lem}

\ssk

\begin{Dem}
For any $(\tau,{\fr t},{\fr n},{\fr e})\in{\bf B}\cup\dt{\bf B}$, the explicit representation of $\Delta_{\varepsilon,p}(\tau,{\fr t},{\fr n},{\fr e})$ is given by (the truncation of) \eqref{*eq:graphicalcoprod}.
For simplicity, we write $\bar{\sigma}=(\sigma,{\fr t}\vert_{E_\sigma},{\fr n}_\sigma+\pi{\fr e}_{\partial\sigma},{\fr e}\vert_{E_\sigma})$ for the left part and $\bar{\eta}=(\tau/\sigma,{\fr t}\vert_{E_\tau\setminus E_\sigma},[{\fr n}-{\fr n}_\sigma]_\sigma,{\fr e}\vert_{E_\tau\setminus E_\sigma}+{\fr e}_{\partial\sigma})$ for the right part.
Then $\bar{\sigma}\in {\bf B}\cup\dt{\bf B}$ and $\bar{\eta}\in {\bf W}_{\varepsilon,p}^+$ by the properties {\bf (\ref{asmpV1})} and {\bf (\ref{asmpV2})} of differentiable sectors.
Obviously, $\bar{\sigma}\preceq\bar{\tau}$ since $\sigma$ is a subtree of $\tau$. The equality $\bar{\sigma}=\bar{\tau}$ is the case only if $\sigma=\tau$ and ${\fr n}_\sigma={\fr n}$, which corresponds to the term $\bar{\tau}\otimes{\bf1}$ in the decomposition of $\Delta\bar{\tau}$. Thus $\bar{\sigma}\prec\bar{\tau}$ for the other terms.
The quotient graph $\bar{\eta}$ is of the form
$$
X^k\prod_{j=1}^n\mcI_{k_j}^{{\fr l}_j}(\sigma_j).
$$
If ${\fr l}_j=\bbOmega$, then $\sigma_j=\bullet$ because of the property {\bf(\ref{asmpleaf})} of differential sectors, but $r_{\varepsilon,p}(\mcI_k^{\bbOmega}(\bullet))=r_0-\varepsilon-|k|_\mfs$ is always negative.
Therefore, all labels ${\fr l}_j$ are $\bbK$ except that at most one ${\fr l}_j$ may be $\bbH$ and then $\sigma_j=\bullet$.
Since $\sigma_j$ is a strict subgraph of $\bar{\tau}$, we have $\sigma_j\prec\bar{\tau}$.
\end{Dem}

\ssk

Recall that ${\bf B}=\{X^k\}_{|k|_\mfs<L}\cup{\bf B}_{\ocircle}$ from the property {\bf(\ref{asmpV0})} of differentiable sectors in Definition \ref{asmpV}. We align the elements of ${\bf B}_{\ocircle}$ as
$$
{\bf B}_{\ocircle} = \big\{\tau_1\preceq\tau_2\preceq\cdots\preceq\tau_N\big\},
$$
with $\tau_1=\ocircle$. Although there may be different choices for this ordering the previous choice has no consequence on Theorem \ref{thm:main} as it is only used as a tool in the proof of this result. %Recall the integer $L$ was introduced in Definition \ref{asmpV}. 
For $0\leq i\leq N$ set
$$
{\bf B}_i \defeq \big\{\tau_j\,;\,j\le i\big\}\cup\{X^k\}_{|k|_\mfs<L}.
$$
Then we observe that ${\bf B}\cap{\bf C}_{\prec\tau_i}\subset{\bf B}_{i-1}$, and if $\tau\in\dt{\bf B}_i$, then ${\bf C}_{\prec\tau}\subset{\bf B}_{i-1}\cup\dt{\bf B}_{i-1}$. The first observation is elementary. The second observation can be shown as follows by a case distinction. Write $\tau=D_e\mu\in\dt{\bf B}_i$ for some $\mu\in{\bf B}_i$ and $\bbOmega$-labeled edge $e\in E_\mu$. If $\sigma\prec\tau$ and $\sigma\in{\bf B}$, then since $|\sigma|_\ocircle \leq |\tau|_\ocircle < |\mu|_\ocircle$, we have $\sigma\in{\bf B}\cap{\bf C}_\mu\subset{\bf B}_{i-1}$. If $\sigma\in\dot{\bf B}$, we can write $\sigma=D_{e'}\mu'$ for some $\mu'\in{\bf B}$ and $\bbOmega$-labeled edge $e'\in E_{\mu'}$. Then $\mu'\prec\mu\in{\bf B}_i$ by the definition of $\prec$.

The following properties can be obtained by applying the above observations to Lemma \ref{leminduction*}.

\ssk

\begin{lem} \label{leminduction}
For $1\le i\le N$, set
$$
V_i \defeq \spa({\bf B}_i),\qquad
W_i\defeq \spa({\bf B}_{i-1}\cup\dt{\bf B}_i).
$$
Moreover, for $0\le i\le N$, $\eps\ge0$, and $p\in[2,\infty]$, define $V_{i,\varepsilon}^+$ as the subalgebra of $V_{\varepsilon}^+$ generated by the symbols 
$$
{\bf V}_{i,\varepsilon}^+\defeq\{X_j\}_{j=1}^d\cup\{\mcI_k^{\bbK}(\tau)\}_{k\in\bbN^d,\, \tau\in{\bf B}_i\setminus\{X^l\}_l,\, r_{\varepsilon,\infty}(\tau)+\beta_0>|k|_\mfs}
$$
and define $W_{i,\varepsilon,p}^+$ as the subalgebra of $W_{\varepsilon,p}^+$ generated by the symbols 
$$
{\bf W}_{i,\varepsilon,p}^+\defeq\{X_j\}_{j=1}^d\cup\{\odot_k\}_{|k|_\mfs<r_{\varepsilon,p}(\bbH)}\cup\{\mcI_k^{\bbK}(\tau)\}_{k\in\bbN^d,\, \tau\in({\bf B}_i\cup\dt{\bf B}_i)\setminus\{X^l\}_l,\, r_{\varepsilon,p}(\tau)+\beta_0>|k|_\mfs}.
$$
Then 
$$
\scV_{i,\varepsilon} \defeq \Big((V_i,\Delta_{\varepsilon,\infty}) , \big(V_{i-1}^+,\Delta_{\varepsilon,\infty}^+\big)\Big),  \qquad  \scW_{i,\varepsilon,p} \defeq \Big((W_i,\Delta_{\varepsilon,p}) , \big(W_{i-1,\varepsilon,p}^+,\Delta_{\varepsilon,p}^+\big)\Big)
$$ 
are concrete regularity-integrability substructures of $\scW_{\varepsilon,p}$. 
\end{lem}

\ssk

We recall the following example which already appeared in the introduction.
%The following spaces would for instance be involved in the study of the three-dimensional parabolic Anderson model equation. 
We adopt below the notation $\{X^k\}_k$ for $\{X^k\}_{|k|_\mfs<L}$.
\begin{align}\label{EqV2W2}
\begin{aligned}
&W_1 = \spa\big(\{\odot\}\cup\{X^{k}\}_{k}\big), 
&&V_1 = \spa\big(\{\ocircle\}\cup\{X^{k}\}_{k}\big),   \\
&W_2 = \spa\big(\big\{\ocircle, \odot, \begin{tikzpicture}[scale=0.3,baseline=0.05cm] \node at (0,0)  [noise] (1) {}; \node at (0,1.1)  [noise] (2) {}; \node at (0,1.1)  [dot] {}; \draw[K] (1) to (2); \end{tikzpicture}\hspace{0.03cm},\hspace{0.03cm} \begin{tikzpicture}[scale=0.3,baseline=0.05cm] \node at (0,0)  [noise] (1) {}; \node at (0,1.1)  [noise] (2) {}; \node at (0,0)  [dot] {}; \draw[K] (1) to (2); \end{tikzpicture}\hspace{0.03cm}\big\}\cup\{X^{k}\}_{k}\big),
&&V_2 = \spa\big(\big\{\ocircle, \begin{tikzpicture}[scale=0.3,baseline=0.05cm] \node at (0,0)  [noise] (1) {}; \node at (0,1.1)  [noise] (2) {}; \draw[K] (1) to (2); \end{tikzpicture}\hspace{0.03cm}\big\}
\cup\{X^{k}\}_{k}
\big),
\\
&W_3 = \spa\big(\big\{\ocircle, \begin{tikzpicture}[scale=0.3,baseline=0.05cm] \node at (0,0)  [noise] (1) {}; \node at (0,1.1)  [noise] (2) {}; \draw[K] (1) to (2); \end{tikzpicture}\hspace{0.03cm}, \odot, \begin{tikzpicture}[scale=0.3,baseline=0.05cm] \node at (0,0)  [noise] (1) {}; \node at (0,1.1)  [noise] (2) {}; \node at (0,1.1)  [dot] {}; \draw[K] (1) to (2); \end{tikzpicture}\hspace{0.03cm},\hspace{0.03cm} \begin{tikzpicture}[scale=0.3,baseline=0.05cm] \node at (0,0)  [noise] (1) {}; \node at (0,1.1)  [noise] (2) {}; \node at (0,0)  [dot] {}; \draw[K] (1) to (2); \end{tikzpicture}\hspace{0.03cm}, \begin{tikzpicture}[scale=0.3,baseline=0.05cm] \node at (0,0)  [noise] (1) {}; \node at (-0.5,1.1) [noise] (2) {}; \node at (0.5,1.1) [noise] (3) {}; \node at (-0.5,1.1)  [dot] {}; \draw[K] (1) to (2); \draw[K] (1) to (3); \end{tikzpicture}\hspace{0.03cm} , \begin{tikzpicture}[scale=0.3,baseline=0.05cm] \node at (0,0)  [noise] (1) {}; \node at (-0.5,1.1) [noise] (2) {}; \node at (0.5,1.1) [noise] (3) {}; \node at (+0.5,1.1)  [dot] {}; \draw[K] (1) to (2); \draw[K] (1) to (3); \end{tikzpicture}\hspace{0.03cm}, \begin{tikzpicture}[scale=0.3,baseline=0.05cm] \node at (0,0)  [noise] (1) {}; \node at (-0.5,1.1) [noise] (2) {}; \node at (0.5,1.1) [noise] (3) {}; \node at (0,0)  [dot] {}; \draw[K] (1) to (2); \draw[K] (1) to (3); \end{tikzpicture}\hspace{0.03cm}\big\}\cup\{X^{k}\}_{k}\big),  \quad
&&V_3 = \spa\big(\big\{\ocircle, \begin{tikzpicture}[scale=0.3,baseline=0.05cm] \node at (0,0)  [noise] (1) {}; \node at (0,1.1)  [noise] (2) {}; \draw[K] (1) to (2); \end{tikzpicture}\hspace{0.03cm}, \begin{tikzpicture}[scale=0.3,baseline=0.05cm] \node at (0,0)  [noise] (1) {}; \node at (-0.5,1.1) [noise] (2) {}; \node at (0.5,1.1) [noise] (3) {}; \draw[K] (1) to (2); \draw[K] (1) to (3); \end{tikzpicture}\hspace{0.03cm}\big\}\cup\{X^{k}\}_{k}\big).
\end{aligned}
\end{align}

\bigskip

{\it \S2 The flow of the proof.} We need some notations to describe the flow of the proof. Let $\big\{{\sf M}^{\xi_n,h_n,R_n;\varepsilon,p}\big\}_{n\geq 0}$ be the sequence of random models on $\scW_{\varepsilon,p}$ under consideration. Set
\begin{align}
\label{*eq:tobebounded1}
\|{\sf M}^{\xi_n,R_n;\varepsilon}\|_{\textbf{\textsf{M}}(\scV_{i,\varepsilon})_{w_c}} 
&\defeq \|{\sf\Pi}^{\xi_n,R_n;\varepsilon} \res {\bf B}_i\|_{w_c}  
+ \|{\sf g}^{\xi_n,R_n;\varepsilon} \res {\bf V}_{i-1,\varepsilon}^+\|_{w_c},\\
\label{*eq:tobebounded2}
\|{\sf M}^{\xi_n,h_n,R_n;\varepsilon,p}\|_{\textbf{\textsf{M}}(\scW_{i,\varepsilon,p})_{w_c}} 
&\defeq \|{\sf\Pi}^{\xi_n,h_n,R_n;\varepsilon,p} \res {\bf B}_{i-1}\cup\dt{\bf B}_i\|_{w_c}  
+ \|{\sf g}^{\xi_n,h_n,R_n;\varepsilon,p} \res {\bf W}_{i-1,\varepsilon,p}^+\|_{w_c}. 
\end{align}
For any fixed $p\in[2,\infty]$ we write below $\textsf{\textbf{bd}}(\scW,i,p)$ to mean that
\begin{equation*} %\label{EqExpectaionSupWi}
\sup_{n\in\bbN}\bbE\bigg[\sup_{\Vert h\Vert_H\leq 1} \Vert {\sf M}^{\xi_n,h_n,R_n;\varepsilon,p} \Vert_{\textbf{\textsf{M}}(\scW_{i,\varepsilon,p})_{w_c}}^q\bigg]<\infty
\end{equation*}
for any small $\varepsilon>0$, large $c>|\mfs|$ and $q\in[1,\infty)$. Similarly, we write $\textsf{\textbf{cv}}(\scW,i,p)$ to mean that 
$$
\lim_{n,m\to\infty}\bbE\bigg[\sup_{\Vert h\Vert_H\leq 1} \Vert{\sf M}^{n;\varepsilon,p} \hspace{-0.03cm}:\hspace{-0.03cm} {\sf M}^{m;\varepsilon,p}\Vert_{\textbf{\textsf{M}}(\scW_{\varepsilon,p})_{w_c}}^q\bigg]=0
$$
for any small $\varepsilon>0$, large $c>|\mfs|$ and all $q\in[1,\infty)$. This implies that the restriction of models ${\sf M}^{\xi_n,h_n,R_n;\varepsilon,p}$ on $\scW_{i,\varepsilon,p}$ converges in $L^q\big(\Omega,\bbP ; \textbf{\textsf{M}}(\scW_{i,\varepsilon,p})_{w_c}\big)$ as $n$ goes to $\infty$. We also write
$$
\big\{\textsf{\textbf{bd}}(\scW,i,p)\big\}_{2\leq p\leq \infty}, \quad  \textrm{resp.}  \quad \big\{\textsf{\textbf{cv}}(\scW,i,p)\big\}_{2\leq p\leq \infty},
$$ 
to mean each statement holds for any $p\in[2,\infty]$. We also write 
$$
\textsf{\textbf{bd}}(\scV,i), \quad  \textrm{resp.}  \quad \textsf{\textbf{cv}}(\scV,i),
$$ 
to mean the similar statements to $\textsf{\textbf{bd}}(\scW,i,p)$, resp. $\textsf{\textbf{cv}}(\scW,i,p)$, with $\scV_{i,\varepsilon}$ in place of $\scW_{i,\varepsilon,p}$. These statements do not depend on $p$ so we do not record this parameter in the notation.

\medskip

We first show in Section \ref{sec:ind0} the initial statement $\{\textsf{\textbf{cv}}(\scW,1,p)\}_{2\leq p\leq \infty}$, that is the convergence result for the single tree $\odot$. Then the induction proceeds in three steps that can schematically be described as follows. We give below each step what are the main type of arguments involved.   \vspace{0.15cm}

\ssk

	\textbf{{Step 1}}: $\textsf{\textbf{cv}}(\scW,i,\infty) \dashrightarrow \textsf{\textbf{cv}}(\scV,i)$
	\begin{itemize}
		\item[] \hspace{0.4cm} Reconstruction theorem and $\big($Lemma \ref{Lem:dPi=PiD}, the Spectral gap assumption and \eqref{EqControlExpectationQ1}$\big)$  \vspace{0.15cm}
	\end{itemize} 

	\textbf{{Step 2}}: $\textsf{\textbf{cv}}(\scW,i,p),\textsf{\textbf{cv}}(\scV,i) \longrightarrow \text{$\sf g$-part of }\textsf{\textbf{cv}}(\scW,i+1,p)$
	\begin{itemize}
		\item[] \hspace{0.4cm} Multilevel Schauder estimate for an appropriate modelled distribution   \vspace{0.15cm}
	\end{itemize} 

%\vfill \pagebreak		

	\textbf{{Step 3}}:  \\
	
	$\left\{ \begin{aligned} 
	\textsf{\textbf{cv}}(\scW,i,2),\textsf{\textbf{cv}}(\scV,i) \longrightarrow
	\text{$\sf\Pi$-part of } \textsf{\textbf{cv}}(\scW,i+1,2)&   \\
	\big\{\textsf{\textbf{cv}}(\scW,i,p)\big\}_{2\leq p\leq \infty}&
	\end{aligned}\right\}$
	$\longrightarrow$ 
	$\begin{aligned}
	&\text{$\sf \Pi$-part of}\\
	&\big\{\textsf{\textbf{cv}}(\scW,i+1,p)\big\}_{2\le p\le\infty}
	\end{aligned}$
		\begin{itemize}
		\item[] \hspace{0.4cm}  Reconstruction theorem and Lemma \ref{LemComparisonFormula}   \vspace{0.2cm}
	\end{itemize} 

The {\it dashed line} in Step 1 is used to emphasize that this is a {\it probabilistic} step: We obtain some stochastic estimates for the models on $\scV_{i,\varepsilon}$ in terms of some stochastic estimates for the models on $\scW_{i,\delta,\infty}$ for some $\delta<\varepsilon$. On the other hand, we use some {\it solid lines} in Step 2 and Step 3 to emphasize that they are {\it deterministic} steps: We obtain some $\omega$-wise estimates for the models on $\scW_{i+1,\varepsilon,p}$ in terms of some $\omega$-wise estimates for the models on $\{\scW_{i,\delta,q_a}\}_a$ and $\scV_{i,\delta}$ for some finite set $\{q_a\}\subset[2,\infty]$ and $\delta<\varepsilon$. Step 1 builds the model on `trees' from the model on `derivative trees'. Step 2 and 3 build the model on new `derivate trees'. 

In the setting of example \eqref{EqV2W2} we would first construct the deterministic model on $\odot$ (the initial case), then the random model on $\ocircle$ (Step 1), and then successively on $\begin{tikzpicture}[scale=0.3,baseline=0.05cm] \node at (0,0)  [noise] (1) {}; \node at (0,1.1)  [noise] (2) {}; \node at (0,1.1)  [dot] {}; \draw[K] (1) to (2); \end{tikzpicture}\hspace{0.03cm},\hspace{0.03cm} \begin{tikzpicture}[scale=0.3,baseline=0.05cm] \node at (0,0)  [noise] (1) {}; \node at (0,1.1)  [noise] (2) {}; \node at (0,0)  [dot] {}; \draw[K] (1) to (2); \end{tikzpicture}\;$ (Steps 2 and 3), then $\begin{tikzpicture}[scale=0.3,baseline=0.05cm] \node at (0,0)  [noise] (1) {}; \node at (0,1.1)  [noise] (2) {}; \draw[K] (1) to (2); \end{tikzpicture}\;$ (Step 1), then $\begin{tikzpicture}[scale=0.3,baseline=0.05cm] \node at (0,0)  [noise] (1) {}; \node at (-0.5,1.1) [noise] (2) {}; \node at (0.5,1.1) [noise] (3) {}; \node at (-0.5,1.1)  [dot] {}; \draw[K] (1) to (2); \draw[K] (1) to (3); \end{tikzpicture}\hspace{0.03cm} , \begin{tikzpicture}[scale=0.3,baseline=0.05cm] \node at (0,0)  [noise] (1) {}; \node at (-0.5,1.1) [noise] (2) {}; \node at (0.5,1.1) [noise] (3) {}; \node at (+0.5,1.1)  [dot] {}; \draw[K] (1) to (2); \draw[K] (1) to (3); \end{tikzpicture}\hspace{0.03cm}, \begin{tikzpicture}[scale=0.3,baseline=0.05cm] \node at (0,0)  [noise] (1) {}; \node at (-0.5,1.1) [noise] (2) {}; \node at (0.5,1.1) [noise] (3) {}; \node at (0,0)  [dot] {}; \draw[K] (1) to (2); \draw[K] (1) to (3); \end{tikzpicture}\;$ (Steps 2 and 3), 
finally on $\begin{tikzpicture}[scale=0.3,baseline=0.05cm] \node at (0,0)  [noise] (1) {}; \node at (-0.5,1.1) [noise] (2) {}; \node at (0.5,1.1) [noise] (3) {}; \draw[K] (1) to (2); \draw[K] (1) to (3); \end{tikzpicture}\;$ (Step 1). 

Step 2 builds the ${\sf g}$-part of the model on the new derivative trees. Step 3 builds the $\sf \Pi$-part of the model on these `derivative trees' $\tau$. It is split in two substeps. Step 3(a) builds ${\sf \Pi}_x^{R;2}(\tau)$; Step 3(b) builds ${\sf \Pi}_x^{R;p}(\tau)$ for all $1\leq p\leq \infty$ from ${\sf \Pi}_x^{R;2}(\tau)$, the comparison result of Lemma \ref{LemComparisonFormula} and the inductive assumptions.

\medskip

The mechanics of Step $1$ is explained below in Section \ref{sec:ind1}. Step $2$ is detailed in Section \ref{sec:ind2}, while Section \ref{sec:ind3} is dedicated to Step $3$. We show that the limit model depends on the finite collection of numbers $\big\{\lim_{n\to\infty} \bbE\big[\mcQ_1(0, {\sf\Pi}_0^{\xi_n,R_n;0}\tau)\big] \big\}_{\tau\in{\bf B},\,r_{0,\infty}(\tau)\le0}$ in Section \ref{sec:ind5}. The proof of Corollary \ref{cor:BPHZ} is in Section \ref{sec:ind6}.

We only describe in the sections \ref{sec:ind1} to \ref{sec:ind3} the flow of the proof of Theorem \ref{thm:main} and defer the proof of a number of lemmas to Section {{\ref{SectionLemmas}}}. To lighten the notations we will suppress from the notations the exponents $\xi_n,h_n,R_n$ in the remainder of this section. So we simply write 
$$
{\sf M}^{n;\varepsilon,p}=({\sf \Pi}^n,{\sf g}^{n;\varepsilon,p})
$$
rather than ${\sf M}^{\xi_n,h_n,R_n;\varepsilon,p}$. Occasionally we write  ${\sf M}^{n,h;\varepsilon,p} = ({\sf \Pi}^{n,h},{\sf g}^{n,h;\varepsilon,p})$ to emphasize the dependence on $h$. 

\ssk

We first prove the $n$-uniform bound \eqref{EqExpectationSupH} by proving the $\textsf{\textbf{bd}}$-versions of Steps 1-3 -- this is the content of Sections {{\ref{sec:ind0}}}-{{\ref{sec:ind3}}}. We use these uniform bounds together with some local Lipschitz estimates satisfied by the reconstruction and integration maps to prove in Section {{\ref{sec:ind4}} the convergence result. The following notations will be useful in the course of the proof. 

\medskip

\noindent \textbf{Notation.} {\it Set for all $x,y\in\bbR^d$
$$
{\sf\Gamma}_{yx}^{n;\varepsilon,p} \defeq \big(\id\otimes\mathop{{\sf g}_{yx}^{n;\varepsilon,p}}\big) \Delta_{\varepsilon,p}.
$$
These operators leave each space $V_i$ and $W_i$ stable. The pair of $({\sf\Pi}^{n;\varepsilon,p},{\sf\Gamma}^{n;\varepsilon,p})$ is called a model in the original terminology of \cite{Hai14}. For any $\rho\in{\bf B}\cup\overset{\text{\tiny$\mybullet$}}{{\bf B}}$ denote by $P_\rho : W\rightarrow \bbR$ the canonical projection to the $\rho$-coefficient. 
For any $\tau\in{\bf B}\cup\dt{\bf B}$ we define the quantity
\begin{align} \label{EqDefnSizeGammaVi}
\|{\sf\Gamma}^{n;\varepsilon,p} \res \tau \|_{w_c} &\defeq \max_{\sigma\prec\tau} \sup_{y\in\bbR^d\setminus\{0\}} \frac{\big\|P_\sigma{\sf\Gamma}_{(x+y)x}^{n;\varepsilon,p}\tau\big\|_{L_x^{i_p(\tau):i_p(\sigma)}(w_c)}}{w_{c}(y)^{-1}\|y\|_{\fr s}^{r_{\varepsilon,p}(\tau) - r_{\varepsilon,p}(\sigma)}}.
\end{align}
Here the exponent $p:q$ for $1\le p\le q\le\infty$ is defined at the beginning of Appendix \ref{SectionAppendix}: $i_p(\tau):i_p(\sigma)=p$ if $\tau\in\dt{\bf B}$ and $\sigma\in{\bf B}$, otherwise $i_p(\tau):i_p(\sigma)=\infty$.
Moreover, for any ${\bf C}\subset{\bf B}\cup\dt{\bf B}$, write 
$$
\|{\sf\Gamma}^{n;\varepsilon,p} \res {\bf C} \|_{w_c}\defeq\max_{\tau\in{\bf C}}\|{\sf\Gamma}^{n;\varepsilon,p} \res \tau \|_{w_c}
$$
With $m_{\bf B}\defeq\max\{|E_\tau|\,;\,\tau\in{\bf B}\}$, we obtain from Lemma \ref{leminduction*} and H\"older inequality the estimates
\begin{equation} \label{EqControlGamma} \begin{split}
\big\|{\sf\Gamma}^{n;\varepsilon} \res {\bf B}_i\big\|_{w_{cm_{\bf B}}} 
&\lesssim \big(1 + \|{\sf g}^{n;\varepsilon} \res {\bf V}_{i-1,\eps}^+\|_{w_c}\big)^{m_{\bf B}},   \\
\big\|{\sf\Gamma}^{n;\varepsilon,p} \res \dt{\bf B}_i\big\|_{w_{cm_{\bf B}}} &\lesssim \big(1 + \|{\sf g}^{n;\varepsilon,p} \hspace{-0.03cm}:\hspace{-0.03cm} {\bf W}_{i-1,\varepsilon,p}^+\|_{w_c}\big)^{m_{\bf B}}.
\end{split} \end{equation}
}

%\ssk

%%------------------------------------------------------------------%%
\subsection{The initial case $\{\textsf{bd}(\scW,1,p)\}_{2\leq p\leq \infty}$$\boldmath .$ \hspace{0.3cm}}
\label{sec:ind0}
%%------------------------------------------------------------------%%

For the initial case one has $W_1 = \spa\big(\{\odot\}\cup\{X^{k}\}_{|k|_\mfs<L}\big)$ and ${\bf W}_{0,\varepsilon,p}^+ =\{X_j\}_{j=1}^d\cup\{\odot_k\}_{|k|_\mfs<r_{\varepsilon,p}(\bbH)}$.
By the direct calculation, we have the formulas
\begin{equation}\label{*Eq:ModelXH} \begin{aligned}
({\sf\Pi}_x^{n;\varepsilon,p}X^k)(y)&=(y-x)^k,\qquad
&({\sf\Pi}_x^{n;\varepsilon,p}\odot)(y)&=h_n(y)-\sum_{|k|_\mfs<r_{\varepsilon,p}(\bbH)}\frac{(y-x)^k}{k!}\partial^kh_n(x),   \\
\end{aligned} \end{equation}

\begin{equation} %\label{*Eq:ModelXH} 
\begin{aligned}
{\sf g}_{yx}^{n;\varepsilon,p}(X_j)&=y_j-x_j,
&{\sf g}_{yx}^{n;\varepsilon,p}(\odot_k)&=\partial^kh_n(y)-\sum_{|k+l|_\mfs<r_{\varepsilon,p}(\bbH)}\frac{(y-x)^l}{l!}\partial^{k+l}h_n(x).
\end{aligned} \end{equation}
Since the bounds on $X^k$ are easy, so it is sufficient to consider $\odot$ and $\odot_k$.
When $r_{\varepsilon,p}(\bbH)\le0$, only ${\sf\Pi}_x^{n;\varepsilon,p}\odot=h_n$ is the matter. Then by the embedding \eqref{*eq:interpolateHOmega},
$$
\|{\sf\Pi}^{n;\varepsilon,p}\res\{\odot\}\|_{w_c}\le
\|{\sf\Pi}^{n;\varepsilon,p}\res\{\odot\}\|_{w_0}
=\|h_n\|_{B_{p,\infty}^{r_0+|\mfs|/p-\varepsilon}(w_0)}
\lesssim\|h_n\|_{H}\le1
$$
for any $c\ge0$, $\varepsilon\ge0$, and $p\in[2,\infty]$.
When $r_{\varepsilon,p}(\bbH)>0$, we also obtain the required bounds by the following facts on Besov spaces. We believe that it is a well known fact at least for isotropic cases, but we provide a proof in Section \ref{SectionBesovTaylor} for self-containedness.

\ssk

\begin{lem}\label{*lem:BesovTaylor}
Let $p\in[1,\infty]$ and $r\in(0,\infty)\setminus\bbN[\mfs]$, where $\bbN[\mfs]\defeq\{|k|_\mfs\,;\,k\in\bbN^d\}$. For any $h\in B_{p,\infty}^r(w_0)$, we have
$$
\Bigg\|h(x+y)-\sum_{|k|_\mfs<r}\frac{y^k}{k!}\partial^kh(x)\Bigg\|_{L_x^p(w_0)}
\lesssim\|h\|_{B_{p,\infty}^r(w_0)}\|y\|_\mfs^r.
$$
\end{lem}

%\ssk

%%---------------------------------------------------------------------------------------------------------------------------------------------------------------------------------------------------%%
\subsection{Step 1: From $\textsf{bd}(\scW,i,\infty)$ to $\textsf{bd}(\scV,i)\boldmath .$ \hspace{0.3cm}}
\label{sec:ind1}
%%---------------------------------------------------------------------------------------------------------------------------------------------------------------------------------------------------%%

We need to estimate
$$
\|{\sf M}^{n;\varepsilon}\|_{\textbf{\textsf{M}}(\scV_{i,\varepsilon})_{w_c}} = \|{\sf\Pi}^{n;\varepsilon} \res {\bf B}_i\|_{w_c} + \|{\sf g}^{n;\varepsilon} \res {\bf V}_{i-1,\varepsilon}^+\|_{w_c}.
$$
Since ${\bf B}_i={\bf B}_{i-1}\cup\{\tau_i\}$ and ${\bf V}_{i-1}^+\subset {\bf W}_{i-1,\varepsilon,\infty}^+$, it is sufficient to prove the appropriate bound for
$$
\sup_{n\in\bbN}\bbE\Big[\big\Vert {\sf\Pi}^{n;\varepsilon,\infty} \hspace{-0.03cm}:\hspace{-0.03cm} \tau_i\big\Vert_{w_c}^q\Big]
$$
in terms of the assumed boundedness results on $\scW_{i,\varepsilon,\infty}$, that is \eqref{*eq:tobebounded2} with $p=\infty$. We use different arguments depending on the sign of $r_{0,\infty}(\tau_i)$. It turns out that if $r_{0,\infty}(\tau_i) > 0$ we have an $\omega$-wise bound of $\Vert {\sf\Pi}^{n;\varepsilon,\infty} \hspace{-0.03cm}:\hspace{-0.03cm} \tau_i\Vert_{w_c}$ in terms of $\omega$-wise bounds of ${\sf M}^{n;\varepsilon,\infty}$ on $\scW_{i,\varepsilon,\infty}$. If $r_{0,\infty}(\tau_i) \leq 0$ we only have an $L^q(\Omega)$ control of $\Vert {\sf\Pi}^{n;\varepsilon,\infty}  \hspace{-0.03cm}:\hspace{-0.03cm} \tau_i\Vert_{w_c}$. Lemmas \ref{lem:Biposi} to \ref{lem:Binega} for this step below are proved in Section \ref{SectionLastProofsLemmas}.

\ssk

-- {\sl First case.} If $r_{0,\infty}(\tau_i) > 0$, then $r_{\eps,\infty}(\tau_i) > 0$ for any $\varepsilon\in[0,\eps_0)$, hence the proof reduces to an application of the reconstruction theorem (see Theorem \ref{thm:besovreconstruction}). Indeed the modelled distribution
$$
f^{n;\varepsilon,\infty}_{\tau_i}(x) \defeq \big(\id\otimes\mathop{{\sf g}^{n;\varepsilon,\infty}_x}\big) \Delta_{\varepsilon,\infty}\tau_i - \tau_i
$$
belongs in that case to $D^{(r_{\varepsilon,\infty}(\tau_i),\infty)}(V_{i-1} ; {\sf \Gamma}^{n;\varepsilon,\infty})_{w_{c_{\bf B}}}$, with a positive regularity exponent. Then its unique reconstruction coincides with ${\sf\Pi}^n\tau_i$. 
This fact is important and it will appear several times later in slightly different settings, so we present it here in a general form.

\ssk

\begin{lem} \label{*argumentforreconst}
Let $\tau\in{\bf B}\cup\dt{\bf B}$, $\varepsilon\ge0$, and $p\in[2,\infty]$. The function
$$
f^{n;\varepsilon,p}_{\tau}(x) \defeq \big(\id\otimes\mathop{{\sf g}^{n;\varepsilon,p}_x}\big) \Delta_{\varepsilon,p}\tau - \tau
$$
is a modelled distribution in the class $D^{(r_{\varepsilon,p}(\tau),i_p(\tau))}(C_{\prec\tau} ; {\sf \Gamma}^{n;\varepsilon,p})_{w_{cm_{\bf B}}}$, provided that ${\sf M}^{n;\varepsilon,p}$ is a model in the class $\textbf{\textsf{M}}(\scW_{\varepsilon,p})_{w_c}$. Moreover ${\sf\Pi}^n\tau$ is a reconstruction of $f_{\tau}^{n;\varepsilon,p}$.
\end{lem}

\ssk

As an immediate consequence of Lemma \ref{*argumentforreconst} we have the following result.

\ssk

\begin{lem} \label{lem:Biposi}
Let $r_{0,\infty}(\tau_i) > 0$. For any $\varepsilon\in[0,\eps_0)$ and $c>0$, there exists a positive constant $C$ which is independent of $n$ and $\omega$, one has
$$
\|{\sf\Pi}^{n;\varepsilon,\infty} \hspace{-0.03cm}:\hspace{-0.03cm} \tau_i\|_{w_{c(m_{\bf B}+1)}} 
\le C \|{\sf\Pi}^{n;\varepsilon,\infty} \hspace{-0.03cm}:\hspace{-0.03cm} {\bf B}_{i-1}\|_{w_c} \, \|{\sf\Gamma}^{n;\varepsilon,\infty} \hspace{-0.03cm}:\hspace{-0.03cm} {\bf B}_i\|_{w_{cm_{\bf B}}}.
$$
\end{lem}

\ssk

It should be noted that the constant $C$ is the one appearing in the reconstruction theorem (Theorem \ref{thm:besovreconstruction}), so it is universal and depends neither on $n$ nor on $\omega$. Moreover, recall from \eqref{EqControlGamma} that the bounds of $\|{\sf\Gamma}^{n;\varepsilon,\infty} \res {\bf B}_i\|_{w_{cm_{\bf B}}}$ is obtained deterministically from the bounds on $\|{\sf g}^{n;\varepsilon,\infty} \res {\bf V}_{i-1,\eps}^+\|_{w_c}$, which is contained in the assumption $\textsf{\textbf{bd}}(\scW,i,\infty)$. 

\ssk

-- {\sl Alternative case.} Consider next the case where $r_{0,\infty}(\tau_i) \leq 0$, so $r_{\varepsilon,\infty}(\tau_i) < 0$ if $\varepsilon>0$. We cannot use the above reconstruction argument as $f_{\tau_i}^{n;\varepsilon,\infty}$ no longer has a unique reconstruction. Instead we use the spectral gap inequality and the algebraic identity 
\begin{align*}
\big(d_\omega{\sf\Pi}^{n;\varepsilon,\infty}_x\tau_i\big)(h) &= \big(d_{\xi_n}{\sf\Pi}^{\xi_n,R_n;\varepsilon,\infty}_x\tau_i\big)\big((d_\omega \xi_n)(h)\big)\\
&= {\sf\Pi}_x^{\xi_n,h_n,R_n;\varepsilon,\infty}(D\tau_i)
= {\sf\Pi}^{n,h;\varepsilon,\infty}_x({D}\tau_i),
\end{align*}
which follows from Lemma \ref{Lem:dPi=PiD} and the chain rule. Then for any finite exponent $q=2^r$ one gets from the inequality \eqref{EqSGInequalityModified} that
\begin{align}\label{*eq:resultofSGineq}
\bbE\big[\big|\mcQ_t(x, {\sf\Pi}^{n;\varepsilon,\infty}_x\tau_i)\big|^q\big] 
&\lesssim_q \big|\bbE\big[\mcQ_t(x, {\sf\Pi}^{n;\varepsilon,\infty}_x\tau_i)\big]\big|^q + \bbE\bigg[\sup_{\|h\|_{H}\le1}\big|\mcQ_t\big(x, {\sf\Pi}_x^{n,h;\varepsilon,\infty}({D}\tau_i)\big)\big|^q\bigg].
\end{align}
The following result holds for the expectation part. Define the quantity
$$
E_i^n \defeq \bbE\big[\mcQ_1\big(0,{\sf\Pi}^{n;\varepsilon,\infty}_0(\tau_i)\big)\big]= \bbE\big[\mcQ_1\big(0,{\sf\Pi}^{n;0,\infty}_0(\tau_i)\big)\big];
$$
it is uniformly bounded over $n$ by assumption.

\ssk

\begin{lem} \label{lem:stationary}
Let $r_{0,\infty}(\tau_i) \le 0$. For any $\varepsilon>0$, there exists a positive constant $C$ which is independent of $n$, one has
\begin{align*}
\big\|\bbE\big[\mcQ_t(x,{\sf\Pi}^{n;\varepsilon,\infty}_x\tau_i)\big]\big\|_{L_x^\infty(w_{c(m_{\bf B}+1)})} \le
C\Big( \vert E_i^n\vert+t^{\frac{r_{\varepsilon,\infty}(\tau_i)}\ell} \, \bbE\big[\|{\sf\Pi}^{n;\varepsilon,\infty} \res {\bf B}_{i-1}\|_{w_c} \, \|{\sf\Gamma}^{n;\varepsilon,\infty} \res \tau_i\|_{w_{cm_{\bf B}}}\big] \Big).
\end{align*}
\end{lem}

\ssk
We emphasize that the proof of this lemma is the only place in this work where we use the assumption that the law of the noise is invariant by translation. From the estimate of the lemma and the spectral gap inequality we obtain the bound
\begin{align*}
\bbE\big[\big|&\mcQ_t(x, {\sf\Pi}^{n;\varepsilon,\infty}_x\tau_i)\big|^qw_{c(m_{\bf B}+2)}(x)^q\big]   
\\
&\lesssim w_c(x)^q\bigg(\big|\bbE\big[\mcQ_t(x,{\sf\Pi}^{n;\varepsilon,\infty}_x\tau_i)\big]\big|^qw_{c(m_{\bf B}+1)}(x)^q 
+ \bbE\bigg[\sup_{\|h\|_{H}\le1}\big|\mcQ_t\big(x,{\sf\Pi}_x^{n,h;\varepsilon,\infty}({D}\tau_i)\big)\big|^qw_c(x)^q\bigg]\bigg)   
\\
&\lesssim w_c(x)^q\bigg(|E_i^n|^q + t^{\frac{r_{\varepsilon,\infty}(\tau_i)}\ell q} \, \bbE\big[\|{\sf\Pi}^{n;\varepsilon,\infty}\res {\bf B}_{i-1}\|_{w_c} \, \|{\sf\Gamma}^{n;\varepsilon,\infty} \res \tau_i\|_{w_{cm_{\bf B}}}\big]^q   \\
&\hspace{6.8cm}+ t^{\frac{r_{\varepsilon,\infty}(D\tau_i)}\ell q} \, \bbE\bigg[\sup_{\|h\|_{H}\le1} \|{\sf\Pi}^{n,h;\varepsilon,\infty} \hspace{-0.03cm}:\hspace{-0.03cm} \dt{\bf B}_i\|_{w_c}^q\bigg]\bigg).
\end{align*}
Recall that $D$ preserves the $r_{\varepsilon,\infty}$ degree: $r_{\varepsilon,\infty}(D\tau_i)=r_{\varepsilon,\infty}(\tau_i)$. If $c>|\mfs|$ then, by integrating the above estimate over $x$, we get 
the bound
\begin{align}\label{eq:Binega} \begin{aligned}
&\bbE\big[\big\|\mcQ_t(x, {\sf\Pi}^{n;\varepsilon,\infty}_x\tau_i)\big\|_{L_x^q(w_{c(m_{\bf B}+2)})}^q\big]   \\
&\lesssim
|E_i^n|^q + t^{\frac{r_{\varepsilon,\infty}(\tau_i)}\ell q} \bigg( \bbE\big[\|{\sf\Pi}^{n;\varepsilon,\infty} \res {\bf B}_{i-1}\|_{w_c} \, \|{\sf\Gamma}^{n;\varepsilon,\infty} \res \tau_i\|_{w_{cm_{\bf B}}}\big]^q + \bbE\bigg[\sup_{\|h\|_{H}\le1} \|{\sf\Pi}^{n,h;\varepsilon,\infty} \res \dt{\bf B}_i\|_{w_c}^q\bigg]\bigg).
\end{aligned} \end{align}
To trade the $L_x^q(w_{c(m_{\bf B}+2)})$ norm for an $L_x^\infty(w_{c(m_{\bf B}+2)})$ norm in the above estimate we use an argument that is reminiscent of the proof of Besov embedding. Here we need a slight change of the parameter $\varepsilon$ to make that trade work.

\ssk

\begin{lem} \label{lem:Binega}
Let $r_{0,\infty}(\tau_i) \le 0$. For any $\varepsilon\in(0,\eps_0)$, any $c>|\mfs|$, and $q\in[1,\infty)$, there exists a positive constant $C$ which is independent of $n$, one has
\begin{align*} \begin{aligned}
\bbE\big[\|{\sf \Pi}^{n;\eps,\infty} &\res \tau_i\|_{w_{c(m_{\bf B}+2)}}^q\big]   \\
&\le C\bigg(|E_i^n|^q + \bbE\big[\|{\sf\Pi}^{n;\frac\eps2,\infty} \res {\bf B}_{i-1}\|_{w_c}^q \, \|{\sf\Gamma}^{n;\frac{\eps}2,\infty} \res \tau_i\|_{w_{cm_{\bf B}}}^q\big] 
+ \bbE\bigg[\sup_{\|h\|_{H}\le1} \|{\sf\Pi}^{n,h;\frac{\eps}2,\infty}\res\dt{\bf B}_i\|_{w_c}^q\bigg]\bigg).
\end{aligned} \end{align*}
\end{lem}

\ssk

Lemma \ref{lem:Biposi} and Lemma \ref{lem:Binega} provide an $L^q(\Omega)$ bound of $\|{\sf \Pi}^{n;\eps,\infty} \res \tau_i \|_{w_{c(m_{\bf B}+2)}}$ in terms of the moment of $\Vert {\sf M}^{n;\varepsilon/2,p}\Vert_{\textbf{\textsf{M}}(\scW_{i,\varepsilon/2,\infty})_{w_c}}$, which is a part of $\textsf{\textbf{bd}}(\scW,i,\infty)$.

%%------------------------------------------------------------------------------------------------------------------------------------------------------------------------------------------------------------------------------------------------------------------------%%
\subsection{Step 2: From $\textsf{bd}(\scW,i,p)$ and $\textsf{bd}(\scV,i)$ to the $\sf g$-part of $\textsf{bd}(\scW,i+1,p) \boldmath .$  \hspace{0.3cm}}
\label{sec:ind2}
%%------------------------------------------------------------------------------------------------------------------------------------------------------------------------------------------------------------------------------------------------------------------------%%

In the remaining steps of the proof of \textsf{\textbf{bd}}, all the inequalities are consequences of some analytic results stated in Appendix \ref{SectionAppendix}. It should then be noted that {\bf all the proportional constants $C$ below are deterministic and do not depend on $n$ nor $\omega$.} Recall again that Step 1 is the only place where the spectral gap assumption on the law of the noise is used.

\medskip

We fix some parameters $\varepsilon\in(0,\varepsilon_0)$ and $p\in[2,\infty]$ in this step. The work is to bound the $\sf g$-part of
\begin{equation}\label{*step2:tobebounded}
\|{\sf M}^{n;\varepsilon,p}\|_{\textbf{\textsf{M}}(\scW_{i+1,\varepsilon,p})_{w_c}} 
= \|{\sf\Pi}^{n;\varepsilon,p} \res {\bf B}_{i}\cup\dt{\bf B}_{i+1}\|_{w_c}  
+ \|{\sf g}^{n;\varepsilon,p} \res {\bf W}_{i,\varepsilon,p}^+\|_{w_c}. 
\end{equation}
It is sufficient for this purpose to show an $\omega$-wise bounds on
$$
\sup_{n\in\bbN}\Vert {\sf g}^{n;\varepsilon,p} \res \mcI_k^{\bbK}(\tau)\Vert_{w_c}
$$
for any $\tau\in{\bf B}_i\cup\dt{\bf B}_i$ and $|k|_\mfs<r_{\eps,p}(\tau)+\beta_0$ in terms of the assumed boundedness results \eqref{*eq:tobebounded1} on $\scV_{i,\varepsilon}$ and \eqref{*eq:tobebounded2} on $\scW_{i,\varepsilon,\infty}$. 

\ssk

We need temporally some notations. Denote by $\overline{W}_{i-1,\varepsilon,p}$ the linear subspace of $W_{i-1,\eps,p}^+$ spanned by the symbols
$$
\overline{\bf W}_{i-1,\varepsilon,p}\defeq\{X^k\}_{k\in\bbN^d}\cup\{\mcI_0^{\bbK}(\tau)\}_{\tau\in({\bf B}_{i-1}\cup\dt{\bf B}_{i-1})\setminus\{X^l\}_l,\, r_{\varepsilon,p}(\tau)+\beta_0>0}.
$$
Since $\Delta_{\eps,p}^+(\overline{W}_{i-1,\eps,p})\subset\overline{W}_{i-1,\eps,p}\otimes W_{i-1,\eps,p}^+$, the pair
$$
\overline{\scW}_{i-1,\eps,p}\defeq\big((\overline{W}_{i-1,\eps,p},\Delta_{\eps,p}^+),(W_{i-1,\varepsilon,p}^+,\Delta_{\eps,p}^+)\big)
$$
is an concrete regularity structure, and the pair of operators $\overline{\sf M}^{n;\varepsilon,p}=(\overline{\sf\Pi}^{n;\varepsilon,p},\overline{\sf\Gamma}^{n;\varepsilon,p})$ given by
$$
(\overline{\sf\Pi}_x^{n;\varepsilon,p}\mu)(x)={\sf g}_{yx}^{n;\varepsilon,p}(\mu),\qquad
\overline{\sf\Gamma}_{yx}^{n;\varepsilon,p} \defeq \big(\id\otimes\mathop{{\sf g}_{yx}^{n;\varepsilon,p}}\big) \Delta^+_{\varepsilon,p}
$$
is a model on $\overline{\scW}_{i-1,\eps,p}$. Their norms are bounded by
\begin{align*}
\|\overline{\sf\Pi}^{n;\eps,p}\res\overline{\bf W}_{i-1,\eps,p}\|_{w_c}
&\lesssim\|{\sf g}^{n;\eps,p}\res{\bf W}_{i-1,\eps,p}^+\|_{w_c},\\
\|\overline{\sf\Gamma}^{n;\eps,p}\res\overline{\bf W}_{i-1,\eps,p}\|_{w_{cm_{\bf B}}}
&\lesssim\big(1+\|{\sf g}^{n;\eps,p}\res{\bf W}_{i-1,\eps,p}^+\|_{w_c}\big)^{m_{\bf B}}.
\end{align*}
We further define the linear map $\overline{\mcI}^{\varepsilon,p}:U_{i-1}\to \overline{W}_{i-1,\varepsilon,p}$ by setting
$$
\overline{\mcI}^{\varepsilon,p}(\sigma) \defeq 
\begin{cases}
\mcI_0^{\bbK}(\sigma), & \textrm{ if } r_{\varepsilon,p}(\sigma)+\beta_0>0,   \\
0, & \textrm{ if } r_{\varepsilon,p}(\sigma)+\beta_0\le0,
\end{cases}
$$
which is an abstract integration map of order $\beta_0$ -- see Definition \ref{defTbarT}. Then in the proof of the following lemma, it turns out that the models ${\sf M}^{n;\varepsilon,p}$ and $\overline{\sf M}^{n;\varepsilon,p}$ are compatible for $\overline{\mcI}^{\varepsilon,p}$. 

\ssk

Set $U_{i-1}\defeq\spa({\bf B}_{i-1}\cup\dt{\bf B}_{i-1})$ and define, as in Lemma \ref{*argumentforreconst}, a modelled distribution $f_\tau^{n;\varepsilon,p}$ setting
\begin{align*}
f_\tau^{n;\varepsilon,p}(x) = \big(\id\otimes\mathop{{\sf g}^{n;\varepsilon,p}_x}\big) \Delta_{\varepsilon,p}\tau-\tau
\in
D^{(r_{\varepsilon,p}(\tau),i_p(\tau))}(U_{i-1};{\sf\Gamma}^{n;\varepsilon,p})_{w_{cm_{\bf B}}},
\end{align*}
We note that ${\sf\Pi}^n\tau$ is a reconstruction of $f_\tau^{n;\varepsilon,p}$, since the appropriate bound on ${\sf\Pi}_x^{n;\varepsilon,p}\tau={\sf\Pi}^n\tau-{\sf\Pi}_x^{n;\varepsilon,p}(f_\tau^{n;\varepsilon,p})$ is a part of $\textsf{\textbf{bd}}(\scW,i,p)$ and $\textsf{\textbf{bd}}(\scV,i)$. As in the usual setting of regularity structures, one needs some additional terms to turn $\overline{\mcI}^{\varepsilon,p}(f_\tau^{n;\varepsilon,p})$ into a modelled distribution -- here in a regularity-integrability setting. We define the modelled distribution $\overline{\mcK}^{n;\varepsilon,p}(f_\tau^{n;\varepsilon,p})$ by the formula \eqref{def:MS} with the model ${\sf M}^{n;\varepsilon,p}$.

\ssk

\begin{lem}\label{lem:abstint}
One has the expansion
$$
\overline{\mcK}^{n;\varepsilon,p}\big(f_\tau^{n;\varepsilon,p}\big)(y) - \overline{\sf \Gamma}_{yx}^{n;\varepsilon,p}\Big(\overline{\mcK}^{n;\varepsilon,p}\big(f_\tau^{n;\varepsilon,p}\big)(x)\Big)
=\sum_{|k|_\mfs<r_{\eps,p}(\tau)+\beta_0}\frac{X^k}{k!} \, {\sf g}_{yx}^{n;\varepsilon,p}(\mcI_{k}^{\bbK}\tau).
$$
\end{lem}

\ssk

The multilevel Schauder estimate from Theorem \ref{thm:MS} then gives an upper bound on ${\sf g}_{yx}^{n;\varepsilon,p}(\mcI_{k}^{\bbK}\tau)$ in terms of the norm of $f_\tau^{n;\varepsilon,p}$. They take the following form, proved in Section {{\ref{SectionProofLemmaSection33}}}.

\ssk

\begin{lem} \label{lem:Pitog}
For any $\varepsilon\in(0,\varepsilon_0)$, $p\in[2,\infty]$ and $c>0$ there exists a positive constant $C$ which is independent of $n$, $\omega$, and $h\in H$ with $\|h\|_H\le1$, such that one has
\begin{align*}
\big\Vert {\sf g}^{n;\varepsilon,p} \res \mcI_{k}^{\bbK}(\tau) \big\Vert_{w_{c(2m_{\bf B}+1)}} 
&\le \begin{cases}
C\big(1 + \Vert {\sf M}^{n;\varepsilon,\infty}\Vert_{\textbf{\textsf{M}}(\scV_{i,\varepsilon})_{w_c}}\big)^{2m_{\bf B}+1},%+ \|{\sf\Pi}^{n;\varepsilon,\infty} \hspace{-0.03cm}:\hspace{-0.03cm} \tau\|_{w_c},  
\quad 
&\textrm{if } \tau\in{\bf B}_i,   \\
C\big(1 + \Vert {\sf M}^{n;\varepsilon,p}\Vert_{\textbf{\textsf{M}}(\scW_{i,\varepsilon,p})_{w_c}}\big)^{2m_{\bf B}+1},  &\textrm{if } \tau\in\dt{\bf B}_i.
\end{cases} 
\end{align*}
\end{lem}

\ssk

Re-inserting the $h$ in the notations ${\sf M}^{n,h}$ and ${\sf \Pi}^{n,h;\varepsilon,\infty}$ for ${\sf M}^n$ and ${\sf \Pi}^{n;\varepsilon,\infty}$ to emphasize the dependence on $h$ of these objects, the $n$-uniform control of 
$$
\bbE\bigg[ \sup_{\Vert h\Vert_{H}\leq 1}  \big\|{\sf g}^{n,h;\varepsilon,p} \res {\bf W}_{i,\varepsilon,p}^+ \big\|_{w_c}^q \bigg]
$$
follows from Lemma \ref{lem:Pitog} and the joint assumptions $\textbf{\textsf{bd}}(\scW,i,p)$ and $\textbf{\textsf{bd}}(\scV,i)$.

%\bigskip

%%-----------------------------------------------------------------------------------------------------------------------------------------------------------------------------------------------------------------------------------------------------------------------%%
\subsection{Step 3: From $\big\{\textsf{bd}(\scW,i,p)\big\}_{2\le p\le\infty}$ and $\textsf{bd}(\scV,i)$ to the $\sf \Pi$-part of $\big\{\textsf{bd}(\scW,i+1,p)\big\}_{2\le p\le\infty} \boldmath .$  \hspace{0.3cm}}
\label{sec:ind3}
%%-----------------------------------------------------------------------------------------------------------------------------------------------------------------------------------------------------------------------------------------------------------------------%%

It remains to show the bounds of the $\sf\Pi$-part of \eqref{*step2:tobebounded}. For this purpose, it is sufficient to show the bound on
$$
\sup_{n\in\bbN}\Vert {\sf \Pi}^{n;\varepsilon,p} \hspace{-0.03cm}:\hspace{-0.03cm} \tau\Vert_{w_c}
$$
for any $\tau\in\dt{\bf B}_{i+1}$ in terms of the assumed boundedness results \eqref{*eq:tobebounded1} on $\scV_{i,\varepsilon}$ and \eqref{*eq:tobebounded2} on $\scW_{i,\varepsilon,p}$. We first establish the result for $p=2$ and next extend it into all $p\in[2,\infty]$ by using the formula \eqref{EqSyntheticDecompositionFormulaDerivativePi}.

\medskip

-- Let $\tau$ be of the form $D_e\sigma$ with some $\sigma\in{\bf B}_{i+1}$ and $e\in{\fr t}^{-1}(\bbOmega)$. Since $\sigma\neq\ocircle$ the assumption of Theorem \ref{thm:main} gives here 
$$
r_{0,2}(\tau)=r_{0,\infty}(\tau)+\frac{|\mfs|}2=r_{0,\infty}(\sigma)+\frac{|\mfs|}2>0.
$$
Thus, as in Lemma \ref{lem:Biposi}, one gets the following statement from an application of the reconstruction theorem to the modelled distribution
\begin{align*}
f_\tau^{n;\varepsilon,2}(x) =
\big(\id\otimes\mathop{{\sf g}^{n;\varepsilon,2}_x}\big) \Delta_{\varepsilon,2}\tau-\tau\in D^{(r_{\varepsilon,2}(\tau),2)}\big(\spa({\bf B}_i\cup\dt{\bf B}_i);{\sf\Gamma}^{n;\varepsilon,2}\big)_{w_c}.
\end{align*}

\ssk

\begin{lem} \label{LemContinuousExtension}
For any $\varepsilon\in(0,\varepsilon_0)$ and $c>0$, there exists a positive constant $C$ which is independent of $n$, $\omega$, and $h\in H$ with $\Vert h\Vert_H\leq 1$ one has
\begin{equation*} %\label{EqControlPEqual2}
\|{\sf\Pi}^{n;\varepsilon,2} \res \tau\|_{w_{c(m_{\bf B}+1)}} 
\leq C \, \|{\sf\Pi}^{n;\varepsilon,2} \res {\bf B}_i\cup\dt{\bf B}_i\|_{w_c} \,  \|{\sf\Gamma}^{n;\varepsilon,2}\res {\bf B}_i\cup\dt{\bf B}_{i+1}\|_{w_{cm_{\bf B}}}.
\end{equation*}
\end{lem}

\ssk

The proof of the preceding lemma is almost the same as that of Lemma \ref{lem:Biposi} and left to the reader. Recall from \eqref{EqControlGamma} that the bound on $\Vert{\sf\Gamma}^{n;\varepsilon,2}\res {\bf B}_i\cup\dt{\bf B}_{i+1}\Vert_{w_{cm_{\bf B}}}$ follows from the assumption $\textbf{\textsf{bd}}(\scW,i,2)$ and the bounds of $\|{\sf g}^{n;\varepsilon,2} \res {\bf W}_{i,\varepsilon,2}^+\|_{w_c}$ obtained in the last section. Thus we have completed the proof of $\textbf{\textsf{bd}}(\scW,i+1,2)$.

\ssk

-- We now want to infer from the bound of Lemma \ref{LemContinuousExtension} a similar bound on ${\sf\Pi}_x^{n;\varepsilon,p}\tau$. We use for that purpose the algebraic formula \eqref{EqSyntheticDecompositionFormulaDerivativePi}
\begin{equation*}
{\sf\Pi}_x^{n;\varepsilon,p}(\tau) = {\sf\Pi}_x^{n;\varepsilon,2}(\tau) + \big({\sf\Pi}^{n;\varepsilon,p}_x\otimes\lambda_x^{n;\varepsilon,p}\big) \Delta_{\varepsilon,2}(\tau).
\end{equation*}
We can estimate as follows the size of the $\lambda_x^{n;\varepsilon,p}$ terms. Recall from Section \ref{SubsectionDifferentialSectors} the definitions of the exponent $p_\varepsilon(\mu)$ and the floor function $\lfloor p\rfloor_{{\sf I}_\varepsilon}$.

\ssk

\begin{lem}\label{lem:Lpfx}
For any $\mu\in{\bf W}_{i,\eps,p}^+$, $p\in[2,\infty]$ and $p_\varepsilon(\mu) > q \ge\lfloor p_\varepsilon(\mu)\rfloor_{{\sf I}_\varepsilon}$ one has
$$
\|\lambda_x^{n;\varepsilon,p}(\mu)\|_{L_x^q(w_{c(m_{\bf B}+1)})}\lesssim\|{\sf\Pi}^{n;\varepsilon,q} \res{\bf B}_{i-1}\cup\dt{\bf B}_i\|_{w_c}\big(1+\|{\sf\Gamma}^{n;\varepsilon,q} \res\dt{\bf B}_i \|_{w_c}\big).
$$
\end{lem}

\ssk

The proof is given in Section {{\ref{SectionProofLemmas34}}}. Formula \eqref{EqSyntheticDecompositionFormulaDerivativePi} then leads to the following estimate.

\ssk

\begin{lem} \label{lem:goal}
Let $\tau\in\dt{\bf B}_{i+1}$, $p\in[2,\infty]$ and $\varepsilon\in(0,\eps_0)\setminus {\sf J}_p$ be given. Then, for each $q\in\{2\}\cup {\sf I}_\eps$, one can choose $\eps'=\eps'(q)\in(0,\eps_0)\setminus {\sf J}_q$ and a positive constant $C$, both independent of $n$, $\omega$ and $h\in H$ with $\Vert h\Vert_H\leq 1$, such that one has for any $c>0$
$$
\big\| {\sf\Pi}^{n;\eps,p}\res\tau \big\|_{w_{c(m_{\bf B}+2)}}  
\le C\Bigg(1+\sum_{q\in \{2\}\cup {\sf I}_{\eps}} \big\Vert {\sf M}^{n;\eps',q} \big\Vert_{\textbf{\textsf{M}}(\scW_{i,\eps',q})_{w_c}} + \big\Vert {\sf M}^{n;\eps,\infty} \big\Vert_{\textbf{\textsf{M}}(\scV_{i,\eps})_{w_c}}\Bigg)^{m_{\bf B}+2}.
$$
\end{lem}

\ssk

Thus we see from Lemma \ref{lem:goal} that the $n$-uniform control of 
$$
\bbE\bigg[\sup_{\Vert h\Vert_{H}\leq 1} \big\| {\sf\Pi}^{n,h;\varepsilon,p}\hspace{-0.03cm}:\hspace{-0.03cm}\tau \big\|_{w_{c(m_{\bf B}+2)}}^q\bigg] 
$$
is given by the joint assumptions $\big\{\textbf{\textsf{bd}}(\scW,i,p)\big\}_{2\leq p\leq \infty}$ and $\textbf{\textsf{bd}}(\scV,i)$.

\ssk

A finite number of iterations of Steps 1-3 starting from the initial case of Section {{\ref{sec:ind0}}} proves the $n$-uniform bounds of ${\sf M}^n$ in $L^q\big(\Omega,\bbP ; \textbf{\textsf{M}}(\scW_{\varepsilon,p})_{w_c}\big)$ for any $p\in[2,\infty]$, $q\in[1,\infty)$, any small $\varepsilon>0$ and any large $c>|\mfs|$.

\medskip

\noindent \textbf{Remark --} \textit{For readers familiar with the work \cite{LOTT} of Linares, Otto, Tempelmayr \& Tsatsoulis our Step 2 is somewhat the equivalent of their algebraic \& `three point' arguments. Out Step 3(a) is the equivalent of their `Reconstruction III' step and our Step 3(b) is the equivalent of their `Averaging' step.}

%\bigskip

%%--------------------------------------------------------------------------------------------%%
\subsection{From uniform boundedness to convergence results$\boldmath .$  \hspace{0.3cm}}
\label{sec:ind4}
%%--------------------------------------------------------------------------------------------%%

Since the models ${\sf M}^n$ stay in a bounded set of $L^q\big(\Omega,\bbP ; \textbf{\textsf{M}}(\scW_{\varepsilon,p})_{w_c}\big)$ for any $p,q,\varepsilon,c$ we can use the local Lipschitz estimates satisfied by the reconstruction operator and the $\overline{\mcK}^{n;\varepsilon,p}$ maps -- see Theorem \ref{thm:besovreconstruction} and Theorem \ref{thm:MS} in Appendix \ref{SectionAppendix}, to prove the Cauchy property
\begin{align}\label{EqCauchyMn}
\lim_{n,m\to\infty}\bbE\bigg[\sup_{\Vert h\Vert_H\leq 1} \Vert{\sf M}^{n;\varepsilon,p} \hspace{-0.03cm}:\hspace{-0.03cm} {\sf M}^{m;\varepsilon,p}\Vert_{\textbf{\textsf{M}}(\scW_{\varepsilon,p})_{w_c}}^q\bigg]=0
\end{align}
for any $p,q,\varepsilon,c$. Since $\textbf{\textsf{M}}(\scW_{\varepsilon,p})_{w_c}$ is complete, this implies the convergence in $L^q\big(\Omega,\bbP ; \textbf{\textsf{M}}(\scW_{\varepsilon,p})_{w_c}\big)$ of $({\sf M}^{n;\varepsilon,p})_{n\geq 0}$. We can prove the Cauchy property \eqref{EqCauchyMn} starting from the convergence result of the initial case corresponding to Section \ref{sec:ind0} and following the same induction steps as above. We only collect below the statements corresponding to the different lemmas from sections \ref{sec:ind1}, \ref{sec:ind2} and \ref{sec:ind3}. A statement corresponding to Lemma $k$ in one of the previous sections is numbered here Lemma $k$'. We denote by $Q_n$ some non-negative random variables which depend polynomially on the quantities $\sup_{\Vert h\Vert_H\leq 1} \Vert {\sf M}^{n;\varepsilon,p} \Vert_{\textbf{\textsf{M}}(\scW_{\varepsilon,p})_{w_c}}$ where the parameters $\varepsilon,p,c$ run over an $n$-independent finite set. We know from the results of Section \ref{sec:ind1} to \ref{sec:ind3} that
$$
\sup_{n\in\bbN}\bbE\big[Q_n^q\big]<\infty
$$
for all $q\in[1,\infty)$.

\ssk

%%%---------------------------------------------------------------------------------------------------------------------------------------------------------------------------------------------------------------------------------------------%%%
\textit{\S0 The initial case $\{\textsf{\textbf{cv}}(\scW,1,p)\}_{2\leq p\leq \infty}\boldmath .$}   \hspace{0.1cm}
%%%---------------------------------------------------------------------------------------------------------------------------------------------------------------------------------------------------------------------------------------------%%%
By the same argument as Section \ref{sec:ind0} where $h_n$ is replaced by $h_n-h_m$, the convergence for the initial case follows from the convergence $\lim_{m,n\to\infty}\|h_n-h_m\|_H=0$, which is uniform over $\|h\|_H\le1$.

\medskip

%%%---------------------------------------------------------------------------------------------------------------------------------------------------------------------------------------------------------------------------------------------%%%
\textit{\S1 Step 1: From $\textsf{\textbf{cv}}(\scW,i,\infty)$ to $\textsf{\textbf{cv}}(\scV,i)\boldmath .$}   \hspace{0.1cm}
%%%---------------------------------------------------------------------------------------------------------------------------------------------------------------------------------------------------------------------------------------------%%%
We assume the convergence result $\textsf{\textbf{cv}}(\scW,i,\infty)$ and prove the local Lipschitz estimate of 
$$
\sup_{n\in\bbN}\bbE\Big[\Vert {\sf\Pi}^{n;\varepsilon,\infty},{\sf\Pi}^{m;\varepsilon,\infty} \res \tau_i\Vert_{w_c}^q\Big]
$$
for $n,m\in\bbN$. We use different inequalities depending on the sign of $r_{0,\infty}(\tau_i)$.

\ssk

\begin{lem11}
Let $r_{0,\infty}(\tau_i)>0$. For any $\eps\in[0,\eps_0)$ and $c>0$, one has the $\omega$-wise estimate
$$
\|{\sf\Pi}^{n;\varepsilon,\infty},{\sf\Pi}^{m;\varepsilon,\infty}\res\tau_i\|_{w_{c(m_{\bf B}+2)}}
\le (Q_{n}+Q_{m})\|{\sf M}^{n;\varepsilon,\infty} \res {\sf M}^{m;\varepsilon,\infty}\|_{\textbf{\textsf{M}}(\scV_{i-1,\varepsilon,p})_{w_c}}
$$
for any $n,m\in\bbN$.
\end{lem11}

\ssk

\begin{lem12}
Let $r_{0,\infty}(\tau_i) \leq 0$. For any $\eps\in[0,\eps_0)$ and $c>|\mfs|$, and $q\in[1,\infty)$, there exists a positive constant $C$ which is independent of $n,m\in\bbN$ and one has the moment estimate
\begin{align*}
&\bbE\Big[ \big\| {\sf\Pi}^{n;\varepsilon,\infty},{\sf\Pi}^{m;\varepsilon,\infty}\res\tau_i \big\|_{w_{c(m_{\bf B}+2)}}^q \Big] \leq   \\
&C \bigg(|E_i^n-E_i^m|^q
+ \bbE\Big[ \big\| {\sf M}^{n;\frac\varepsilon2,\infty} \res {\sf M}^{m;\frac{\varepsilon}2,\infty} \big\|_{\textbf{\textsf{M}}(\scV_{i-1,\frac\varepsilon2,p})_{w_c}}^{2q}\Big]^{\frac{1}{2}}
+ \bbE\bigg[\sup_{\|h\|_{H}\le1} \big\|{\sf\Pi}^{n,h;\frac\varepsilon2,\infty},{\sf\Pi}^{m,h;\frac\varepsilon2,\infty}\res \dt{\bf B}_i \big\|_{w_c}^q\bigg]
\bigg)
\end{align*}
\end{lem12}

\medskip

%%%---------------------------------------------------------------------------------------------------------------------------------------------------------------------------------------------------------------------------------------------%%%
\textit{\S2 Step 2: From $\textsf{\textbf{cv}}(\scW,i,p)$ and $\textsf{\textbf{cv}}(\scV,i)$ to the $\sf g$-part of $\textsf{\textbf{cv}}(\scW,i+1,p) \boldmath .$}   \hspace{0.1cm} 
%%%---------------------------------------------------------------------------------------------------------------------------------------------------------------------------------------------------------------------------------------------%%%
We assume the convergence results $\textsf{\textbf{cv}}(\scW,i,p)$ and $\textsf{\textbf{cv}}(\scV,i)$ and prove the ($\omega$-wise) local Lipschitz estimate of 
$$
\big\Vert {\sf g}^{n;\varepsilon,p}, {\sf g}^{m;\varepsilon,p} \res \mcI_k^{\bbK}(\tau) \big\Vert_{w_c}
$$
for $\tau\in{\bf B}_i\cup\dt{\bf B}_i$ and $k\in\bbN^d$.

\ssk

\begin{lem2}
For any $\varepsilon\in(0,\varepsilon_0)$, $p\in[2,\infty]$, and $c>0$, one has the $\omega$-wise estimate
\begin{align*}
\big\Vert {\sf g}^{n;\varepsilon,p}, {\sf g}^{m;\varepsilon,p} \res \mcI_k^{\bbK}(\tau) \big\Vert_{w_c}
&\le \begin{cases}
(Q_{n}+Q_{m}) \, \big\| {\sf M}^{n;\varepsilon,\infty} \res {\sf M}^{m;\varepsilon,\infty} \big\|_{\textbf{\textsf{M}}(\scV_{i,\varepsilon})_{w_c}},
\quad 
&\textrm{if } \tau\in{\bf B}_i,   \\
(Q_{n}+Q_{m}) \, \big\| {\sf M}^{n;\varepsilon,p} \res {\sf M}^{m;\varepsilon,p} \big\|_{\textbf{\textsf{M}}(\scW_{i,\varepsilon,p})_{w_c}},&\textrm{if } \tau\in\dt{\bf B}_i.
\end{cases} 
\end{align*}
\end{lem2}

\medskip

%%%---------------------------------------------------------------------------------------------------------------------------------------------------------------------------------------------------------------------------------------------%%%
\textit{\S3 Step 3: From $\{\textsf{\textbf{cv}}(\scW,i,p)\}_{2\le p\le\infty}$ and $\textsf{\textbf{cv}}(\scV,i)$ to the $\sf \Pi$-part of $\{\textsf{\textbf{cv}}(\scW,i+1,p)\}_{2\le p\le\infty} \boldmath .$}   \hspace{0.1cm}
%%%---------------------------------------------------------------------------------------------------------------------------------------------------------------------------------------------------------------------------------------------%%%
We assume the convergence results $\{\textsf{\textbf{cv}}(\scW,i,p)\}_{2\le p\le\infty}$ and $\textsf{\textbf{cv}}(\scV,i)$ and prove the ($\omega$-wise) local Lipschitz estimate of 
$$
\big\Vert {\sf \Pi}^{n;\varepsilon,p}, {\sf \Pi}^{m;\varepsilon,p} \hspace{-0.03cm}:\hspace{-0.03cm} \tau\big\Vert_{w_c}
$$
for $\tau\in\dt{\bf B}_{i+1}$.

\ssk

\begin{lem3}
For any $p\in[2,\infty]$ and $\eps\in(0,\eps_0)\setminus {\sf J}_p$, with the same choice of $\eps'=\eps(q)\in(0,\eps_0)\setminus {\sf J}_q$ for each $q\in\{2\}\cup {\sf I}_\eps$ as in Lemma \ref{lem:goal}, one has
\begin{equation*} \begin{split}
\big\Vert {\sf \Pi}^{n;\varepsilon,p}, {\sf \Pi}^{m;\varepsilon,p} &\res \tau \big\Vert_{w_c}   \\
&\leq (Q_n+Q_m)
\Bigg(
\sum_{q\in \{2\}\cup {\sf I}_\eps}
\big\| {\sf M}^{n;\varepsilon',q} \res {\sf M}^{m;\varepsilon',q} \big\|_{\textbf{\textsf{M}}(\scW_{i,\varepsilon',q})_{w_c}}
+
\big\| {\sf M}^{n;\varepsilon,\infty} \res {\sf M}^{m;\varepsilon,\infty} \big\|_{\textbf{\textsf{M}}(\scV_{i,\varepsilon})_{w_c}}
\Bigg).
\end{split} \end{equation*}
\end{lem3}

%%--------------------------------------------------------------------------------------------%%
\subsection{Universality of the limit model$\boldmath .$  \hspace{0.3cm}}
\label{sec:ind5}
%%--------------------------------------------------------------------------------------------%%

Finally we show that the limit model ${\sf M}^{\eps,p}=\lim_{n\to\infty}{\sf M}^{n;\eps,p}$ depends only on the finite family $\{\lim_{n\to\infty}E_\tau^{n,R_n}\}_{\tau\in{\bf B},\,r_{0,\infty}(\tau)\le0}$.

We write ${\sf M}^{n;\eps,p}(\varrho,R)$ to emphasize the dependence on $\varrho$ and $R$. For another choices of mollifiers $\widetilde{\varrho}=(\widetilde{\varrho}_n)_{n\geq 0}$ and $\widetilde{R}=(\widetilde{R}_n)_{n\geq 0}$ such that $\lim_{n\to\infty}E_\tau^{n,\widetilde{R}_n}$, by using some argument similar to the arguments used in Section \ref{sec:ind4}, we have that
\begin{equation}
\lim_{n\to\infty}\bbE\bigg[\sup_{\Vert h\Vert_H\leq 1} \Vert{\sf M}^{n;\varepsilon,p}(\varrho,R) \res {\sf M}^{n;\varepsilon,p}(\widetilde{\varrho},\widetilde{R})\Vert_{\textbf{\textsf{M}}(\scW_{\varepsilon,p})_{w_c}}^q\bigg]=0
\end{equation}
whenever
$$
\lim_{n\to\infty} \sup_{\|h\|_H\le1} \big\| h*\varrho_n-h*\widetilde{\varrho}_n \big\|_H=0, \qquad
\lim_{n\to\infty} \max_{\tau\in{\bf B},\,r_{0,\infty}(\tau)\le0} \big| E_\tau^{n,R_n}-E_\tau^{n,\widetilde{R}_n} \big| = 0.
$$
Thus $\lim_{n\to\infty}{\sf M}^{n;\eps,p}(\varrho,R)=\lim_{n\to\infty}{\sf M}^{n;\eps,p}(\widetilde{\varrho},\widetilde{R})$.
	
%\ssk

%%--------------------------------------------------------------------------------------------%%
\subsection{Proof of Corollary \ref{cor:BPHZ}$\boldmath .$  \hspace{0.3cm}}
\label{sec:ind6}
%%--------------------------------------------------------------------------------------------%%

We show that the renormalization procedure $(R_n)_{n\geq 0}$ associated with the BPHZ models satisfies the assumption of Theorem \ref{thm:main}.
First we show the finiteness of
$$
\sup_{n}|E_\tau^{n,R_n}|=\sup_n\big|\bbE\big[\mcQ_1(0, {\sf\Pi}_0^{n;\eps}\tau)\big]\big|.
$$
for any $\tau\in{\bf B}$. Note that this finiteness is used only in Step 1, from $\textsf{\textbf{bd}}(\scW,i,\infty)$ to $\textsf{\textbf{bd}}(\scV,i)$. Therefore, we can show the result by an induction as follows.
\begin{enumerate} \setlength{\itemsep}{0.1cm}
\renewcommand{\labelenumi}{\textrm{\bf (\arabic{enumi})}}
	\item For the initial case $\tau=\ocircle$, since ${\sf\Pi}_0^{n;\eps}(\tau_1) = {\sf\Pi}^{n}(\tau_1) = \xi_n$, we have $E_{\tau_1}^{n,R_n}=0$.
	
	\item Assuming $\textsf{\textbf{bd}}(\scW,i,\infty)$ and $\sup_n|E_{\tau_i}^{n,R_n}|<\infty$, we have $\textsf{\textbf{bd}}(\scV,i)$. 
Moreover, by the inequality \eqref{*eq:resultofSGineq}, we have
\begin{equation}\label{*eq:cor:proofBPHZ}
\sup_n\bbE\big[\big|\mcQ_1(0, {\sf\Pi}^{n;\varepsilon}_0(\tau_i))\big|^q\big] <\infty
\end{equation}
for any $q\ge1$.
After that, we have $\{\textsf{\textbf{bd}}(\scW,i+1,p)\}_{2\leq p\leq \infty}$ following Step 2 and Step 3. 

	\item If $r_{0,\infty}(\tau_{i+1})>0$ the bound of ${\sf\Pi}_x^{n;\eps}(\tau_{i+1})$ is obtained from $\textsf{\textbf{bd}}(\scV,i)$ as a consequence of Lemma \ref{lem:Biposi}. We have in particular $\sup_n|E_{\tau_{i+1}}^{n,R_n}|<\infty$.

	\item If $r_{0,\infty}(\tau_{i+1})\le0$, since $({\sf\Pi}_0^{n;\eps}\otimes{\sf g}_0^{n;\eps})\Delta_\eps\tau_{i+1}={\sf\Pi}^n\tau_{i+1}$, we can write in synthetic form
$$
{\sf\Pi}_0^{n;\eps}(\tau_{i+1}) = {\sf\Pi}^n(\tau_{i+1}) - \sum_{\sigma_1,\sigma_2}{\sf g}_0^{n;\eps}(\sigma_2) \, {\sf\Pi}_0^{n;\eps}(\sigma_1),
$$
where $\sigma_1\in V_i$ and $\sigma_2\in V_{i,\eps}^+$. By applying the operator $\mcQ_1$ and taking the expectation, we have
$$
\bbE\big[\mcQ_1\big(0, {\sf\Pi}_0^{n;\eps}(\tau_{i+1})\big)\big]=-\sum_{\sigma_1,\sigma_2}\bbE\big[{\sf g}_0^{n;\eps}(\sigma_2) \, \mcQ_1\big(0,{\sf\Pi}_0^{n;\eps}(\sigma_1)\big)\big].
$$
from the definition of BPHZ model.
By the properties \eqref{*eq:cor:proofBPHZ} which is satisfied by all $\tau_j$ ($j\le i$) by induction and by $\textsf{\textbf{bd}}(\scW,i+1,\infty)$, the right hand side has an $n$-uniform bound.
Thus $\sup_n|E_{\tau_{i+1}}^{n,R_n}|<\infty$.
\end{enumerate}

\ssk

It remains to show the convergence of $(E_\tau^{n,R_n})_{n\geq 0}$ and that the limit $\lim_{n\to\infty}E_\tau^{n,R_n}$ does not depend on the choice of $(\varrho_n)_{n\geq 0}$. It is sufficient for that purpose to modify the above arguments as in Section \ref{sec:ind4} and Section \ref{sec:ind5}.

\medskip

%------------------------------------%
\section{Proofs of the lemmas}
\label{SectionLemmas}
%------------------------------------%

We give in this section the proofs of the lemmas used in the proof of Theorem \ref{thm:main} detailed in Section \ref{SectionMechanics}.

%\ssk

%%--------------------------------------------------------------------------------------%%
\subsection{Trees conforming to rules$\boldmath .$ \hspace{0.03cm}}
\label{ProofStronglyConform}
%%--------------------------------------------------------------------------------------%%

We briefly recall some notions from Section 5 of \cite{BHZ}. Denote by $\varrho_\tau$ the root of $\tau\in{\bf T}$.

\ssk

\begin{defn}
Let $\mcN$ be the collection of all finite multisets $\{\!\{r_1,\dots,r_n\}\!\}$ whose elements $r_i$ are elements of $\{\bbOmega,\bbK\}\times\bbN^d$. An empty multiset $\{\!\{\}\!\}$ is denoted by $\varnothing$ in this paper.
\begin{itemize}\setlength{\itemsep}{0.1cm}
	\item A map $\scrR:\{\bbOmega,\bbK\}\to 2^{\mcN}$ is called a \emph{rule}. A rule $\scrR$ is said to be \emph{normal} if $M\in R({\fr t})$ whenever $M\subset N\in \scrR({\fr t})$.

	\item For each $\tau\in{\bf T}^{(0)}$ and $v\in N_\tau$, define $\mcN(v)=\{\!\{({\fr t}(e_i),{\fr e}(e_i))\}\!\}_{i=1}^n$, where $e_1,\dots,e_n\in E_\tau$ are all edges outgoing from $v$. Given a rule $\scrR$, we say that $\tau\in{\bf T}^{(0)}$ \emph{strongly conforms to $\scrR$} if $\mcN(\varrho_\tau)\in \scrR({\fr t})$ for some ${\fr t}\in\{\bbOmega,\bbK\}$ and $\mcN(v)\in \scrR({\fr t}(e))$ for each $v\in N_\tau\setminus\{\varrho_\tau\}$, where $e$ is the unique edge incoming into $v$. We write $T(\scrR)$ for the collection of all trees strongly conforming to $\scrR$, and write $T_-(\scrR)$ for the subset of $T(\scrR)$ whose elements have negative $r_{0,\infty}$-degrees.
	
	\item A normal rule $\scrR$ is said to be \emph{complete} if, in the expansion \eqref{*eq:graphicalcoprod} 
	$$
	(\tau/\sigma,{\fr t}\vert_{E_\tau\setminus E_\sigma},[{\fr n}-{\fr n}_\sigma]_\sigma,{\fr e}\vert_{E_\tau\setminus E_\sigma}+{\fr e}_{\partial\sigma})\in T(\scrR)
	$$ 
	whenever $(\tau,{\fr t},{\fr n},{\fr e})\in T(\scrR)$ and $(\sigma,{\fr t}\vert_{E_\sigma},{\fr n}_\sigma+\pi{\fr e}_{\partial\sigma},{\fr e}\vert_{E_\sigma})\in T_-(\scrR)$. (See Section 5.3 of \cite{BHZ} for a more detailed definition.)
\end{itemize}
\end{defn}

\ssk

\begin{prop}\label{*prop:ExampleDifferentiable}
Let $\scrR$ be a rule such that $\scrR(\bbOmega)=\{\varnothing\}$, $\{\!\{(\bbOmega,0)\}\!\}\in \scrR(\bbK)$, and if
$$
\{\!\{({\fr l}_1,k_1),\dots,({\fr l}_n,k_n)\}\!\}\in \scrR(\bbK),
$$
then ${\fr l}_i$ can be $\bbOmega$ at most one $i$ and then $k_i=0$.
\begin{enumerate}\setlength{\itemsep}{0.1cm}
\renewcommand{\labelenumi}{\textrm{\bf (\arabic{enumi})}}
	
	\item
Let ${\bf B}_\ocircle$ be the set of all $\tau\in T(\scrR)$ such that ${\fr n}(v)=0$ for all leaves $e=(u,v)\in{\fr t}^{-1}(\bbOmega)$ and that $1\le|{\fr t}^{-1}(\bbOmega)|<M$ for a fixed $M>1$. If the rule $\scrR$ is normal, then for some $L>0$, the set ${\bf B}=\{X^k\}_{|k|_\mfs<L}\cup{\bf B}_\ocircle$ is a differentiable sector.

	\item
Assume that ${\bf1}$ is the only element $\tau\in{\bf B}\cup\dt{\bf B}$ such that $r_{0,\infty}(\tau)=0$ or $r_{0,2}(\tau)=0$. (This becomes the case under appropriate choices of $r_0$ and $\beta_0$.)
Let ${\bf B}_-$ the set of all $\tau\in{\bf B}_\ocircle$ such that $r_{0,\infty}(\tau)<0$ and $\tau\notin\{\ocircle\}\cup\{\mcI_k^{\bbK}(\sigma)\}_{k\in\bbN^d,\,\sigma\in{\bf B}}$.
If the rule $\scrR$ is complete, for any linear map $c:W\to\bbR$ such that $c(\tau)=0$ for $\tau\in({\bf B}\cup\dt{\bf B})\setminus{\bf B}_-$, the linear map $R_c:W\to W$ defined by
$$
R_c(\tau) = \tau+(c\otimes\id)\Delta\tau\qquad(\tau\in{\bf B}\cup\dt{\bf B})
$$
is a preparation map.
\end{enumerate}
\end{prop}

\ssk

\begin{Dem}
The first assertion is obvious from the definition.
For the second assertion, the completeness ensures that $R_c(V)\subset V$ and $R_c(W)\subset W$.
The properties {\bf(\ref{Prep1})}, {\bf(\ref{PrepTriangle})}, and {\bf(\ref{EqRandI})} are immediately obtained from the definition.
The property {\bf(\ref{EqCommutationRDelta})} is obtained from the coassociativity of $\Delta$ as follows.
\begin{align*}
\big((R_c-\id)\otimes\id\big)\Delta_{0,2} &= \big( c\otimes\id\otimes P_{0,2}^+\big) (\Delta\otimes\id)\Delta   \\
&= \big( c\otimes\id\otimes P_{0,2}^+\big) (\id\otimes\Delta)\Delta = \Delta_{0,2}(R_c-\id).
\end{align*}
It remains to show the property {\bf(\ref{Prep5})}. Since the change of labels from $\bbOmega$ to $\bbH$ does not have change the expansion rule \eqref{*eq:graphicalcoprod}, the equality
$$
\Delta D = \big(D\otimes\id+\id\otimes D\big)\Delta
$$
holds. Since $c\circ D=0$, we have from the coassociativity of $\Delta$,
\begin{align*}
(R_c-\id)D&=(c\otimes\id)\Delta D = (c\otimes\id) \big(D\otimes\id+\id\otimes D\big) \Delta   \\
&=(c\otimes D)\Delta=D(R_c-\id).
\end{align*}
This concludes the proof.
\end{Dem}

%%--------------------------------------------------------------------------------------%%
\subsection{Proofs of algebraic identities$\boldmath .$ \hspace{0.03cm}}
\label{ProofLemmaComparisonFormula}
%%--------------------------------------------------------------------------------------%%

We prove two algebraic identities stated in Section \ref{SubsectionRenormalizedModels}. In the proofs we use the fact proved by Bruned in Proposition 3.15 of \cite{Bru18} that one can factorize the renormalized interpretation map by
\begin{equation} \label{EqMultiplicativeFactorization}
{\sf\Pi}_x^{R;\varepsilon,p} = \widehat{\sf\Pi}_x^{R;\varepsilon,p}R,
\end{equation}
where $\widehat{\sf\Pi}_x^{R;\varepsilon,p}$ is the linear and {\it multiplicative} map defined by
$$
\widehat{\sf\Pi}_x^{R;\varepsilon,p}(X^k) = (\cdot-x)^k,\quad
\widehat{\sf\Pi}_x^{R;\varepsilon,p}(\mcI_k\tau)=\partial^k\mcK\big(\cdot,{\sf\Pi}_x^{R;\varepsilon,p}(\tau)\big) - \hspace{-0.2cm}\underset{|l|_\mfs<r_{\varepsilon,p}(\mcI_k\tau)}{\sum_{l\in\bbN^d}}\frac{(\cdot-x)^l}{l!} \, \partial^{k+l}\mcK\big(x,{\sf\Pi}_x^{R;\varepsilon,p}(\tau)\big).
$$

\ssk

\noindent \textbf{{Proof of Lemma {{\ref{Lem:dPi=PiD}}} --}}
In addition to \eqref{Eq:dPi=PiD}, we also prove the similar identity
\begin{equation}\label{Eq:dPi=PiD'}
d_\xi\big(\widehat{\sf\Pi}_x^{\xi,R;\varepsilon,\infty}(\tau)\big)(h) = \widehat{\sf\Pi}_x^{\xi,h,R;\varepsilon,\infty}(D\tau)
\end{equation}
simultaneously. The proof is an induction on the preorder $\preceq$ defined by \eqref{EqDefnOrder}.
Both \eqref{Eq:dPi=PiD} and \eqref{Eq:dPi=PiD'} are obvious for the initial cases $\tau\in\{\ocircle,X^k\}$. 
Next we prove that \eqref{Eq:dPi=PiD'} holds for any $\tau$ by assuming that both \eqref{Eq:dPi=PiD} and \eqref{Eq:dPi=PiD'} hold for any trees with less edges than $\tau$. First let $\tau$ be a planted tree of the form $\mcI_k(\sigma)$ with $\sigma\in{\bf B}$. If $\sigma$ satisfies \eqref{Eq:dPi=PiD}, then by the definition of $\widehat{\sf\Pi}_x^{\xi,R;\varepsilon,\infty}$ operator,
\begin{align*}
d_\xi\big(\widehat{\sf\Pi}_x^{\xi,R;\varepsilon,\infty}(\tau)\big) (h)
&=d_\xi\bigg(\partial^k\mcK\big(\cdot,{\sf\Pi}_x^{\xi,R;\varepsilon,\infty}(\sigma)\big) - \sum_{|l|_\mfs<r_{\varepsilon,\infty}(\mcI_k\sigma)}\frac{(\cdot-x)^l}{l!} \, \partial^{k+l}\mcK\big(x,{\sf\Pi}_x^{\xi,R;\varepsilon,\infty}(\sigma)\big) \bigg)(h)   \\
&=\partial^k\mcK\big(\cdot,{\sf\Pi}_x^{\xi,h,R;\varepsilon,\infty}(D\sigma)\big)-\sum_{|l|_\mfs<r_{\varepsilon,\infty}(\mcI_k\sigma)}\frac{(\cdot-x)^l}{l!} \, \partial^{k+l}\mcK\big(x,{\sf\Pi}_x^{\xi,h,R;\varepsilon,\infty}(D\sigma)\big)   \\
&=\widehat{\sf\Pi}_x^{\xi,h,R;\varepsilon,\infty}(\mcI_kD\sigma) = \widehat{\sf\Pi}_x^{\xi,h,R;\varepsilon,\infty}(D\tau).
\end{align*}
In the third equality, we use the fact that $D$ preserves the $r_{\varepsilon,\infty}$-degree: $r_{\varepsilon,\infty}(\mcI_k(D\sigma))=r_{\varepsilon,\infty}(\mcI_k\sigma)$. Thus $\tau$ satisfies \eqref{Eq:dPi=PiD'}.
Next we consider a non-planted tree $\tau$ factorized by $\tau=\prod_{i=0}^N\eta_i$ with $\eta_0=X^k$ and planted trees $\eta_i$ ($1\le i\le N$). 
If each $\eta_i$ satisfies \eqref{Eq:dPi=PiD'}, then we use the multiplicativity of $\widehat{\sf\Pi}_x^{\xi,R;\varepsilon,\infty}$ and Leibniz rules for $d_\xi$ and $D$ to derive
\begin{align*}
d_\xi\big(\widehat{\sf\Pi}_x^{\xi,R;\varepsilon,\infty}(\tau)\big)&(h)
= d_\xi\bigg(\prod_i\big(\widehat{\sf\Pi}_x^{\xi,R;\varepsilon,\infty}(\eta_i)\big) \bigg) (h)
= \sum_i d_\xi\big(\widehat{\sf\Pi}_x^{\xi,R;\varepsilon,\infty}(\eta_i)\big) (h) \prod_{j\neq i}\big(\widehat{\sf\Pi}_x^{\xi,R;\varepsilon,\infty}(\eta_j)\big)   \\
&= \sum_i \big(\widehat{\sf\Pi}_x^{\xi,h,R;\varepsilon,\infty}(D\eta_i)\big) (h) \prod_{j\neq i}\big(\widehat{\sf\Pi}_x^{\xi,R;\varepsilon,\infty}(\eta_j)\big)   \\
&= \widehat{\sf\Pi}_x^{\xi,h,R;\varepsilon,\infty}\bigg(\sum_i(D\eta_i)\prod_{j\neq i}\eta_j\bigg) = \widehat{\sf\Pi}_x^{\xi,h,R;\varepsilon,\infty}(D\tau).
\end{align*}
Thus $\tau$ satisfies \eqref{Eq:dPi=PiD'}.
Finally, we prove that $\tau$ also satisfies \eqref{Eq:dPi=PiD} by using the commutation \eqref{Prep5} between $R$ and $D$. Indeed, we have
\begin{align*}
d_\xi\big({\sf\Pi}_x^{\xi,R;\varepsilon,\infty}(\tau)\big)(h)
&=d_\xi\big(\widehat{\sf\Pi}_x^{\xi,R;\varepsilon,\infty}(R\tau)\big)(h)
=\widehat{\sf\Pi}_x^{\xi,h,R;\varepsilon,\infty}\big(D(R\tau)\big)   \\
&=\widehat{\sf\Pi}_x^{\xi,h,R;\varepsilon,\infty}\big(R(D\tau)\big)
={\sf\Pi}_x^{\xi,h,R;\varepsilon,\infty}(D\tau).
\end{align*}
In the third equality, we use the triangular property \eqref{PrepTriangle} of $R$ and the inductive assumption that all trees $\sigma$ with $|\sigma|_\ocircle<|\tau|_\ocircle$ satisfy \eqref{Eq:dPi=PiD'}. Thus $\tau$ satisfies \eqref{Eq:dPi=PiD}.

\medskip

\noindent \textbf{{Proof of Lemma {{\ref{LemComparisonFormula}}} --}} To simplify the notations we suppress the $R,\varepsilon$ dependence in the notation and we write ${\sf\Pi}_x^p,\lambda_x^p,\Delta_p,r_p(\tau),p(\mu)$ for ${\sf\Pi}_x^{R;\varepsilon,p},\lambda_x^{R;\varepsilon,p},\Delta_{\varepsilon,p},r_{\varepsilon,p}(\tau),p_\varepsilon(\mu)$ respectively.
Since $\lambda_x^p$ vanishes on non-planted trees, we can write \eqref{EqSyntheticDecompositionFormulaDerivativePi} in the explicit form
\begin{equation}\label{*EqSyntheticDecompositionFormulaDerivativePi}
{\sf\Pi}_x^p(\tau)={\sf\Pi}_x^2(\tau)
+\sum_{e\in E_\tau}\sum_{{k}\in\bbN^d} \frac1{{k}!} \, \lambda_x^{p}\big(\mcI_{{\fr e}(e)+k}^{{\fr t}(e)}\tau_+^e\big) \, {\sf\Pi}_x^{p}\big(\uparrow_{e_-}^k\tau_-^e\big),
\end{equation}
by \eqref{*eq:graphicalcoprod}.
Here, for any arbitrary tree $\tau\in{\bf T}$ and $e=(e_-,e_+)\in E_\tau$ denote by $\{\tau_+^e,\tau_-^e\}$ the connected components of $\tau\setminus\{e\}$ such that $\tau_-^e$ contains the root $\varrho_\tau$ and the node $e_-$. Moreover, for any $\sigma\in{\bf T}$, ${k}\in\bbN^d$, and $v\in N_\sigma$ denote by $\uparrow^{k}_v\sigma$ the decorated tree $\sigma$ with the same decorations as $\sigma$ except that the node $v$ has now decoration ${\frak n}(v)+{k}$.
Note that $\lambda_x^p$ does not vanish only at trees with an $\bbH$-labeled edge. If $\mcI_{{\fr e}(e)+k}^{{\fr t}(e)}\tau_+^e$ contains such an edge, then $\uparrow_{e_-}^k\tau_-^e$ does not. Therefore, ${\sf\Pi}_x^{p}(\uparrow_{e_-}^k\tau_-^e)$ does not depend on $p$ and we can write ${\sf\Pi}_x^{p}(\uparrow_{e_-}^k\tau_-^e)={\sf\Pi}_x(\uparrow_{e_-}^k\tau_-^e)$.
We prove \eqref{*EqSyntheticDecompositionFormulaDerivativePi} and the auxiliary formula
\begin{equation}\label{EqSyntheticDecompositionFormulaDerivativePi'}
\widehat{\sf\Pi}_x^p(\tau) = \widehat{\sf\Pi}_x^2(\tau) + \sum_{e\in E_\tau}\sum_{{k}\in\bbN^d} \frac1{{k}!} \, \lambda_x^{p}\big(\mcI_{{\fr e}(e)+k}^{{\fr t}(e)}\tau_+^e\big) \, \widehat{\sf\Pi}_x(\uparrow_{e_-}^k\tau_-^e)
\end{equation}
for any $\tau\in\dt{\bf B}$ by a similar induction to the proof of Lemma \ref{Lem:dPi=PiD}.

\ssk

For the initial case $\tau=\mcI_0^{\bbH}(\bullet)$, by the second formula of \eqref{*Eq:ModelXH},
\begin{align*}
{\sf\Pi}_x^p(\tau) &= {\sf\Pi}_x^2(\tau) + \sum_{r_{p}(\bbH)\le |k|_\mfs<r_{2}(\bbH)}
\frac{(y-x)^k}{k!} \, \partial^kh_n(x)\\
&={\sf\Pi}_x^2(\tau) + \sum_{k\in\bbN^d}
\frac1{k!} \, \lambda_x^p(\mcI_k^{\bbH}(\bullet)) \, {\sf\Pi}_x(X^k).
\end{align*}
Since $\widehat{\sf\Pi}_x^p(\tau)={\sf\Pi}_x^p(\tau)$ and $\widehat{\sf\Pi}_x^p(X^k)={\sf\Pi}_x^p(X^k)$, both \eqref{*EqSyntheticDecompositionFormulaDerivativePi} and \eqref{EqSyntheticDecompositionFormulaDerivativePi'} hold for $\tau=\mcI_0^{\bbH}(\bullet)$. Let $\tau$ be a planted tree $\mcI_k^{\bbK}(\sigma)$ with $\sigma\in\dt{\bf B}$. By the definition of $\widehat{\sf\Pi}_x^p$ and the inductive assumption that \eqref{*EqSyntheticDecompositionFormulaDerivativePi} holds for $\sigma\in\dt{\bf B}$, we have
\begin{align*}
\widehat{\sf\Pi}_x^p(\tau) &= \partial^{k}\mcK\big(\cdot,{\sf\Pi}_x^p(\sigma)\big) - \sum_{|{l}|_{\fr s} < r_{p}(\tau)}\frac{(\cdot-x)^l}{{l}!} \, \partial^{{k} + {l} }\mcK\big(x,{\sf\Pi}_x^{p}(\sigma)\big)   \\
&= \bigg(\partial^{k} \mcK\big(\cdot,{\sf\Pi}_x^{2}(\sigma)\big) - \sum_{|{l}|_{\fr s}<r_{p}(\tau)} \frac{(\cdot-x)^l}{{l}!} \, \partial^{{k} + {l} }\mcK\big(x,{\sf\Pi}_x^{2}(\sigma)\big)\bigg)   \\
&\quad+ \sum_{e,{m}}\frac1{{m}!} \, \lambda_x^{p}(\sigma_+^{e,m})
\bigg(\partial^{k} \mcK\big(\cdot,{\sf\Pi}_x(\sigma_-^{e,m})\big) -  \sum_{|{l}|_{\fr s}<r_{p}(\tau)} \frac{(\cdot-x)^l}{{l}!} \,\partial^{{k} + {l} }\mcK\big(x,{\sf\Pi}_x(\sigma_-^{e,m})\big) \bigg),
\end{align*}
where we write $\sigma_+^{e,m}=\mcI_{m+{\fr e}(e)}^{{\fr t}(e)}\sigma_+^e$ and $\sigma_-^{e,m}=\,\uparrow_{e_-}^{m} \sigma_-^e$ for simplicity.
Here we can compare the two quantities in the big parentheses with $\widehat{\sf\Pi}_x^{2}(\tau)$ and $\widehat{\sf\Pi}^{p}_x(\mcI_k^{\bbK}(\sigma_-^{e,m}))$ except different domains for ${l}$ in both cases.
Recall that $r_{p}(\tau)\le r_{2}(\tau)$. Also, since $r_p(\sigma_+^{e,m})\le0$ for $m$ such that $\lambda_x^{p}(\sigma_+^{e,m})\neq0$, we have 
$$
r_{p}(\tau) = r_p(\mcI_{k}^{\bbK}\sigma) = r_p(\sigma_+^{e,m})+r_p(\mcI_{k}^{\bbK}(\sigma_-^{e,m}))
\le r_p(\mcI_{k}^{\bbK}(\sigma_-^{e,m})).
$$ 
One therefore has
\begin{align}\label{Eq:ProofComparison}
\begin{aligned}
\widehat{\sf\Pi}_x^{p}(\tau) &= \widehat{\sf\Pi}_x^2(\tau) + \sum_{r_{p}(\tau)\le|{l}|_{\fr s}<r_{2}(\tau)}\frac{(\cdot-x)^l}{{l}!}\,\partial^{{k} + {l}}\, \mcK\big(x,{\sf\Pi}_x^{2}(\sigma)\big)   \\
&\quad+\sum_{e,m}\frac1{m!}\,\lambda_x^{p}(\sigma_+^{e,m}) \bigg( \widehat{\sf\Pi}_x\big(\mcI_{k}^{\bbK}(\sigma_-^{e,m})\big) + \sum_{r_{p}(\tau)\le |{l}|_{\fr s}< r_p(\mcI_{k}^{\bbK}(\sigma_-^{e,m}))} \frac{(\cdot-x)^l}{{l}!} \, \partial^{{k} + {l} }\mcK\big(x,{\sf\Pi}_x(\sigma_-^{e,m})\big)\bigg).
\end{aligned}
\end{align}
For the second term in the right hand side, $r_{p}(\tau)\le|{l}|_{\fr s}<r_{2}(\tau)$ is equivalent to $r_p(\mcI_{k+l}^{\bbK}\sigma)\le0<r_2(\mcI_{k+l}^{\bbK}\sigma)$. By using the assumption \eqref{*EqSyntheticDecompositionFormulaDerivativePi} on $\sigma$ with the integrability exponent $p(\mcI_{{k} + {l} }^{\bbK}\sigma)-\delta$ for sufficiently small $\delta>0$ (fixed later), we have that this term is equal to
%\vfill \pagebreak
\begin{equation*} \begin{split}
&\sum_{r_p(\mcI_{k+l}^{\bbK}\sigma)\le0<r_2(\mcI_{k+l}^{\bbK}\sigma)} \frac{(\cdot-x)^l}{{l}!}\, \partial_x^{{k} + {l} }\mcK\big(x,{\sf\Pi}_x^{p(\mcI_{k+l}^{\bbK}\sigma)-\delta}(\sigma)\big)  \\
&\hspace{1.4cm}- \sum_{r_p(\mcI_{k+l}^{\bbK}\sigma)\le0<r_2(\mcI_{k+l}^{\bbK}\sigma)}\frac{(\cdot-x)^l}{{l}!}\sum_{e,{m}}\frac1{{m}!} \, \lambda_x^{p(\mcI_{{k} + {l} }^{\bbK}\sigma)-\delta}(\sigma_+^{e,m})  \, \partial^{{k} + {l} }\mcK\big(x,{\sf\Pi}_x(\sigma_-^{e,m})\big).
\end{split} \end{equation*}
The first term is of the form $\sum_l\frac1{l!} \, \lambda_x^p(\mcI_{k+l}^{\bbK}\sigma)\,\widehat{\sf\Pi}_x^p(X^l)$, which appears in \eqref{EqSyntheticDecompositionFormulaDerivativePi'} in the case that $e$ is the unique edge connected to the root of $\tau=\mcI_k^{\bbK}\sigma$. The other terms in \eqref{EqSyntheticDecompositionFormulaDerivativePi'} appear in the third term of the right hand side of \eqref{Eq:ProofComparison} as the form $\sum_{e,m}\frac1{m!}\,\lambda_x^{p}(\sigma_+^{e,m}) \, \widehat{\sf\Pi}^{p}_x(\mcI_{k}^{\bbK}(\sigma_-^{e,m}))$.
Therefore to show that \eqref{EqSyntheticDecompositionFormulaDerivativePi'} holds for $\tau$ it is sufficient to show the identity
\begin{equation*} \begin{split}
&\sum_{e,{m}} \frac1{{m}!} \, \lambda_x^{p}(\sigma_+^{e,m}) \sum_{r_p(\tau)\le |{l}|_{\fr s}<r_p(\mcI_{k}^{\bbK}(\sigma_-^{e,m}))} \frac{(\cdot-x)^l}{{l}!} \, \partial^{{k} + {l} }\mcK\big(x,{\sf\Pi}_x(\sigma_-^{e,m})\big)   \\
&=
\sum_{r_p(\mcI_{k+l}^{\bbK}\sigma)\le0<r_2(\mcI_{k+l}^{\bbK}\sigma)}\frac{(\cdot-x)^l}{{l}!}\sum_{e,{m}}\frac1{{m}!} \, \lambda_x^{p(\mcI_{{k} + {l} }^{\bbK}\sigma)-\delta}(\sigma_+^{e,m})  \, \partial^{{k} + {l} }\mcK\big(x,{\sf\Pi}_x(\sigma_-^{e,m})\big),
\end{split} \end{equation*} 
that is, to show that the following two conditions are equivalent when we choose any small $\delta>0$.   \vspace{0.1cm}
\begin{itemize}  \setlength{\itemsep}{0.15cm}
	\item[(a)] $r_p(\sigma_+^{e,m})\le 0 < r_2(\sigma_+^{e,m})$ and $r_p(\tau)\le|{l}|_{\fr s}<r_p(\mcI_{k}^{\bbK}(\sigma_-^{e,m}))$,
	\item[(b)] $r_p(\mcI_{k+l}^{\bbK}\sigma)\le0<r_2(\mcI_{k+l}^{\bbK}\sigma)$ and $r_{p(\mcI_{k+l}^{\bbK}\sigma)-\delta}(\sigma_+^{e,m}) \leq 0 < r_{2}(\sigma_+^{e,m})$.   \vspace{0.1cm}
\end{itemize}
First we assume (a). Since $r_2(\mcI_{k}(\sigma_-^{e,m}))=r_p(\mcI_{k}(\sigma_-^{e,m}))>|l|_{\fr s}$,
\begin{align*}
r_2(\mcI_{k+l}^{\bbK}\sigma)=
r_2(\sigma_+^{e,m}) + r_2(\mcI_{k+l}^{\bbK}(\sigma_-^{e,m})) > 0.
\end{align*}
Also, since $r_{p(\mcI_{k+l}^{\bbK}\sigma)-\delta}(\mcI_{k+l}^{\bbK}\sigma) \to0$ as $\delta\to0$ and $r_p(\mcI_{k+l}^{\bbK}(\sigma_-^{e,m}))$ is a $p$-independent positive number, if $\delta>0$ is small enough then
\begin{align*}
r_{p(\mcI_{k+l}^{\bbK}\sigma)-\delta}(\sigma_+^{e,m}) = r_{p(\mcI_{k+l}^{\bbK}\sigma)-\delta}(\mcI_{k+l}^{\bbK}\sigma) - r_{p(\mcI_{k+l}^{\bbK}\sigma)-\delta}(\mcI_{k+l}^{\bbK}(\sigma_-^{e,m}))\leq 0.
\end{align*}
Thus (b) holds. Next we assume (b). Note that $r_p$-degree is monotonically decreasing in $p$.
Since $r_p(\mcI_{k+l}^{\bbK}\sigma)\le0$ implies $p\ge p(\mcI_{k+l}^{\bbK}\sigma)$, we have
$$
r_p(\sigma_+^{e,m})\le r_{p(\mcI_{k+l}^{\bbK}\sigma)-\delta}(\sigma_+^{e,m}) \le0.
$$
Moreover, since $r_{p(\mcI_{k+l}^{\bbK}\sigma)-\delta}(\mcI_{k+l}^{\bbK} \sigma)>0$ for any $\delta>0$, one also has
\begin{align*}
r_p(\mcI_{k+l}^{\bbK}(\sigma_-^{e,m})) = 
r_{p(\mcI_{k+l}^{\bbK}\sigma)-\delta}(\mcI_{k+l}^{\bbK}(\sigma_-^{e,m})) =
r_{p(\mcI_{k+l}^{\bbK}\sigma)-\delta}(\mcI_{k+l}^{\bbK} \sigma) - r_{p(\mcI_{k+l}^{\bbK}\sigma)-\delta}(\sigma_+^{e,m}) > 0
\end{align*}
Thus (a) holds. Since (a) and (b) are equivalent, \eqref{EqSyntheticDecompositionFormulaDerivativePi'} holds for $\tau$.

\ssk

We now consider non-planted trees of the form $\tau=\sigma\mu$ where $\sigma\in{\bf B}$ and $\mu$ is a planted tree with an $\bbH$-labeled edge. Note that the cut at $e$ in \eqref{EqSyntheticDecompositionFormulaDerivativePi'} makes sense only if $e\in E_\mu$ since $\lambda_x^p$ vanishes on trees without $\bbH$ label. Assuming that \eqref{EqSyntheticDecompositionFormulaDerivativePi'} holds for $\mu$, and recalling the multiplicativity of $\widehat{\sf\Pi}_x^p$, we have
\begin{align*}
\widehat{\sf\Pi}_x^{p}(\tau)
&= \widehat{\sf\Pi}_x(\sigma) \,\widehat{\sf\Pi}_x^{p}(\mu)   \\
&= \widehat{\sf\Pi}_x(\sigma) \bigg(\widehat{\sf\Pi}_x^2(\mu) + \sum_{e\in E_\mu}\sum_{{k}\in\bbN^d} \frac1{{k}!} \, \lambda_x^{p}(\mcI_{{\fr e}(e)+k}^{{\fr t}(e)}\mu_+^e) \, \widehat{\sf\Pi}_x(\uparrow_{e_-}^k\mu_-^e)\bigg)
\\
&= \widehat{\sf\Pi}_x^2(\sigma\mu) + \sum_{e\in E_\mu}\sum_{{k}\in\bbN^d} \frac1{{k}!} \, \lambda_x^{p}\big(\mcI_{{\fr e}(e)+k}^{{\fr t}(e)}\mu_+^e\big) \, \widehat{\sf\Pi}_x\big(\sigma(\uparrow_{e_-}^k\tau_-^e)\big)  \\
&= \widehat{\sf\Pi}_x^2(\tau) + \sum_{e\in E_\tau}\sum_{{k}\in\bbN^d} \frac1{{k}!} \, \lambda_x^p\big(\mcI_{{\fr e}(e)+k}^{{\fr t}(e)}\tau_+^e\big) \, \widehat{\sf\Pi}_x(\uparrow_{e_-}^k\tau_-^e).
\end{align*}
Thus \eqref{EqSyntheticDecompositionFormulaDerivativePi'} holds for $\tau$.

\ssk

Finally, by using the commutation {\bf(\ref{EqCommutationRDelta})} between $R$ and $\Delta_2$ and the triangular property {\bf(\ref{PrepTriangle})} of $R$ with respect to $|\cdot|_\ocircle$, we have
\begin{align*}
({\sf\Pi}_x^p - {\sf\Pi}_x^2)(\tau)
&= (\widehat{\sf\Pi}_x^{p} - \widehat{\sf\Pi}_x^{2})(R\tau )
= (\widehat{\sf\Pi}_x^{p} \otimes \lambda_x^{p})\Delta_{2}(R\tau)    \\
&= (\widehat{\sf\Pi}_x^{p}R\otimes\lambda_x^{p})\Delta_{2}(\tau) = ({\sf\Pi}_x^p\otimes\lambda_x^{p})\Delta_{2}(\tau).
\end{align*}
Thus \eqref{*EqSyntheticDecompositionFormulaDerivativePi} holds for $\tau$.

\medskip

%%-----------------------------------------------------------------%%
\subsection{Proof of the Lemma \ref{*lem:BesovTaylor}$\boldmath .$ \hspace{0.03cm}}
\label{SectionBesovTaylor}
%%-----------------------------------------------------------------%%

Fix $p\in[1,\infty]$, $r>0$ and $h\in B_{p,\infty}^r=B_{p,\infty}^r(w_0)$. Since $\mcQ_t=e^{t\mcL}$ we can write
$$
h = \mcQ_1(h) - \int_0^1\mcL\mcQ_t(h) dt.
$$
We only focus on the estimate of
$$
\int_0^1\bigg\{\mcL\mcQ_t(h)(x+y)-\sum_{|k|_\mfs<r}\frac{y^k}{k!}\partial^k\mcL\mcQ_t(h)(x)\bigg\}dt\eqdef\int_0^1H_t(x,y)dt,
$$
since the estimate of $\mcQ_1h$ is easier. By Proposition A.1 of \cite{Hai14}, there exists a family of Borel probability measures $\{M^k\}_{k\in\bbN^d}$ on $[0,1]^d$ and
$$
\mcL\mcQ_t(h)(x+y)-\sum_{k\in A}\frac{y^k}{k!}\partial^k\mcL\mcQ_t(h)(x)
=\sum_{l\in\partial A}\frac{x^l}{l!}\int_{[0,1]^d}\partial^l\mcL\mcQ_t(h)(x+\theta y) \, M^k(d\theta),
$$
where $\theta y=(\theta_jy_j)_{j=1}^d$, $A$ is any finite subset of $\bbN^d$ such that $l\in A$ whenever $k\in A$ and $l\le k$, and $\partial A=\{k\in\bbN^d\,;\,k\notin A,\ k-e_{m(k)}\in A\}$ with $m(k)=\min\{j\,;\,k_j\neq0\}$.
In particular, if $A=\{k\,;\,|k|_\mfs<r\}$, then $\partial A\subset\{k\,;\,r<|k|_\mfs<r+s\}$ where $s=\max_j\mfs_j$.
Moreover, since $\mcL$ is a differential operator of order $\ell$, we have that for any $|k|_\mfs>r-\ell$,
$$
\|\partial^k\mcL\mcQ_t(h)\|_{L^p(\bbR^d)}\lesssim t^{\frac{r-\ell-|k|_\mfs}\ell}\|h\|_{B_{p,\infty}^r}.
$$
First let $\|y\|_\mfs\le1$ and consider the integral over $t\in(\|y\|_\mfs^\ell,1]$.
\begin{align*}
\Bigg\|\int_{\|y\|_\mfs^\ell}^1
H_t(x,y)dt
\Bigg\|_{L_x^p}
&\le
\int_{\|y\|_\mfs^\ell}^1dt\sum_{r<|l|_\mfs<r+s}\|y\|_\mfs^{|l|}\int_{[0,1]^d}\|\partial^l\mcL\mcQ_t(h)(x+\theta y)\|_{L_x^p} \, M^k(d\theta)   \\
&\lesssim\|h\|_{B_{p,\infty}^r}\int_{\|y\|_\mfs^\ell}^1\sum_{r<|l|_\mfs<r+s}\|y\|_\mfs^{|l|}t^{\frac{r-\ell-|l|_\mfs}\ell}dt
\lesssim\|h\|_{B_{p,\infty}^r}\|y\|_\mfs^r.
\end{align*}
For the integral over $t\in\big(0,1\wedge\|y\|_\mfs^\ell\big]$, we decompose
\begin{align*}
H_t(x,y)&=\bigg\{\mcL\mcQ_t(h)(x+y)-\sum_{|k|_\mfs<r-s}\frac{y^k}{k!}\partial^k\mcL\mcQ_t(h)(x)\bigg\}
-\sum_{r-s<|k|_\mfs<r}\frac{y^k}{k!}\partial^k\mcL\mcQ_t(h)(x)\\
&\eqdef H_t^1(x,y)+H_t^2(x,y).
\end{align*}
For the integral of $H_t^1$, by noting that $r-s>r-\ell$, we have
\begin{align*}
\Bigg\|\int_0^{1\wedge\|y\|_\mfs^\ell}
H_t^1(x,y)dt
\Bigg\|_{L_x^p}
&\le\int_0^{1\wedge\|y\|_\mfs^\ell}dt
\sum_{r-s<|l|_\mfs<r}\|y\|_\mfs^{|l|}\int_{[0,1]^d}\|\partial^l\mcL\mcQ_t(h)(x+\theta y)\|_{L_x^p}M^k(d\theta)\\
&\lesssim\|h\|_{B_{p,\infty}^r}\int_0^{1\wedge\|y\|_\mfs^\ell}
\sum_{r-s<|l|_\mfs<r}\|y\|_\mfs^{|l|}t^{\frac{r-\ell-|l|_\mfs}\ell}dt
\lesssim\|h\|_{B_{p,\infty}^r}\|y\|_\mfs^r.
\end{align*}
We have a similar bound for the integral of $H_t^2$.
\hfill$\rhd$

\medskip

%%-----------------------------------------------------------------%%
\subsection{Proof of the lemmas of Section {{\ref{sec:ind1}}}$\boldmath .$ \hspace{0.03cm}}
\label{SectionLastProofsLemmas}
%%-----------------------------------------------------------------%%

We have four statements to prove.

\ssk

\noindent \textbf{{Proof of Lemma {{\ref{*argumentforreconst}}} --}}
Recall from Definition \ref{DefnMD} in Appendix \ref{SectionAppendix} the definition of the quantities $\llparenthesis f\rrparenthesis_{{\bf c};w_b}$ and $\Vert f\Vert^{\sf\Gamma}_{{\bf c};w_b}$ involved in the definition of modelled distribution. We can write $f_{\tau}^{n;\varepsilon,p}(x) = {\sf\Gamma}_x^{n;\varepsilon,p}(\tau)-\tau$, with ${\sf\Gamma}_x^{n;\varepsilon,p}=(\id\otimes\mathop{{\sf g}^{n;\varepsilon,p}_x}) \Delta_{\varepsilon,p}$.
By Lemma \ref{leminduction*}, $f_\tau^{n;\varepsilon,p}$ takes values in $C_{\prec\tau}$ and 
$$
\lp f_{\tau}^{n;\varepsilon,p} \rp_{(r_{\varepsilon,p}(\tau),i_p(\tau)) ; w_{cm_{\bf B}}} 
\lesssim \big( 1+\|{\sf g}^{n;\varepsilon,p} \res {\bf C}_{\prec\tau,\varepsilon,p}\|_{w_c}\big)^{m_{\bf B}} < \infty.
$$
Since ${\sf\Gamma}_{yx}^{n;\varepsilon,p}{\sf\Gamma}_{x}^{n;\varepsilon,p}={\sf\Gamma}_{y}^{n;\varepsilon,p}$ the relation
\begin{align*}
{\sf\Gamma}_{yx}^{n;\varepsilon,p} \big(f_{\tau}^{n;\varepsilon,p}(x)\big) - f_{\tau}^{n;\varepsilon,p}(y) 
={\sf\Gamma}_{yx}^{n;\varepsilon,p}({\sf\Gamma}_{x}^{n;\varepsilon,p}-\id)\tau-({\sf\Gamma}_{y}^{n;\varepsilon,p}-\id)\tau= (\id-\mathop{{\sf\Gamma}_{yx}^{n;\varepsilon,p}})\tau
\end{align*}
holds. Hence $\| f_{\tau}^{n;\varepsilon,p} \|_{(r_{\varepsilon,p}(\tau),p);w_{cm_{\bf B}}}^{{\sf \Gamma}^{n;\varepsilon,p}} \le \|{\sf\Gamma}^{n;\varepsilon,p}\res\tau\|_{w_{cm_{\bf B}}}<\infty$. Moreover the relation ${\sf\Pi}_x^{n;\varepsilon,p}{\sf\Gamma}_x^{n;\varepsilon,p}={\sf\Pi}^n$ yields
$$
{\sf\Pi}^{n;\varepsilon,p}_x(\tau) = {\sf\Pi}^n(\tau) - {\sf\Pi}^{n;\varepsilon,p}_x \big(f_{\tau}^{n;\varepsilon,p}(x)\big).
$$
Since the left hand side has a (possibly $(n,\omega)$-dependent) bound 
$$
\big\| \mcQ_t\big(x,{\sf\Pi}^{n;\varepsilon,p}_x(\tau)\big)\big\|_{L_x^{i_p(\tau)}(w_c)}\lesssim t^{\frac{r_{\varepsilon,p}(\tau)}\ell},
$$ 
this bound means that ${\sf\Pi}^n\tau$ is a reconstruction of $f_{\tau}^{n;\varepsilon,p}$.
\hfill$\rhd$

\ssk

\noindent \textbf{{Proof of Lemma {{\ref{lem:Biposi}}} --}}
We apply Lemma \ref{*argumentforreconst} to $\tau=\tau_i$. Then $C_{\prec\tau}$ can be replaced by $V_{i-1}$. Since $r_{\varepsilon,\infty}(\tau_i)$, the reconstruction of $f_{\tau_i}^{n;\varepsilon,\infty}\in D^{(r_{\varepsilon,\infty}(\tau_i),\infty)}(V_{i-1} ; {\sf \Gamma}^{n;\varepsilon,\infty})_{w_{cm_{\bf B}}}$ is unique and coincides with ${\sf \Pi}^n\tau_i$.
From the inequalities obtained in the proof of Lemma \ref{*argumentforreconst}, we have
\begin{align*}
\|{\sf\Pi}^{n;\varepsilon,\infty} \hspace{-0.03cm}:\hspace{-0.03cm} \tau_i\|_{w_{c(m_{\bf B}+1)}} 
&=\lb{\mcR}^{{\sf M}^{n;\varepsilon,\infty}} f_{\tau_i}^{n;\varepsilon,\infty}\rb_{(r_{\varepsilon,\infty}(\tau_i),\infty);w_{c(m_{\bf B}+1)}}\\
&\le C \|{\sf\Pi}^{n;\varepsilon,\infty} \hspace{-0.03cm}:\hspace{-0.03cm} {\bf B}_{i-1}\|_{w_c} \, 
\| f_{\tau}^{n;\varepsilon,\infty} \|_{(r_{\varepsilon,\infty}(\tau_i),\infty);w_{cm_{\bf B}}}^{{\sf \Gamma}^{n;\varepsilon,p\infty}}\\
&\le C \|{\sf\Pi}^{n;\varepsilon,\infty} \hspace{-0.03cm}:\hspace{-0.03cm} {\bf B}_{i-1}\|_{w_c} \, \|{\sf\Gamma}^{n;\varepsilon,\infty} \hspace{-0.03cm}:\hspace{-0.03cm} {\bf B}_i\|_{w_{cm_{\bf B}}}.
\end{align*}
The constant is the one appearing in the reconstruction theorem (Theorem \ref{thm:besovreconstruction}), so it is universal and does not depend on $n$ nor $\omega$.
\hfill$\rhd$

\ssk

We prove some useful technical results before going to the proof of Lemma \ref{lem:stationary}.

\ssk

\begin{lem}\label{lem:wyoung}
There exists a finite positive constant $C_c$ depending only on $(Q_t)_{t>0}$ and $c>0$ such that, for any $q\ge p\in[1,\infty]$, $t\in(0,1]$, and $f\in L^p(w_c)$, one has
\begin{align*}
\|G_t*f\|_{L^q(w_c)}\le C_c\,t^{\frac{|\mfs|}\ell(\frac1q-\frac1p)}\|f\|_{L^p(w_c)}.
\end{align*}
\end{lem}

\ssk

\begin{Dem}
Since $|(G_t*f)w_c|\le(G_tw_{c}^{-1})*(|f|w_c)$, the result follows from usual Young inequality. The proportional constant is $\|G_tw_{c}^{-1}\|_{L^r(\bbR^d)}$, where $\frac1r=1+\frac1q-\frac1p$. By \eqref{ineq:Gw} we have the bound
$$
\|G_tw_{c}^{-1}\|_{L^r(\bbR^d)}\lesssim t^{-\frac{|\mfs|}\ell(1-\frac1r)}=t^{\frac{|\mfs|}\ell(\frac1q-\frac1p)}
$$
which concludes the proof.
\end{Dem}

\ssk

\begin{lem}\label{LemRemainder}
For any $\varepsilon\ge0$, $p\in[2,\infty]$, and $c,c'>0$, there exists a positive constant $C$, which is independent of $n$, $\omega$, and $h\in H$ with $\Vert h\Vert_H\leq 1$, one has
$$
\bigg\Vert\int_{\bbR^d}G_t(x-y)\big\vert \mcQ_t\big(y,{\sf\Pi}_x^{n;\varepsilon,p}\tau-{\sf\Pi}_y^{n;\varepsilon,p}\tau\big)\big\vert dy\bigg\Vert_{L_x^{i_p(\tau)}(w_{c+c'})}
\le Ct^{\frac{r_{\varepsilon,p}(\tau)}\ell}
\|{\sf\Gamma}^{n;\varepsilon,p}\res\tau\|_{w_{c'}}  \|{\sf\Pi}^{n;\varepsilon,p} \res  {\bf C}_{\prec\tau}\|_{w_c}
$$
for any $\tau\in{\bf B}\cup\dt{\bf B}$.
\end{lem}

\ssk

\begin{Dem}
To simplify the notations we omit here the symbols $n,\varepsilon,p$ from the model. By the change of variable $y\to x-z$ and the expansion of ${\sf\Gamma}_{yx}$, we can write the quantity inside the $L_x^{i_p(\tau)}$ norm as
\begin{align} \label{EqIntermediate}
\int_{\bbR^d}G_t(z)\big\vert \mcQ_t\big(x-z,{\sf\Pi}_{x-z}({\sf\Gamma}_{(x-z)x}-\id)\tau\big)\big\vert dz
=\sum_{\sigma\prec\tau}
\int_{\bbR^d}G_t(z)\vert P_\sigma{\sf\Gamma}_{(x-z)x}\tau\vert \big\vert \mcQ_t(x-z,{\sf\Pi}_{x-z}(\sigma))\big\vert dz.
\end{align}
By H\"older inequality and the elementary inequality $\|f(x-z)\|_{L_x^p(\bbR^d)}\le w_{c}(z)^{-1}\|f\|_{L_x^p(\bbR^d)}$, we can bound above the $L_x^{i_p(\tau)}$ norm of \eqref{EqIntermediate} by
\begin{align*}
&\sum_{\sigma\prec\tau}
\int_{\bbR^d}G_t(z)\Vert P_\sigma{\sf\Gamma}_{(x-z)x}\tau\Vert_{L_x^{i_p(\tau):i_p(\sigma)}} \Vert \mcQ_t\big(x-z,{\sf\Pi}_{x-z}(\sigma)\big)\Vert_{L_x^{i_p(\sigma)}} dz\\
&\le\|{\sf\Gamma}\res  \tau\|_{w_{c'}} \|{\sf\Pi} \res  {\bf C}_{\prec\tau}\|_{w_c}
\sum_{\sigma\prec\tau}\int_{\bbR^d} G_t(z)w_{c+c'}(z)^{-1}\|z\|_\mfs^{r_{\varepsilon,p}(\tau)-r_{\varepsilon,p}(\sigma)}
t^{\frac{r_{\varepsilon,p}(\sigma)}\ell}
dz\\
&\lesssim \|{\sf\Gamma}\res  \tau\|_{w_{c'}} \|{\sf\Pi} \res  {\bf C}_{\prec\tau}\|_{w_c}
\sum_{\sigma\prec\tau}
t^{\frac{r_{\varepsilon,p}(\tau)-r_{\varepsilon,p}(\sigma)}\ell}t^{\frac{r_{\varepsilon,p}(\sigma)}\ell}
\lesssim
t^{\frac{r_{\varepsilon,p}(\tau)}\ell}
\|{\sf\Gamma}\res  \tau\|_{w_{c'}} \|{\sf\Pi} \res  {\bf C}_{\prec\tau}\|_{w_c}.
\end{align*}
This proves the statement of the lemma.
\end{Dem}

\ssk

\noindent \textbf{{Proof of Lemma \ref{lem:stationary} --}} A similar statement was proved by Linares, Otto, Tempelmayr \& Tsatsoulis in Proposition 4.6 of \cite{LOTT}. We follow here their reasoning. To simplify the notations we write ${\sf\Pi}_x,{\sf\Gamma},E_i$ for ${\sf\Pi}_x^{n;\varepsilon,\infty},{\sf \Gamma}^{n;\varepsilon,\infty},E_i^n$, respectively.

Note that $\bbE\big[\mcQ_t(x, {\sf\Pi}_x(\tau_i))\big]$ is independent of $x$ by stationarity of the law of the random noise under the translations of $\bbR^d$. One therefore has for any $0<s<t$
\begin{align*}
\partial_t\bbE\big[\mcQ_t(x,{\sf\Pi}_x(\tau_i))\big]
&=\int_{\bbR^d}\partial_t Q_{t-s}(x,y) \, \bbE\big[\mcQ_s(y,{\sf\Pi}_x(\tau_i))\big]dy   \\
&=\int_{\bbR^d}\partial_t Q_{t-s}(x,y) \, \bbE\big[\mcQ_s\big(y,{\sf\Pi}_x(\tau_i)-{\sf\Pi}_y(\tau_i)\big)\big]dy.
\end{align*}
Picking $s=\frac{t}2$ we apply Lemma \ref{LemRemainder} and have
\begin{align*}
\big\|\partial_t\bbE\big[\mcQ_t(x, {\sf\Pi}_x(\tau_i))\big]\big\|_{L_x^\infty(w_{c(m_{\bf B}+1)})}   
\lesssim t^{\frac{r_{\varepsilon,\infty}(\tau_i)}\ell-1}
\bbE\big[\|{\sf\Gamma}\res  \tau_i\|_{w_{cm_{\bf B}}}  \|{\sf\Pi} \res {\bf B}_{i-1}\|_{w_c}\big].
\end{align*}
Since $r_{\varepsilon,\infty}(\tau_i) < 0$ and
$$
\bbE\big[\mcQ_1(x, {\sf\Pi}_x(\tau_i))\big] = \bbE\big[\mcQ_1(0, {\sf\Pi}_0\tau)\big] =E_i
$$
by the stationarity of the law of the noise, we obtain the result by writing
\begin{align*}
\big\|\bbE\big[\mcQ_t(x, {\sf\Pi}_x(\tau_i))\big]&\big\|_{L_x^\infty(w_{c(m_{\bf B}+1)})}   \\
&\leq \big\|\bbE\big[\mcQ_1(x, {\sf\Pi}_x(\tau_i))\big]\big\|_{L_x^\infty(w_{c(m_{\bf B}+1)})} + \int_t^1\big\|\partial_s\bbE\big[\mcQ_s(x,{\sf\Pi}_x(\tau_i))\big]\big\|_{L_x^\infty(w_{c(m_{\bf B}+1)})} ds   \\
&\lesssim |E_i| + t^{\frac{r_{\varepsilon,\infty}(\tau_i)}\ell} \, \bbE\big[\|{\sf\Gamma}\res  \tau_i\|_{w_{cm_{\bf B}}}  \|{\sf\Pi} \res {\bf B}_{i-1}\|_{w_c}\big].
\end{align*}
This concludes the proof.   \hfill $\rhd$

\medskip

\noindent \textbf{{Proof of Lemma \ref{lem:Binega} --}} 
Here we omit the parameters $n,\infty,R$ similarly to the proof of Lemma \ref{lem:stationary}, but we leave the parameter $\varepsilon$ to make explicit the $\varepsilon$-dependence of the quantities $\|(\cdot)^{n;\varepsilon,\infty}\res\tau_i\|_{w_c}$.
Note that
$$
{\sf \Pi}_x^\varepsilon(\tau_i) = {\sf\Pi}_x^{2\varepsilon}(\tau_i)
$$
for any $\eps\in(0,\frac{\eps_0}2)$ by Lemma \ref{*epsilonuseless}, because ${\sf M}_x^{n;\eps,\infty}\vert_{\scV_\eps}$ is determined by ${\sf\Pi}^n$ and the structure of $(V_\eps^+,\Delta_{\eps,\infty}^+)$.
By using the integral formula
\begin{align*}
\mcQ_t\big(x, {\sf\Pi}_x^{\varepsilon}(\tau_i)\big) = \mcQ_1\big(x,{\sf\Pi}_x^{\varepsilon}(\tau_i)\big) - \int_t^1\partial_s \mcQ_s\big(x,{\sf\Pi}_x^{\varepsilon}(\tau_i)\big)ds,
\end{align*}
which is also used in the proof of Lemma \ref{lem:stationary}, we have
\begin{align*}
&\bbE\big[\|{\sf\Pi}^{2\varepsilon}\hspace{-0.03cm}:\hspace{0.07cm} \tau_i\|_{w_{c(m_{\bf B}+2)}}^q\big] 
= \bbE\bigg[\sup_{0<t\le1}t^{-\frac{r_{2\varepsilon,\infty}(\tau_i)}\ell q} \big\| \mcQ_t\big(x, {\sf\Pi}_x^{\varepsilon}(\tau_i)\big)\big\|_{L^\infty_x(w_{c(m_{\bf B}+2)})}^q\bigg]   \\
&\le\bbE\Big[\|\mcQ_1\big(x,{\sf\Pi}_x^{\varepsilon}(\tau_i)\big)\big\|_{L_x^\infty(w_{c(m_{\bf B}+2)})}^q\Big] \hspace{-0.03cm} + \hspace{-0.03cm} \bbE\left[\sup_{0<t\le1}\bigg|\int_t^1\hspace{-0.03cm}s^{-\frac{r_{2\varepsilon,\infty}(\tau_i)}\ell} \big\|\partial_s \mcQ_s\big(x, {\sf\Pi}_x^{\varepsilon}(\tau_i)\big)\big\|_{L^\infty_x(w_{c(m_{\bf B}+2)})} ds\bigg|^q\right]   \\
&\eqdef A_1 + A_2.
\end{align*}
We focus on the term $A_2$. By the semigroup property of $Q_s$ we have
\begin{align*}
\partial_s \mcQ_s\big(x,{\sf\Pi}_x^{\varepsilon}(\tau_i)\big) 
&= \int_{\bbR^d} (\partial_sQ)_{\frac{s}2}(x,y) \mcQ_{\frac{s}2}\big(y,{\sf\Pi}_x^{\varepsilon}(\tau_i)\big)dy   \\
&= \int_{\bbR^d} (\partial_sQ)_{\frac{s}2}(x,y) \mcQ_{\frac{s}2}\big(y,{\sf\Pi} _y^{\varepsilon}\tau_i)\big)dy + 
\int_{\bbR^d}(\partial_sQ)_{\frac{s}2}(x,y)\mcQ_{\frac{s}2}\big(y,{\sf\Pi}_x^{\varepsilon}(\tau_i)-{\sf\Pi} _y^{\varepsilon}(\tau_i)\big)dy.
\end{align*}
In the last line, the $L_x^\infty(w_{c(m_{\bf B}+1)})$ norm of the second term is bounded above  by
$$
s^{\frac{r_{\varepsilon,\infty}(\tau_i)}\ell - 1} \, (\star)_\varepsilon,\qquad
(\star)_\varepsilon\defeq
\|{\sf\Gamma}^{\varepsilon}\res  \tau_i\|_{w_{cm_{\bf B}}}  \|{\sf\Pi}^{\varepsilon} \res {\bf B}_{i-1}\|_{w_c}
$$
by Lemma \ref{LemRemainder}.
For the remaining term we can use Lemma \ref{lem:wyoung} and the bound $|\partial_tQ_t(x)|\lesssim t^{-1}G_t(x)$ to get
$$
\bigg\|\int_{\bbR^d} (\partial_sQ)_{\frac{s}2}(x,y) \mcQ_{\frac{s}2}\big(y, {\sf\Pi}_y^{\varepsilon}(\tau_i)\big)dy\bigg\|_{L_x^\infty(w_{c(m_{\bf B}+2)})} \lesssim s^{-\frac{|\mfs|}\ell \frac1q-1} \big\|\mcQ_{\frac{s}2}\big(y, {\sf\Pi}_y^{\varepsilon}(\tau_i)\big)\big\|_{L_y^q(w_{c(m_{\bf B}+2)})}.
$$
Write $q'$ for the conjugate exponent of $q$. If we choose sufficiently large $q$ such that 
$$
q>\frac{\ell+|\mfs|}{r_{\varepsilon,\infty}(\tau_i) - r_{2\varepsilon,\infty}(\tau_i)},
$$ 
then since $\big(\frac{r_{\varepsilon,\infty}(\tau_i) - r_{2\varepsilon,\infty}(\tau_i)}\ell-\frac{|\mfs|}\ell\frac1q-1\big)q' > -1$ we get from H\"older inequality in $t$ the estimate

\begin{align*}
\int_0^1s^{-\frac{r_{2\varepsilon,\infty}(\tau_i)}\ell} &\big\|\partial_s\mcQ_s\big(x,{\sf\Pi} _x^{\varepsilon}(\tau_i)\big)\big\|_{L_x^\infty(w_{c(m_{\bf B}+2)})}ds  \\
&\lesssim\int_0^1s^{\frac{r_{\varepsilon,\infty}(\tau_i) - r_{2\varepsilon,\infty}(\tau_i)}\ell - \frac{\vert\mfs\vert}\ell \frac1q-1}s^{-\frac{r_{\varepsilon,\infty}(\tau_i)}\ell} \big\|\mcQ_{\frac{s}2}\big(y,{\sf\Pi}_y^{\varepsilon}(\tau_i)\big)\big\|_{L_y^q(w_{c(m_{\bf B}+2)})}ds   \\
&\qquad+ (\star)_\varepsilon \int_0^1s^{\frac{r_{\varepsilon,\infty}(\tau_i) - r_{2\varepsilon,\infty}(\tau_i)}\ell-1}ds   \\
&\lesssim\bigg(\int_0^1s^{-\frac{r_{\varepsilon,\infty}(\tau_i)}\ell q} \, \big\|\mcQ_{\frac{s}2}\big(y,{\sf\Pi} _y^{\varepsilon}(\tau_i)\big)\big\|_{L_y^q(w_{c(m_{\bf B}+2)})}^qds\bigg)^{\frac1q} + (\star)_\varepsilon.
\end{align*}
We conclude from the estimate \eqref{eq:Binega} that
\begin{align*}
A_2&\lesssim\int_0^1s^{-\frac{r_{\varepsilon,\infty}(\tau_i)}\ell q} \, \bbE\big[\big\|\mcQ_{\frac{s}2}\big(y, {\sf\Pi}_y^{\varepsilon}(\tau_i)\big)\big\|_{L_y^q(w_{c(m_{\bf B}+2)})}^q\big]ds + \bbE[(\star)^q]   \\
&\lesssim |E_i|^q + \bbE\bigg[\sup_{\|h\|_{H}\le1} \|{\sf\Pi}^{n,h;\varepsilon,\infty} \hspace{-0.03cm}:\hspace{-0.03cm} \dt{\bf B}_i\|_{w_c}^q\bigg] + \bbE[(\star)^q].
\end{align*}
The proof of the estimate of $A_1$ is similar and left to the reader.  % \hfill $\rhd$

%  \vfill \pagebreak

%%-------------------------------------------------------------------%%
\subsection{Proofs of the Lemmas of Section \ref{sec:ind2}$\boldmath .$ \hspace{0.03cm}}
\label{SectionProofLemmaSection33}
%%-------------------------------------------------------------------%%

We prove the following two statements. We suppress the useless exponents $n,\varepsilon,p$ from the objects for more readability.

\ssk

\noindent \textbf{{Proof of Lemma \ref{lem:abstint} --}} 
Arguments yielding recursive representations of ${\sf g}_{yx}$ of this kind appear, for example, in Proposition 8.27 of \cite{Hai14} and Section 3 of \cite{Bru18}. Since the formulations there do not exactly match the one needed here, we provide a self-contained proof.
%The proof is classical and recalled here for completeness. 
Similar arguments can be found at Section 4 of \cite{BH21}.

Recall that ${\sf g}_x={\sf g}_x^{n;\eps,p}$ and ${\sf f}_x={\sf f}_x^{n;\eps,p}$ are defined by the formulas \eqref{EqConstructionGfromPi} (in terms of the renormalized model ${\sf\Pi}^{\xi_n,h_n,R_n}$ instead of the naive one ${\sf\Pi}^{\xi,h_n}$). 
In this proof, we extend ${\sf g}_x$ and ${\sf f}_x$ to all planted trees of the form $\mcI_l^{\bbK}(\sigma)$ by setting ${\sf g}_x(\mcI_l^{\bbK}(\sigma))={\sf f}_x(\mcI_l^{\bbK}(\sigma))=0$ if $r_{\eps,p}(\mcI_l^{\bbK}(\sigma))\le0$.
By applying the operator ${\sf g}_x\otimes{\sf g}_x^{-1}\otimes{\sf g}_x$ to the identity
\begin{align*}
\sum_{m\in\bbN^d}\frac{X^m}{m!}\otimes\Delta_{\eps,p}^+(\mcI_{k+m}^{\bbK}\sigma)
=\sum_{m\in\bbN^d}\frac{X^m}{m!}\otimes(\mcI_{k+m}^{\bbK}\otimes\id)\Delta_{\eps,p}(\sigma)
+\sum_{m,l\in\bbN^d}\frac{X^m}{m!}\otimes\frac{X^l}{l!}\otimes\mcI_{k+l+m}^{\bbK}\sigma,
\end{align*}
and using the fact $({\sf g}_x^{-1}\otimes{\sf g}_x)\Delta_{\eps,p}^+(\mcI_{k+m}^{\bbK}\sigma)=0$ and the binomial theorem, we have the explicit representation of ${\sf g}_x$
\begin{align}\label{Eq:g=fg}
{\sf g}_x(\mcI_k^{\bbK}\sigma)=({\sf f}_x\otimes{\sf g}_x)(\mcI_k^{\bbK}\otimes\id)\Delta_{\eps,p}(\sigma).
\end{align}

We prove the compatibility between ${\sf M}={\sf M}^{n;\varepsilon,p}$ and $\overline{\sf M}=\overline{\sf M}^{n;\varepsilon,p}$ for the abstract integration map $\overline{\mcI}=\overline{\mcI}^{\eps,p}$.
We denote by
$$
\mcJ(x)\sigma\defeq\sum_{|k|_\mfs<r_{\eps,p}(\tau)+\beta_0}\frac{X^k}{k!}\partial^k\mcK(x,{\sf\Pi}_x\tau)
=\sum_{k\in\bbN^d}\frac{X^k}{k!}{\sf f}_x(\mcI_k^{\bbK}(\sigma))
\qquad(\sigma\in{\bf B}_{i-1}\cup\dt{\bf B}_{i-1})
$$
the linear map $\mcJ(x):U_{i-1}\to\spa\{X^k\}$ introduced in Definition \ref{defTbarT}.
As an application of \eqref{Eq:g=fg} with $k=0$, we have the formula
\begin{align*}
{\sf g}_{yx}\big((\overline{\mcI}+\mcJ(x))\sigma\big)={\sf f}_y\big(\overline{\mcI}({\sf\Gamma}_{yx}\sigma)\big).
\end{align*}
Indeed one has
\begin{align*}
{\sf f}_y\big(\overline{\mcI}({\sf\Gamma}_{yx}\sigma)\big)&=({\sf f}_y\otimes{\sf g}_y\otimes{\sf g}_x^{-1})(\mcI_0^{\bbK}\otimes\id\otimes\id)(\id\otimes\Delta_{\eps,p}^+)\Delta_{\eps,p}(\sigma)   \\
&=({\sf f}_y\otimes{\sf g}_y\otimes{\sf g}_x^{-1})(\mcI_0^{\bbK}\otimes\id\otimes\id)(\Delta_{\eps,p}\otimes\id)\Delta_{\eps,p}(\sigma)
=({\sf g}_y\otimes{\sf g}_x^{-1})(\mcI_0^{\bbK}\otimes\id)\Delta_{\eps,p}(\sigma)   \\
&=({\sf g}_y\otimes{\sf g}_x^{-1})\bigg(\Delta_{\eps,p}^+(\mcI_0^{\bbK}\tau) - \sum_{k\in\bbN^d}\frac{X^k}{k!}\otimes \mcI_k^{\bbK}\sigma\bigg)   \\
&={\sf g}_{yx}(\mcI_0^{\bbK}\sigma)+\sum_{k,l\in\bbN^d}\frac{y^k}{k!}\frac{(-x)^l}{l!} \, {\sf f}_x(\mcI_{k+l}^{\bbK}\sigma)\\
&={\sf g}_{yx}(\mcI_0^{\bbK}\sigma)+\sum_{m\in\bbN^d}\frac{(y-x)^m}{m!} \, {\sf f}_x(\mcI_{m}^{\bbK}\sigma)={\sf g}_{yx}\big((\overline{\mcI}+\mcJ^{\sf M}(x))\sigma\big).
\end{align*}
Thus for any $x,y,z\in\bbR^d$ we have
\begin{align*}
{\sf g}_{zy} \Big\{\overline{\sf \Gamma}_{yx}(\overline{\mcI}+\mcJ(x))(\sigma) - \big(\overline{\mcI}+\mcJ(y)\big)\overline{\sf \Gamma}_{yx}(\sigma)\Big\}
&={\sf g}_{zx}\big((\overline{\mcI}+\mcJ(x))\sigma\big) - {\sf f}_z\big(\overline{\mcI}{\sf \Gamma}_{zy}{\sf\Gamma}_{yx}(\sigma)\big)   \\
&={\sf f}_z(\overline{\mcI}{\sf\Gamma}_{zx}\sigma)-{\sf f}_z(\overline{\mcI}{\sf\Gamma}_{zx}\sigma)=0.
\end{align*}
Since the quantity inside ${\sf g}_{zy}$ is in the subspace spanned by $\{X^{k}\}_{k}$, on which ${\sf g}_{zy}$ is injective, we have the compatibility
$$
\overline{\sf \Gamma}_{yx}\big((\overline{\mcI}+\mcJ(x))(\sigma)\big)
=\big(\overline{\mcI}+\mcJ(y)\big)\overline{\sf \Gamma}_{yx}(\sigma).
$$   %  

Finally we prove Lemma \ref{lem:abstint}. Write $f_\tau=f_\tau^{n;\eps,p}$ for simplicity.
We define
\begin{align*}
\mcN\big(x;f_\tau,{\sf\Pi}^n(\tau)\big) &\defeq \sum_{|k|_\mfs<r_{\eps,p}(\tau)+\beta_0}\frac{X^k}{k!} \, \partial^k\mcK\big(x,{\sf\Pi}^n\tau-{\sf\Pi}_x^{n;\eps,p}(f_\tau^{n;\eps,p}(x))\big)   \\
&=\sum_{|k|_\mfs<r_{\eps,p}(\tau)+\beta_0}\frac{X^k}{k!} \, \partial^k\mcK(x,{\sf\Pi}_x^{n;\eps,p}\tau)
=\sum_{k\in\bbN^d}\frac{X^k}{k!} \, {\sf f}_x(\mcI_k^{\bbK}\tau),   \\
\big(\overline{\mcK}f_\tau\big)(x)
&\defeq\overline{\mcI}\big(f_\tau(x)\big)+\mcJ(x)f_\tau(x) + \mcN\big(x;f_\tau,{\sf\Pi}^n(\tau)\big),
\end{align*}
the functions introduced in Definition \ref{defTbarT}.
Since $f_\tau(x)=(\id\otimes{\sf g}_x)\Delta_{\eps,p}\tau-\tau$, by using \eqref{Eq:g=fg} we have
\begin{align*}
(\overline{\mcK}f_\tau)(x)-\overline{\mcI}\big(f_\tau(x)\big)
&=\sum_{k\in\bbN^d}\frac{X^k}{k!}{\sf f}_x\big(\mcI_k^{\bbK}f_\tau(x)\big)
+\sum_{k\in\bbN^d}\frac{X^k}{k!}{\sf f}_x(\mcI_k^{\bbK}\tau)\\
&=\sum_{k\in\bbN^d}\frac{X^k}{k!}({\sf f}_x\otimes{\sf g}_x)(\mcI_k^{\bbK}\otimes\id)\Delta_{\eps,p}\tau
=\sum_{k\in\bbN^d}{\sf g}_x(\mcI_{k}^{\bbK}\tau) \, \frac{X^k}{{k}!}.
\end{align*}
This formula further transformed into
\begin{align*}
(\overline{\mcK}f_\tau)(x) 
&=\overline{\mcI}\big(f_\tau(x)\big)+\sum_{k\in\bbN^d}{\sf g}_x(\mcI_{k}^{\bbK}\tau) \, \frac{X^k}{{k}!}
=(\mcI_0^{\bbK}\otimes{\sf g}_x)\Delta_{\eps,p}\tau-\mcI_0^{\bbK}\tau+\sum_{k\in\bbN^d}{\sf g}_x(\mcI_{k}^{\bbK}\tau) \, \frac{X^k}{{k}!}\\
&=(\id\otimes{\sf g}_x)\Delta_{\eps,p}^+\mcI_0^{\bbK}\tau-\mcI_0^{\bbK}\tau,
\end{align*}
which is a form similar to the definition of $f_\tau$ in Lemma \ref{*argumentforreconst} where $\tau$ is replaced by $\mcI_0^{\bbK}\tau$.
Thus we have the result in a similar argument to the proof of Lemma \ref{*argumentforreconst}.
\hfill $\rhd$

\ssk

\noindent \textbf{{Proof of Lemma \ref{lem:Pitog} --}} We consider the case where $\tau\in\dt{\bf B}_i$ and let the reader treat the case where $\tau\in{\bf B}_i$. We saw in Lemma \ref{*argumentforreconst} that $f_\tau$ belongs to the space $D^{(r_{\varepsilon,p}(\tau),p)}(U_{i-1};{\sf\Gamma})_{w_{cm_{\bf B}}}$ 
with norms
$$
\tri f_\tau\tri_{(r_{\varepsilon,p}(\tau),p);w_{cm_{\bf B}}}^{\sf\Gamma}
\lesssim \big(1+\|{\sf g}\res {\bf W}_{i-1,\varepsilon,p}^+\|_{w_c}\big)^{m_{\bf B}}.
$$
Therefore we have from Theorem \ref{thm:MS} the estimate
\begin{align*}
&\tri \overline{\mcK} f_\tau\tri_{(r_{\varepsilon,p}(\tau) + \beta_0,p);w_{c(2m_{\bf B}+1)}} \\
&\lesssim\|{\sf\Pi}\res {\bf B}_{i-1}\cup\dt{\bf B}_{i-1}\|_{w_c} \big(1+\|{\sf\Gamma}\res {\bf B}_{i-1}\cup\dt{\bf B}_{i-1}\|_{w_{cm_{\bf B}}}\big) \tri f_\tau\tri_{(r_{\varepsilon,p}(\tau),p);w_{cm_{\bf B}}}^{\sf\Gamma} + \lb{\sf\Pi}\tau\rb_{(r_{\varepsilon,p}(\tau),p);w_{c}}^{{\sf\Pi}, f_\tau}   \\
&\lesssim\|{\sf\Pi}\res {\bf B}_{i-1}\cup\dt{\bf B}_{i-1}\|_{w_c} \big(1+\|{\sf g}\res {\bf W}_{i-1,\varepsilon,p}^+\|_{w_c}\big)^{2m_{\bf B}} + \|{\sf\Pi}\res\tau\|_{w_c},
\end{align*}
on which one reads the conclusion.   %  \hfill $\rhd$

%\ssk

%%-----------------------------------------------------------------%%
\subsection{Proofs of the lemmas of Section \ref{sec:ind3}$\boldmath .$ \hspace{0.03cm}}
\label{SectionProofLemmas34}
%%-----------------------------------------------------------------%%

Recall $\tau\in \overset{\text{\tiny$\mybullet$}}{{\bf B}}_{i+1}$. We have two statements to prove.

\ssk

\noindent \textbf{{Proof of Lemma \ref{lem:Lpfx} --}} Recall from \eqref{EqDefnCaracterH} the definition of $\lambda_x^{n;\varepsilon,p}$. The proof leads to the estimate of the character ${\sf f}_x^{n;\eps,p}$ defined by \eqref{EqConstructionGfromPi} where the naive one ${\sf\Pi}^{\xi,h_n}$ is replaced by the renormalized one ${\sf\Pi}^{\xi_n,h_n,R_n}$. This is a content of Lemma 5.6 of \cite{SemigroupMasato}, but the proof is recalled here for completeness. Pick $\sigma\in\dt{\bf B}_i$ with $\mu=\mcI_k^{\bbK}(\sigma)$. Since the behavior of ${\sf \Pi}^{n;\varepsilon,p}$ as a `step function' of $p$, we can write
$$
\lambda_x^{n;\varepsilon,p}(\mu) \defeq {\bf1}_{r_{\varepsilon,p}(\mu)\le0<r_{\varepsilon,2}(\mu)} \, {\sf f}_x^{n;\varepsilon,q}(\mu).
$$
We can write the ${\sf f}_x^{n;\varepsilon,q}(\mu)$ by the form
$$
\int_0^1\partial^{k}\mcK_t\big(x, {\sf\Pi}_x^{n;\varepsilon,q}(\sigma)\big) \, dt
=\int_0^1\int_{\bbR^d}\partial_x^{k}K_{\frac{t}2}(x,y) \mcQ_{\frac{t}2}\big(y,{\sf\Pi}_x^{n;\varepsilon,q}(\sigma)\big)\,dydt.
$$
By Lemma \ref{LemRemainder}, we can replace ${\sf\Pi}_x^{n;\varepsilon,q}\sigma$ by ${\sf\Pi}_y^{n;\varepsilon,q}\sigma$ with an error 
$$
\|{\sf\Gamma}^{n;\varepsilon,q}\res \sigma\|_{w_{cm_{\bf B}}}  \|{\sf\Pi}^{n;\varepsilon,q} \res {\bf B}_{i-1}\cup\dt{\bf B}_{i-1}\|_{w_c},
$$ 
where we note that $r_{\varepsilon,q}(\mu)>r_{\varepsilon,p_\varepsilon(\mu)}(\mu)=0$ and the $t$-dependent factor $t^{\frac{r_{\varepsilon,q}(\mu)}\ell-1}$ is integrable.
After the replacement, we obtain the statement as a direct consequence of the size estimate \eqref{EqDefnSizePiTau} satisfied by $\mcQ_{\frac{t}2}(y,{\sf\Pi}_y^{n;\varepsilon,q}(\sigma))$.   % \hfill $\rhd$

\ssk

\noindent \textbf{{Proof of Lemma \ref{lem:goal} --}} To lighten the notations we suppress in this proof the exponent $n$ from the models and their associated quantities. 
By the semigroup property of $Q_t$ we have
\begin{align*}
\mcQ_t\big(x, {\sf\Pi}_x^{\eps,p}(\tau)\big)
&= \int_{\bbR^d}Q_{\frac{t}2}(x,y) \mcQ_{\frac{t}2}\big(y,{\sf\Pi} _x^{\eps,p}(\tau)\big)dy   \\
&= \int_{\bbR^d}Q_{\frac{t}2}(x,y) \mcQ_{\frac{t}2}\big(y,{\sf\Pi} _y^{\eps,p}(\tau)\big)dy + \int_{\bbR^d}Q_{\frac{t}2}(x,y) \mcQ_{\frac{t}2}\big(y,{\sf\Pi} _x^{\eps,p}(\tau) - {\sf\Pi} _y^{\eps,p}(\tau)\big)dy\\
&\eqdef B_1 + B_2.
\end{align*}
We can apply Lemma \ref{LemRemainder} to bound above the $L_x^p(w_{c(m_{\bf B}+1)})$ norm of $B_2$ by
\begin{equation*}
t^{\frac{r_{\varepsilon,p}(\tau)}\ell} \|{\sf\Gamma}^{\varepsilon,p} \res \tau\|_{w_{cm_{\bf B}}} \|{\sf\Pi}^{\varepsilon,p}\res {\bf B}_i\cup\dt{\bf B}_i\|_{w_c}.
\end{equation*}
For $B_1$, since $\eps\notin {\sf J}_p$, we replace ${\sf\Pi}^{\eps,p}$ by ${\sf\Pi}^{\eps',p}$ for any $\eps'\in(0,\eps)$ such that $(\eps',\eps)\cap {\sf J}_p=\varnothing$ -- recall by Lemma \ref{*epsilonuseless}.
We write \eqref{EqSyntheticDecompositionFormulaDerivativePi} in synthetic form
\begin{equation} \label{EqSyntheticDecompositionFormula}	
{\sf\Pi}_y^{\eps',p}(\tau) = {\sf\Pi} _y^{\eps',2}(\tau) + \sum_{\tau_1,\tau_2}\lambda_y^{\eps',p}(\tau_1) \, {\sf\Pi} _y^{\eps'}(\tau_2)
\end{equation}	
with $\tau_1\in{\bf W}_{i,\eps',p}^+$, $\tau_2\in{\bf B}_i$, $r_{\eps',p}(\tau_1) + r_{\eps',p}(\tau_2) = r_{\eps',p}(\tau)$, and $p\ge p_{\eps'}(\tau_1)$. By by Lemma \ref{*epsilonuseless}, ${\sf\Pi}_y^{\eps',2}(\tau)={\sf\Pi}_y^{\eps,2}(\tau)$ and ${\sf\Pi}_y^{\eps'}(\tau_2)={\sf\Pi}_y^\eps(\tau_2)$. Note the elementary, yet fundamental, identities
\begin{align*}
r_{\varepsilon,2}(\tau) + \vert\frak{s}\vert\Big(\frac{1}{p}-\frac{1}{2}\Big) &= r_{\varepsilon,p}(\tau), \\
r_{\varepsilon,p}(\tau_2) + \vert\frak{s}\vert\Big(\frac{1}{p}-\frac{1}{q_1}\Big) &= r_{\varepsilon,p}(\tau) - r_{\varepsilon,q_1}(\tau_1).
\end{align*}
We can choose $\eps'<\eps$ so that 
$$
\lfloor p_{\eps'}(\tau_1)\rfloor_{{\sf I}{\eps'}}\leq p_\eps(\tau_1)<p_{\eps'}(\tau_1),\qquad
\eps'\notin\bigcup_{\tau_1}{\sf J}_{p_\eps(\tau_1)}.
$$
(The first condition is checked as follows. Since $\eps'<\eps$, we have $p_\eps(\tau_1)<p_{\eps'}(\tau_1)$. For any $\sigma\in{\bf W}_{\eps',p}^+$ such that $p_{\eps'}(\sigma)<p_{\eps'}(\tau_1)$, we have $p_{\eps}(\sigma)<p_{\eps}(\tau_1)$ by Lemma \ref{*epsilonuseless}. Thus $p_{\eps'}(\sigma)<p_\eps(\tau_1)$ if $\eps'$ is sufficiently close to $\eps$.)
We then have from Lemmas \ref{LemContinuousExtension} and \ref{lem:Lpfx} an upper bound for $\|B_1\|_{L_x^p(w_{c(m_{\bf B}+2)})}$ of the form
\begin{align*}
&\bigg\|\int_{\bbR^d}G_t(x-y) \big|\mcQ_{\frac{t}2}\big(y,{\sf\Pi} _{y}^{\eps,2}(\tau)\big)\big| dy
+ \sum_{\tau_1,\tau_2}\int_{\bbR^d}G_t(x-y)|\lambda_y^{\eps',p}(\tau_1)| \big|\mcQ_{\frac{t}2}\big(y,{\sf\Pi} _y^{\eps}(\tau_2)\big)\big| dy\bigg\|_{L_x^p(w_{c(m_{\bf B}+2)})} 
\\
&\lesssim t^{\frac{|\mfs|}\ell (\frac1p-\frac12)} \big\|\mcQ_{\frac{t}2}\big(y,{\sf\Pi}_{y}^{\eps,2}(\tau)\big)\big\|_{L_y^2(w_c)}   \\
&\qquad+ \sum_{\tau_1,\tau_2}t^{\frac{|\mfs|}\ell(\frac1p-\frac1{p_{\eps}(\tau_1)})}\|\lambda_y^{\eps',p}(\tau_1)\|_{L_y^{p_\eps(\tau_1)}(w_{c(m_{\bf B}+1)})} \big\|\mcQ_{\frac{t}2}\big(y,{\sf\Pi} _y^{\eps}(\tau_2)\big)\big\|_{L_y^\infty(w_c)}   
\\
&\lesssim t^{\frac{|\mfs|}\ell (\frac1p-\frac12)}t^{\frac{r_{\eps,2}(\tau)}\ell} \|{\sf\Pi}^{\varepsilon,2}\res\tau\|_{w_c}   \\
&\quad+ \sum_{\tau_1,\tau_2}t^{\frac{|\mfs|}\ell(\frac1p-\frac1{p_{\eps}(\tau_1)})}t^{\frac{r_{\eps,\infty}(\tau_2)}\ell} \big(1 + \Vert {\sf M}^{\eps',p_{\eps}(\tau_1)}\Vert_{\textbf{\textsf{M}}(\scW_{i,\eps',p_{\eps}(\tau_1)})_{w_c}}\big)^{m_{\bf B}+1} \|{\sf\Pi}^{\varepsilon,\infty}\res {\bf B}_i\|_{w_c}   
\\
&=
t^{\frac{r_{\eps,p}(\tau)}\ell}\|{\sf\Pi}^{\varepsilon,2}\res\tau\|_{w_c} 
+ \sum_{\tau_1,\tau_2}t^{\frac{r_{\eps,p}(\tau)-r_{\eps,p_{\eps}(\tau_1)}(\tau_1)}\ell} \big(1 + \Vert {\sf M}^{\eps',p_{\eps}(\tau_1)}\Vert_{\textbf{\textsf{M}}(\scW_{i,\eps',p_{\eps}(\tau_1)})_{w_c}}\big)^{m_{\bf B}+1} \|{\sf\Pi}^{\eps,\infty}\res{\bf B}_i\|_{w_c}  
 \\
&\lesssim 
t^{\frac{r_{\eps,p}(\tau)}\ell}\Bigg(1+
\Vert {\sf M}^{\eps,2}\Vert_{\textbf{\textsf{M}}(\scW_{i,\eps,2})_{w_c}}
+\sum_{\tau_1}\Vert {\sf M}^{\eps',p_{\eps}(\tau_1)}\Vert_{\textbf{\textsf{M}}(\scW_{i,\eps',p_{\eps}(\tau_1)})_{w_c}}+\Vert {\sf M}^{\eps,\infty}\Vert_{\textbf{\textsf{M}}(\scV_{i,\eps})_{w_c}}\Bigg)^{m_{\bf B}+2},
\end{align*}
which implies the result.    \hfill $\rhd$

%--------------------------------------------------------------------------------------------------------------------%
\section{Proofs of Theorem \ref{ThmContinuity} and Theorem \ref{ThmLipschitzContinuity}}
\label{SectionThmContinuity}
%--------------------------------------------------------------------------------------------------------------------%

%%--------------------------------------------------------------------------------------------------%%
\subsection{Proof of Theorem \ref{ThmContinuity}$\boldmath .$ \hspace{0.03cm}} 
\label{SectionThmContinuity1}
%%--------------------------------------------------------------------------------------------------%%
Our proof of Theorem \ref{ThmContinuity} follows the pattern of proof of the corresponding statement in Tempelmayr's recent work \cite{Tempelmayr}. Since this part is independent of the parameters $\varepsilon,p,c$, we suppress them from the corresponding objects. We write for instance $\textbf{\textsf{M}}(\scW)$ for $\textbf{\textsf{M}}(\scW_{\varepsilon,p})_{w_c}$.

Assume now that we are given a sequence $(\bbP_j)_{j\geq 0}$ of stationary Borel probability measures on $\Omega=C^{r_0,Q}(w_c)$ that all satisfy the spectral gap inequality \eqref{EqSGInequality} with the same constant $C$. Also, assume that $\bbP_j$ tend weakly to a probability $\bbP_\infty$. Then $\bbP_\infty$ also satisfies the spectral gap inequality with the same $C$. It is easily checked for any bounded cylindrical functions $F:\Omega\to\bbR$ of the form
$$
F(\omega) = f\big(\varphi_1(\omega),\dots,\varphi_m(\omega)\big),
$$
where $\varphi_1,\dots,\varphi_m\in \Omega^*$ and $f\in C^\infty(\bbR^n)$ and extended to all $F$ by a density argument -- see e.g. Lemma 2.23 in Hairer \& Steele's work \cite{HS23} for the density argument. 
For each $j\in\bbN\cup\{\infty\}$ and $n\in\bbN$, let ${\sf M}_j^n(\omega)={\sf M}^{\xi_n(\omega),R_j^n}$ be the BPHZ model on $\scW$ associated with the random variable $\xi_n(\omega)=\varrho_n*\omega$ and the unique BPHZ-type preparation map $R_j^n$ defined by \eqref{EqdeterminesBHZ}. We also denote by ${\sf M}_j=\lim_{n\to\infty}{\sf M}_j^n$ the $L^q(\bbP_j)$-limit for any $q<\infty$ ensured by Theorem \ref{thm:main}. From its proof it is obvious that the quantities
$$
\sup_{n\in\bbN}\bbE_j\big[\|{\sf M}_j^{n}\|_{\textbf{\textsf{M}}(\scW)}^{q}\big]
$$
have the $j$-uniform upper bound and the convergence
$$
\lim_{n\to\infty}\bbE_j\big[\|{\sf M}_j^{n} \hspace{-0.03cm}:\hspace{-0.03cm} {\sf M}_j\|_{\textbf{\textsf{M}}(\scW)}^{q}\big]=0
$$
is $j$-uniform. By Skorohod representation theorem there is a probability space $(E, \mathbf{\mcE},\bbQ)$ and random variables
$$
\widetilde{\xi}_j:E\to\Omega
$$
for each $j\in\bbN\cup\{\infty\}$ such that the $\bbQ$-law of $\widetilde{\xi}_j$ is equal to $\bbP_j$ and $\widetilde{\xi}_j$ converges $\bbQ$-almost surely to $\widetilde{\xi}_\infty$ as $j\to\infty$.
Once we define the random variables
\begin{equation*} \begin{split}
&\widetilde{\sf M}_j^n\defeq{\sf M}_j^n(\widetilde{\xi}_j) : E \rightarrow \textbf{\textsf{M}}(\scW),  \qquad  \widetilde{\sf M}_j\defeq{\sf M}_j(\widetilde{\xi}_j) : E \rightarrow \textbf{\textsf{M}}(\scW),
\end{split} \end{equation*}
we have that
\begin{enumerate}
	\item[(1)] for each $j\in\bbN\cup\{\infty\}$, the $\bbQ$-law of $\widetilde{\sf M}_j^n$ (resp. $\widetilde{\sf M}_j$) is equal to the $\bbP_j$-law of ${\sf M}_j^n$ (resp. ${\sf M}_j$), and $\widetilde{\sf M}_j^n$ converges to $\widetilde{\sf M}_j$ in $L^q(E,\bbQ;\textbf{\textsf{M}}(\scW))$ for any $q<\infty$ as $n$ goes to $\infty$,
	\item[(2)] for each $n\in\bbN$, the random variables $\widetilde{\sf M}_j^n$ converges $\bbQ$-almost surely to $\widetilde{\sf M}_\infty^n$ in $\textbf{\textsf{M}}(\scW)$ as $j\to\infty$.
\end{enumerate}
Point $(2)$ follows from the facts that $\widetilde{\sf M}_j^n$ is a continuous function of the smooth function $\widetilde{\xi}_j^n=\varrho_n*\widetilde{\xi}_j$ and the $\bbQ$-almost sure convergence of $\widetilde{\xi}_j^n$ to $\widetilde{\xi}_\infty^n$ in the space of smooth functions as $j\to\infty$. Denoting by ${\sf E}_\bbQ[\cdot]$ the expectation operator associated to $\bbQ$, one has for every $j,k\in\bbN$ and $n\in\bbN$,
\begin{equation*} \begin{split}
{\sf E}_\bbQ\big[\Vert \widetilde{\sf M}_j \hspace{-0.03cm}:\hspace{-0.03cm} \widetilde{\sf M}_\infty\Vert_{\textbf{\textsf{M}}(\scW)}^q\big]
\lesssim {\sf E}_\bbQ\big[\Vert \widetilde{\sf M}_j \hspace{-0.03cm}:\hspace{-0.03cm} \widetilde{\sf M}_j^n \Vert_{\textbf{\textsf{M}}(\scW)}^q\big]
+ {\sf E}_\bbQ\big[\Vert \widetilde{\sf M}_j^n \hspace{-0.03cm}:\hspace{-0.03cm} \widetilde{\sf M}_\infty^n \Vert_{\textbf{\textsf{M}}(\scW)}^q\big]
+ {\sf E}_\bbQ\big[\Vert \widetilde{\sf M}_\infty^n \hspace{-0.03cm}:\hspace{-0.03cm} \widetilde{\sf M}_\infty \Vert_{\textbf{\textsf{M}}(\scW)}^q\big].
\end{split} \end{equation*}
By point (1), the first and third terms in the right hand side are bounded above by $(j,k)$-independent constant $C_n$ which vanishes as $n$ goes to $\infty$. 
By point (2) and the $j$-uniformly integrability of $\Vert\widetilde{\sf M}_j^n\Vert_{\textbf{\textsf{M}}(\scW)}$,
the second term vanishes as $j\to\infty$ for each $n$ by Vitali convergence theorem.
Therefore we have
$$
\lim_{j\to\infty}{\sf E}_\bbQ\big[\Vert \widetilde{\sf M}_j \hspace{-0.03cm}:\hspace{-0.03cm} \widetilde{\sf M}_\infty\Vert_{\textbf{\textsf{M}}(\scW)}^q\big]=0,
$$
which implies that the $\bbP_j$-law of ${\sf M}_j$ converges to the $\bbP_\infty$-law of ${\sf M}_\infty$ as $j\to\infty$.

%\medskip

%%---------------------------------------------------------------------------------------------%%
\subsection{Theorem \ref{ThmLipschitzContinuity} and random Fourier series$\boldmath .$ \hspace{0.03cm}} 
\label{SectionThmContinuity2}
%%---------------------------------------------------------------------------------------------%%

Let $(\xi_j)_{j\in\bbN\cup\{\infty\}}$ be a sequence of $C^{r_0,Q}(w_c)$-valued random variables defined on the same probability space $(\Omega,\mcF,\bbP)$. Assume that all $\xi_j$ satisfy the spectral gap inequality \eqref{EqSGInequality} with the same constant $C$ and that $(\xi_j)_{j\in\bbN}$ converges to $\xi_\infty$ in $L^q(\bbP)$ for any $q\in[1,\infty)$.
For each $j\in\bbN\cup\{\infty\}$, denote by ${\sf M}_j={\sf M}(\xi_j)$ the limit BPHZ model ensured by Theorem \ref{thm:main}.
Theorem  \ref{ThmLipschitzContinuity} claims that, for any $q\in[1,\infty)$, there exists $r\in(q,\infty)$ such that we have the qualitative convergence
\begin{equation}\label{*eq:QuantitativeConvergence}
\bbE\big[\Vert {\sf M}_j \res {\sf M}_\infty\Vert_{\textbf{\textsf{M}}(\scW)}^q\big]^{1/q}
\lesssim\bbE[\Vert\xi_j-\xi_\infty\Vert_{C^{r_0,Q}(w_c)}^r]^{1/r}.
\end{equation}
It is a direct consequence of the local Lipschitz estimates from Section \ref{sec:ind4}, where $({\sf M}^{n;\eps,p},{\sf M}^{m;\eps,p})$ is replaced with $({\sf M}_j,{\sf M}_\infty)$. As for the first step, we have the deterministic estimates
$$
\Vert{\sf\Pi}_j,{\sf\Pi}_\infty\res\odot\Vert_{w_c}=0,\qquad
\Vert{\sf\Pi}_j,{\sf\Pi}_\infty\res\ocircle\Vert_{w_c}\le\Vert\xi_j-\xi_\infty\Vert_{C^{r_0,Q}(w_c)},
$$
since $({\sf\Pi}_j)_x(\odot)=h$ and $({\sf\Pi}_j)_x(\ocircle)=\xi_j$ by definition.
Since the quantities
$$
\bbE\big[\|{\sf M}_j\|_{\textbf{\textsf{M}}(\scW)}^{q}\big]
$$
have the $j$-uniform upper bound, by repeating each step described in Section \ref{sec:ind4}, we can establish the estimate \eqref{*eq:QuantitativeConvergence} inductively.

\ssk

Random Fourier series provide an interesting class of random distributions where our main results are relevant. Here is an example where we work with $\bbZ^2$-periodic distributions on $\bbR^2$. The heat semigroup $(Q_t)_{t>0}$ that we use to define the functional setting is that of the usual heat operator on $\bbR^2$ and we write $B^\alpha_{pq}$ for $B^{\alpha,Q}_{pq}(w_0)$ and $C^\alpha$ for the corresponding Besov-H\"older space. Set 
$$
e_k(x) \defeq e^{ik\cdot x}
$$
for any $k\in\bbN^2$ and $x\in\bbT^2$. For any sequence $(f^k)_{k\in\bbN^2}$ of real numbers the Sobolev norms of the distributions 
$$
f(\cdot) = \sum_{k\in\bbN^2} f^k e_k(\cdot)
$$
satisfy for any $\alpha\in\bbR$
$$
\sum_{k\in\bbN^2} \vert k\vert^{2\alpha} (f^k)^2 \lesssim \Vert f\Vert_{B^\alpha_{22}}^2 \lesssim \sum_{k\in\bbN^2} \vert k\vert^{2\alpha} (f^k)^2.
$$
Pick two deterministic sequences $(f_1^k)_{k\in\bbN^2}, (f_2^k)_{k\in\bbN^2}$ such that
$$
\sum_{k\in\bbN^2}  \vert k\vert^{2r_0} \big((f_1^k)^2 + (f_2^k)^2\big) < \infty.
$$
Let $(\gamma_1^k)_{k\in\bbN^2}, (\gamma_2^k)_{k\in\bbN^2}$ be two independent sequences of independent standard normal variables. Then for $i\in\{1,2\}$ the sum
$$
\xi_i(\cdot) \defeq \sum_{k\in\bbN^2} \gamma_i^k f_i^k e_k(\cdot)
$$
converges in $L^q(\Omega,C^{r_0-\eps})$ for any $1\leq q<\infty$ and $\eps>0$, and one has from hypercontractivity
$$
\big\Vert \Vert \xi_1 - \xi_2\Vert_{B^{r_0-2(\frac12-\frac1q)}_{qq}} \big\Vert_{L^q(\bbP)} \lesssim_q \big\Vert \Vert \xi_1 - \xi_2\Vert_{B^{r_0}_{22}} \big\Vert_{L^2(\bbP)}^{\frac{q}{2}} \lesssim_q \left(\sum_{k\in\bbN^2} (f_1^k - f_2^k)^2\vert k\vert^{2r_0} \right)^{\frac{1}{2}}
$$
In particular if $f_2^k = f_1^k {\bf 1}_{\vert k\vert \leq n}$ for some $n\geq 1$, and one writes $\xi_1^{(n)}$ for the corresponding distribution $\xi_2$, then one has 

$$
\big\Vert \Vert \xi_1 - \xi_1^{(n)}\Vert_{C^{r_0-\eps}} \big\Vert_{L^q(\bbP)} \lesssim_q \Bigg(\sum_{\vert k\vert >n} (f_1^k)^2\vert k\vert^{2r_0} \Bigg)^{\frac{1}{2}}.
$$
Denote by ${\sf M}(\xi_1)$ and ${\sf M}(\xi_1^{(n)})$ the BPHZ models associated with $\xi_1$ and $\xi_1^{(n)}$ respectively.

\ssk

\begin{cor}
Pick $1\leq q<\infty$. The models ${\sf M}(\xi_1^{(n)})$ converge to ${\sf M}(\xi_1)$ in $L^q\big(\Omega ; \textsf{\textbf{M}}(\scW)\big)$ as $n$ goes to $\infty$ at speed $\big(\sum_{\vert k\vert >n} (f_1^k)^2\vert k\vert^{2r_0} \big)^{1/2}$.
\end{cor}

\bigskip

\appendix

%-------------------------------------------------------------------------------------------------------------------------%
\section{Reconstruction and multilevel Schauder estimates in regularity-integrability structures}
\label{SectionAppendix}
%-------------------------------------------------------------------------------------------------------------------------%

We give in this section a summary of the results on modelled distributions over a regularity-integrability structure proved in the companion work \cite{SemigroupMasato}. Given $1\leq p\leq q\leq \infty$ we define the exponent $p\hspace{-0.03cm}:\hspace{-0.03cm}q\in[p,\infty]$ from the relation
$$
\frac{1}{p:q} + \frac{1}{q} = \frac{1}{p}.
$$

\ssk

\begin{defn}
Let $\scT=(A,T,G)$ be a regularity-integrability structure of regularity $r_0\le0$.
For any $c>0$, we define $\textbf{\textsf{M}}(\scT)_{w_c}$ as the set of pairs $\big(\{\Pi_x\}_{x\in\bbR^d}, \{\Gamma_{xy}\}_{x,y\in\bbR^d}\big)$ with the following properties.
\begin{itemize}
\item[(1)]
{\sl (Algebraic conditions)} The maps $\Pi_x:T\to C^{r_0,Q}(w_c)$ are linear continuous, $\Gamma_{xy}\in G$, and $\Pi_x\Gamma_{xy}=\Pi_y$, $\Gamma_{xx}=\text{\rm Id}$, and $\Gamma_{xy}\Gamma_{yz}=\Gamma_{xz}$ for any $x,y,z\in\bbR^d$.
\item[(2)]
{\sl (Analytic conditions)}
For any ${\bf a}\in\bbR\times[1,\infty]$, one has
$$
\|\Pi\|_{{\bf a};w_c}\defeq
\max_{\substack{(\alpha,p)\in A   \\   (\alpha,p)<{\bf a}}}
\sup_{0<t\le1}
 t^{-\alpha/\ell}
\Bigg\|
\sup_{\tau\in T_{(\alpha,p)}\setminus\{0\}}
\frac{\big|\mcQ_t(x,\Pi_x(\tau))\big|}{\|\tau\|_{(\alpha,p)}}
\Bigg\|_{L_x^p(w_c)} < \infty
$$
and 
$$
\|\Gamma\|_{{\bf a};w_c}\defeq\max_{\substack{(\alpha,p),(\beta,q)\in A   \\
(\beta,q)<(\alpha,p)<{\bf a}}} \sup_{y\in\bbR^d\setminus\{0\}} 
\Bigg\{\frac{w_c(y)}{\|y\|_\mfs^{\alpha-\beta}}\,
\Bigg\|
\sup_{\tau\in T_{(\alpha,p)}\setminus\{0\}}
\frac{\|\Gamma_{(x+y)x}(\tau)\|_{(\beta,q)}}{\|\tau\|_{(\alpha,p)}}
\Bigg\|_{L_x^{p:q}(w_c)} \Bigg\}< \infty.
$$
\end{itemize}
We write 
$$
\tri {\sf M}\tri_{{\bf a};w_c}\defeq\|\Pi\|_{{\bf a};w_c} + \|\Gamma\|_{{\bf a};w_c}.
$$
Furthermore, for any two models ${\sf M}_i = (\Pi_i,\Gamma_i)\in\textbf{\textsf{M}}(\scT)_{w_c}$ with $i\in\{1,2\}$, define the pseudo-metric 
$$
\tri {\sf M}_1 ; {\sf M}_2\tri_{{\bf a};w_c} \defeq\|\Pi_1 - \Pi_2\|_{{\bf a};w_c} + \|\Gamma_1 - \Gamma_2\|_{{\bf a};w_c}
$$ 
by replacing $\Pi$ and $\Gamma$ above with $\Pi_1 - \Pi_2$ and $\Gamma_1 - \Gamma_2$ respectively.
\end{defn}

\ssk

The topological space $\textbf{\textsf{M}}(\scT)_{w_c}$ is complete with respect to the pseudo-metrics $\tri {\sf M}_1 ; {\sf M}_2\tri_{{\bf a};w_c}$.

\ssk

\begin{defn} \label{DefnMD}
Let ${\sf M}=(\Pi,\Gamma)\in\textbf{\textsf{M}}(\scT)_{w_c}$. For every ${\bf c}=(\gamma,r)\in\bbR\times[1,\infty]$ and $b>0$, we define $D^{\bf c}(\Gamma)_{w_b}$ as the space of all functions $f:\bbR^d\to T_{<{\bf c}}\defeq \bigoplus_{{\bf a}\in A,\, {\bf a}<{\bf c}}T_{{\bf a}}$ such that
\begin{align*}
\lp f\rp_{{\bf c};w_b}&\defeq \max_{(\alpha,p)<{\bf c}}\big\|\|f(x)\|_{(\alpha,p)}\big\|_{L_x^{r:p}(w_b)} < \infty,   \\ 
\| f\|_{{\bf c};w_b}^\Gamma&\defeq\max_{(\alpha,p)<{\bf c}} \sup_{h\in\bbR^d\setminus\{0\}} w_b(y)\frac{\big\|\|\Delta_{x;h}^\Gamma f\|_{(\alpha,p)}\big\|_{L_x^{r:p}(w_b)}}{\|h\|_\mfs^{\gamma-\alpha}} < \infty,
\end{align*}
where 
$$
\Delta_{x;h}^\Gamma f\defeq f(x+h)-\Gamma_{(x+h)x}f(x).
$$
We call each element of $D^{\bf c}(\Gamma)_{w_b}$ a {\sl modelled distribution}.
In addition, we write 
$$
\tri f\tri_{{\bf c};w_b}^\Gamma\defeq\lp f\rp_{{\bf c};w_b} + \|f\|_{{\bf c};w_b}^\Gamma.
$$
Furthermore, for any models ${\sf M}_i=(\Pi_i,\Gamma_i)\in\textbf{\textsf{M}}(\scT)_{w_c}$ and functions $f_i\in D^{\bf c}(\Gamma_i)_{w_b}$ with $i\in\{1,2\}$, we define 
$$
\tri f_1 ; f_2\tri_{{\bf c};w_b}^{\Gamma_1;\Gamma_2} \defeq \lp f_1-f_2\rp_{{\bf c};w_b} + \|f_1;f_2\|_{{\bf c};w_b}^{\Gamma_1;\Gamma_2}
$$ 
by
\begin{align*}
\lp f_1-f_2\rp_{{\bf c};w_b} &\defeq \max_{(\alpha,p)<{\bf c}}\big\|\|f_1(x)-f_2(x)\|_{(\alpha,p)}\big\|_{L_x^{r:p}(w_b)},   \\
\|f_1;f_2\|_{{\bf c};w_b}^{\Gamma_1;\Gamma_2} &\defeq \max_{(\alpha,p)<{\bf c}}
\sup_{h\in\bbR^d\setminus\{0\}}w_b(y)
\frac{\big\|\|\Delta_{x;h}^{\Gamma_1} f_1-\Delta_{x;h}^{\Gamma_2} f_2\|_{(\alpha,p)}\big\|_{L_x^{r:p}(w_b)}}{\|h\|_\mfs^{\gamma-\alpha}}.
\end{align*}
For any subspace $V$ of $T$ which is stable under $G$, we denote by $D^{\bf c}(V;\Gamma)_{w_b}$ the space of all modelled distributions $f\in D^{\bf c}(\Gamma)_{w_b}$ which takes values in $V$
\end{defn}

\ssk

Recall from \eqref{EqDefnQBesovSpaces} the definition of the weighted Besov spaces $B^{\alpha,Q}_{pq}(w)$ associated to $Q$.

\ssk

\begin{defn} \label{def:besovreconst}
Let ${\sf M}=(\Pi,\Gamma)\in\textbf{\textsf{M}}(\scT)_{w_c}$ and ${\bf c}=(\gamma,r)\in\bbR\times[1,\infty]$. Then for any $f\in D^{\bf c}(\Gamma)_{w_b}$, any distribution $\Lambda\in B_{r,\infty}^{r_0,Q}(w_{b+c})$ satisfying
\begin{align*}
\lb\Lambda\rb_{{\bf c};w_{b+c}}^{\Pi,f} \defeq \sup_{0<t\le1}t^{-\gamma/\ell}\big\|\mcQ_t\big(x,\Lambda_x^{\Pi,f}\big)\big\|_{L_x^r(w_{b+c})} < \infty,
\qquad \big(\Lambda_x^{\Pi,f} \defeq \Lambda-\Pi_xf(x)\big)
\end{align*}
is called a \textbf{reconstruction} of $f$ for $\sf M$. Furthermore, for any models ${\sf M}_i=(\Pi_i,\Gamma_i)\in\textbf{\textsf{M}}(\scT)_{w_c}$ and functions $f_i\in D^{\bf c}(\Gamma_i)_{w_b}$ with $i\in\{1,2\}$, if there is a reconstruction $\Lambda_i$ for each $i$, we define
\begin{align*}
\lb\Lambda_1;\Lambda_2\rb_{{\bf c};w_{b+c}} \defeq\sup_{0<t\le1}t^{-\gamma/\ell}\big\|\mcQ_t\big(x, (\Lambda_1)_{x}^{\Pi_1,f_1} - (\Lambda_2)_{x}^{\Pi_2,f_2}\big)\big\|_{L_x^r(w_{b+c})}.
\end{align*}
\end{defn}

\ssk

The next statement is the regularity-integrability counterpart of the reconstruction theorem.

\ssk

\begin{thm} \label{thm:besovreconstruction}
Pick ${\sf M}=(\Pi,\Gamma)\in\textbf{\textsf{M}}(\scT)_{w_c}$ and ${\bf c}=(\gamma,r)\in\bbR\times[1,\infty]$ with $\gamma>0$. There exists a unique reconstruction ${\mcR}^{\sf M} f$ of $f\in D^{\bf c}(\Gamma)_{w_b}$ for ${\sf M}$. Moreover, it holds that
\begin{align*} 
\|{\mcR}^{\sf M} f\|_{B_{r,\infty}^{r_0,Q}(w_{b+c})} &\le C \|\Pi\|_{{\bf c};w_c}\tri f\tri_{{\bf c};w_b}^\Gamma,   \\
\lb{\mcR}^{\sf M} f\rb_{{\bf c};w_{b+c}}^{\Pi,f} &\le C \|\Pi\|_{{\bf c};w_c} \|f\|_{{\bf c};w_b}^\Gamma,
\end{align*}
where the constant $C$ depends only on $\scT$, $\bf c$, $b,c$, and $Q$.
Moreover there is an affine function $C_\lambda>0$ of $\lambda>0$ whose coefficients depend only on $\scT$, $\bf c$, $b,c$, and $Q$, such that
\begin{align*}
\big\| {\mcR}^{{\sf M}_1} f_1 - {\mcR}^{{\sf M}_2} f_2 \big\|_{B_{r,\infty}^{r_0,Q}(w_{b+c})} &\le C_\lambda \Big(\|\Pi_1-\Pi_2\|_{{\bf c};w_c} + \tri f_1;f_2\tri^{\Gamma_1,\Gamma_2}_{{\bf c};w_b}\Big),   \\
\big\lb {\mcR}^{{\sf M}_1} f_1 ; {\mcR}^{{\sf M}_2} f_2 \big\rb_{{\bf c};w_{b+c}} &\le C_\lambda\Big(\|\Pi_1 - \Pi_2\|_{{\bf c};w_c} + \| f_1 ; f_2\|^{\Gamma_1,\Gamma_2}_{{\bf c};w_b}\Big),
\end{align*}
for any models ${\sf M}_i=(\Pi_i,\Gamma_i)\in\textbf{\textsf{M}}(\scT)_{w_c}$ and functions $f_i\in D^{\bf c}(\Gamma_i)_{w_b}$ with $i\in\{1,2\}$ such that $\tri {\sf M}_i\tri_{{\bf c};w_c}\le \lambda$ and $\tri f_i\tri_{{\bf c};w_b}\le \lambda$.
\end{thm}

\ssk

We now lift the operator $\mcK$ into an operator that maps modelled distributions on modelled distributions.

\ssk

\begin{defn} \label{defTbarT}
Let $\scT=(A,T,G)$ be a regularity-integrability structure. In addition, let $\overline\scT=(\overline{A},\overline{T},\overline{G})$ be a regularity-integrability structure satisfying the following properties.
\begin{itemize}
	\item[(1)] $\bbN[\mfs]\times\{\infty\}\subset \overline{A}$, where $\bbN[\mfs] \defeq \big\{|{k}|_{\mfs}\, ;\, {k}\in\bbN^d\big\}$.
	
	\item[(2)] For each $\alpha\in\bbN[\mfs]$, the vector space $\overline{T}_{(\alpha,\infty)}$ contains all $X^{k}$ with $|{k}|_{\mfs}=\alpha$.

	\item[(3)] The linear space $\spa\{X^{k}\}_{{k}\in\bbN^d}$ is closed under the action of the group $\overline{G}$.
\end{itemize}
For ${\sf M}\in\textbf{\textsf{M}}(\scT)_{w_c}$ we define the linear map $\mcJ^{\sf M}(x):T\to\spa\{X^{k}\}\subset\overline{T}$ by setting
$$
\mcJ^{\sf M}(x)\tau \defeq \sum_{|{k}|_{\mfs}<\alpha+\beta_0}\frac{X^{k}}{{k}!} \, \partial^{k} \mcK(x,\Pi_x\tau)
$$
for any $(\alpha,p)\in A$ and $\tau\in T_{(\alpha,p)}$. A continuous linear map $\mcI:T\to\overline{T}$ is called an \textbf{abstract integration map of order $\beta_0\in(0,\ell-\ell_1)$} if 
$$
\mcI:T_{(\alpha,p)}\to\overline{T}_{(\alpha+\beta_0,p)}
$$
for any $(\alpha,p)\in A$. The models ${\sf M}=(\Pi,\Gamma)\in\textbf{\textsf{M}}(\scT)_{w_c}$ and $\overline{\sf M}=(\overline\Pi ,\overline{\Gamma})\in\textbf{\textsf{M}}(\overline{\scT})_{w_c}$ are said to be \textbf{compatible for $\mcI$} if they satisfy the following properties.
\begin{enumerate} 
\renewcommand{\theenumi}{(\roman{enumi})}
\renewcommand{\labelenumi}{(\roman{enumi})}
\item\label{Iadmissible2}
For any ${k}\in\bbN^d$ one has 
$$
(\overline\Pi _xX^{k})(\cdot) = (\cdot-x)^{k},\qquad
\overline\Gamma_{yx} X^{k} = \sum_{{l}\le {k}} \binom{{k}}{{l}} (y-x)^{l} X^{{k}-{l}}.
$$
\item\label{Iadmissible3}
For any $x,y\in\bbR^d$ one has
$$
\overline{\Gamma}_{yx}\circ\big(\mcI+\mcJ^{\sf M}(x)\big) = \big(\mcI + \mcJ^{{\sf M}}(y)\big)\circ\Gamma_{yx}.
$$
\end{enumerate}
\end{defn}

\ssk

For any model ${\sf M} = (\Pi,\Gamma)\in\textbf{\textsf{M}}(\scT)_{w_c}$, $f\in D^{\bf c}(\Gamma)_{w_b}$ with ${\bf c}=(\gamma,r)\in\bbR\times[1,\infty]$, and its reconstruction $\Lambda$, we set
$$
\mcN^{\sf M}(x;f,\Lambda) \defeq \sum_{|{k}|_{\mfs}<\gamma+\beta_0} \frac{X^{k}}{{k}!} \, \partial^{k} \mcK\big(x,\Lambda_x^{\Pi,f}\big).
$$
Finally, we define
\begin{equation} \label{def:MS}
\mcK^{\sf M} f(x)\defeq\mcI(f(x)) + \mcJ^{\sf M}(x) (f(x)) + \mcN^{\sf M}(x;f,\Lambda).
\end{equation}

\ssk

\begin{thm} \label{thm:MS}
Let $\beta_0\in(0,\ell-\ell_1)$ and ${\bf c}=(\gamma,r)\in\bbR\times[1,\infty]$. Then for any ${\sf M} = (\Pi,\Gamma)\in\textbf{\textsf{M}}(\scT)_{w_c}$, $f\in D^{\bf c}(\Gamma)_{w_b}$, and any reconstruction $\Lambda$ of $f$ for $\sf M$, the function $\mcK^{\sf M} f$ belongs to $D^{(\gamma+\beta_0,r)}(\overline\Gamma)_{w_{2c+b}}$, and
\begin{align*}
\lp \mcK^{\sf M} f\rp_{(\gamma+\beta_0,r);w_{2c+b}} 
&\le C\Big\{(1+\|\Pi\|_{{\bf c};w_c}) (1+\|\Gamma\|_{{\bf c};w_c}) \lp f\rp_{{\bf c};w_b}+\lb\Lambda\rb_{{\bf c};w_{c+b}}^{\Pi,f}\Big\},   \\
\|\mcK^{\sf M} f\|_{(\gamma+\beta_0,r);w_{2c+b}}^{\overline\Gamma} 
&\le C\Big\{(1+\|\Pi\|_{{\bf c};w_c}) (1+\|\Gamma\|_{{\bf c};w_c}) \| f\|_{{\bf c};w_b}^{\Gamma}+\lb\Lambda\rb_{{\bf c};w_{c+b}}^{\Pi,f}\Big\},
\end{align*}
where the constant $C$ depends only on $\scT$, $\bar{\scT}$, $\beta_0$, $\bf c$, $b,c$, and $Q$.
Moreover there is a quadratic function $C_\lambda>0$ of $\lambda>0$ whose coefficients depend only on $\scT$, $\bar{\scT}$, $\beta_0$, $\bf c$, $b,c$, and $Q$, such that
\begin{align*}
\tri\mcK^{{\sf M}_1} f_1;\mcK^{{\sf M}_2} f_2\tri_{(\gamma+\beta_0,r);w_{2c+b}}
\le C_\lambda \Big( \tri {\sf M}_1;{\sf M}_2\tri_{{\bf c};w_c} + \tri f_1 ; f_2\tri_{{\bf c};w_b} + \|\Lambda_1 ; \Lambda_2\|_{{\bf c};w_{c+b}} \Big),
\end{align*}
for any models ${\sf M}_i = (\Pi_i,\Gamma_i)\in\textbf{\textsf{M}}(\scT)_{w_c}$ and $\overline{{\sf M}}_i = (\overline\Pi_i , \overline\Gamma_i) \in \textbf{\textsf{M}}(\overline\scT)_{w_c}$ such that ${\sf M}_i$ and $\overline{{\sf M}}_i$ are compatible, and functions $f_i\in D^{\bf a}(\Gamma_i)_{w_b}$ with $i\in\{1,2\}$ such that $\tri {\sf M}_i\tri_{{\bf c};w_c}\le \lambda$ and $\tri f_i\tri_{{\bf c};w_b}\le \lambda$.
\end{thm}

\bigskip

\noindent \textcolor{gray}{$\bullet$} {\sf I. Bailleul} -- Univ Brest, CNRS UMR 6205, Laboratoire de Math\'ematiques de Bretagne Atlantique, France. {\it E-mail}: ismael.bailleul@univ-brest.fr   \\
{\small Partial support from the ANR-11-LABX-0020-01 Labex CHL and the ANR project Smooth ANR-22-CE40-0017 is acknowledged.}

\medskip

\noindent \textcolor{gray}{$\bullet$} {\sf M. Hoshino} -- Department of Mathematics, Institute of Science Tokyo, Japan. \\
{\it E-mail}: hoshino@math.sci.isct.ac.jp.   \\
{\small Support from JSPS KAKENHI Grant Number 23K12987 is acknowledged.}

\end{document}